\documentclass{amsproc}

\usepackage{amsmath}
\usepackage{amsxtra}

\usepackage{amssymb}
\usepackage[mathscr]{eucal}     
\DeclareSymbolFont{pssymbols}     {OMS}{ztmcm}{m}{n}
\DeclareSymbolFontAlphabet{\mathpsscr}   {pssymbols}

\usepackage{amsthm}
\theoremstyle{plain}
\newtheorem{thm}{Theorem}[section]
\newtheorem{conj}[thm]{Conjecture}
\newtheorem{prop}[thm]{Proposition}
\newtheorem{lem}[thm]{Lemma}

\theoremstyle{remark}
\newtheorem*{rem}{Remark}
\newtheorem*{rems}{Remarks}
\newtheorem*{notation}{Notation}

\newcommand{\itemref}[1]{\textup{(\ref{#1})}}

\usepackage[arrow,matrix,tips,frame,curve,arc,web,ps,dvips,color,line,poly]{xy}
\SelectTips{cm}{}
\UseTips
\CompileMatrices
\newdir^{ (}{{}*!/-5pt/@^{(}}
\newcommand{\xyhookrightarrow}{\ar@<.5ex>@{^{ (}->}}
\newxyColor{lightgrey}{.5}{gray}{}

\newcommand{\CC}{{\mathbb C}}
\newcommand{\RR}{{\mathbb R}}
\newcommand{\QQ}{{\mathbb Q}}

\newcommand{\ZZ}{{\mathbb Z}}
\newcommand{\EE}{{\mathbb E}}
\newcommand{\HH}{{\mathbb H}}
\newcommand{\VV}{{\mathbb V}}
\newcommand{\PP}{{\mathbb P}}

\DeclareMathOperator{\codim}{codim}

\DeclareMathOperator{\rank}{rank}
\DeclareMathOperator{\Ker}{Ker}
\renewcommand{\Im}{\operatorname{Im}}
\DeclareMathOperator{\QQrank}{\QQ-rank}
\DeclareMathOperator{\RRrank}{\RR-rank}
\DeclareMathOperator{\CCrank}{\CC-rank}
\DeclareMathOperator{\mS}{SS}
\DeclareMathOperator{\emS}{\mS_{ess}}
\DeclareMathOperator{\Res}{Res}

\renewcommand{\l}{\ell}
\DeclareMathOperator{\Cl}{cl}
\newcommand{\cl}[1]{\Cl(#1)}
\DeclareMathOperator{\diag}{diag}
\newcommand{\geo}{\mathbin{\mathbf o}}
\DeclareMathOperator{\Herm}{Herm}
\DeclareMathOperator{\dom}{Dom}
\newcommand{\dbar}{\bar d}
\DeclareMathOperator{\Fix}{Fix}

\newcommand{\lsb}[1]{{}_{\scriptscriptstyle #1}}
\newcommand{\lsp}[1]{{}^{#1}\!}
\newcommand{\tildearrow}{\xrightarrow{\sim}}
\newcommand{\longtildearrow}{\mathrel{\overset{\textstyle\mathchar"0366}%
{\smash\longrightarrow}}}

\newcommand{\G}{\Gamma}
\DeclareMathOperator{\GL}{GL}
\DeclareMathOperator{\SL}{SL}
\DeclareMathOperator{\SO}{SO}
\renewcommand{\O}{\operatorname{O}}
\newcommand{\Sympl}{\operatorname{Sp}}
\newcommand{\GSympl}{\operatorname{GSp}}

\DeclareMathOperator{\U}{U}
\newcommand{\Pl}{\mathscr P}

\DeclareMathOperator{\Hom}{Hom}

\newcommand{\Abar}{{\bar{A}}}
\newcommand{\Xbar}{\overline{X}}
\newcommand{\Xhat}{\widehat{X}}
\newcommand{\back}{\backslash}
\newcommand{\Dbar}{\overline{D}}
\newcommand{\ebar}{\overline{e}}
\newcommand{\Dstar}{D^*}
\newcommand{\Xstar}{X^*}
\newcommand{\nil}{\mathscr N'}

\renewcommand{\i}{i}		
\newcommand{\ihat}{{\hat{\imath}}}
\renewcommand{\j}{j}		

\newcommand{\IpC}{{\mathcal I_p\mathcal C}}
\newcommand{\Sheaf}{\mathcal S}
\newcommand{\Hsheaf}{H}

\newcommand{\WnC}{\mathcal W^\eta \mathcal C}
\renewcommand{\L}{\mathscr L}
\newcommand{\M}{\mathcal M}
\newcommand{\A}{\mathcal A}
\renewcommand{\sp}{_{\,\text{sp}}}
\newcommand{\Asp}{\A\sp}

\DeclareMathOperator{\Lie}{Lie}
\newcommand{\sa}{{\mathfrak a}}
\newcommand{\h}{{\mathfrak h}}

\newcommand{\n}{{\mathfrak n}}
\newcommand{\levi}{{\mathfrak l}}

\newcommand{\al}{\alpha}
\renewcommand{\b}{\beta}
\newcommand{\D}{\Delta}
\newcommand{\Dhat}{\widehat\D}
\newcommand{\g}{\gamma}
\renewcommand{\u}{\mu}
\renewcommand{\v}{\nu}
\renewcommand{\r}{\rho}
\renewcommand{\t}{\tau}
\newcommand{\kap}{\kappa}
\renewcommand{\o}{\omega}
\newcommand{\z}{\zeta}

\DeclareMathOperator{\IrrRep}{\mathfrak I\mathfrak r\mathfrak r}
\DeclareMathOperator{\Mor}{Mor}

\let\Der\Derived

\begin{document}


\author{Leslie Saper}
\address{Department of Mathematics\\ Duke University\\ Box 90320\\ Durham,
NC 27708\\U.S.A.}
\email{saper@math.duke.edu}
\urladdr{http://www.math.duke.edu/faculty/saper}
\title[Cohomology of Locally Symmetric Spaces]{On the Cohomology of Locally
  Symmetric Spaces and of their Compactifications}
\thanks{An expanded version of two talks  at the Harvard-M.I.T. conference
  ``Current Developments in Mathematics, 2002'', November 15--17.}
\thanks{This material is based upon work supported in part by the National
  Science Foundation under grant DMS-9870162.  The manuscript was typeset
  using AMS-\LaTeX; the \Xy-pic\ package was used for figures and
  diagrams.}
\begin{abstract}
This expository article gives an introduction to the (generalized)
conjecture of Rapoport and Goresky-MacPherson which identifies the
intersection cohomology of a real equal-rank Satake compactification of a
locally symmetric space with that of the reductive Borel-Serre
compactification.  We motivate the conjecture with examples and then give
an introduction to the various topics that are involved: intersection
cohomology, the derived category, and compactifications of a locally
symmetric space, particularly those above.  We then give an overview of the
theory of $\L$-modules and micro-support which was developed to solve the
conjecture but has other important applications as well.  We end with
sketches of the proofs of three main theorems on $\L$-modules that lead to
the resolution of the conjecture.  The text is enriched with many examples,
illustrations, and references to the literature.
\end{abstract}
\maketitle


\setcounter{tocdepth}{1}
\tableofcontents


\section{An example}
\label{sectExample}
We begin with the simplest example of a locally symmetric space.  This is
the moduli of complex elliptic curves $\CC/\Lambda$, where $\Lambda$ is a
lattice in $\CC$.  The moduli space can be expressed as
\begin{align*}
    &\frac{\text{\{elliptic curves\}}}{\text{$\CC$-isomorphism}}\\
    &\qquad\longleftrightarrow \{\text{lattices $\Lambda$ in 
    $\CC$}\}/{\SO(2)\RR^*_+}\\
    &\qquad\longleftrightarrow \SL_2(\ZZ) \back \GL_2(\RR)^0 / \SO(2)
\RR^*_+\\
    &\qquad\longleftrightarrow \SL_2(\ZZ) \back H
\end{align*}
where $H = \{z\in\CC\mid \Im z > 0\}$ is the upper half-plane and $\RR^*_+$
denotes the group of positive real scalar matrices.  Here we let
$\bigl(\begin{smallmatrix} a & b \\ c & d
\end{smallmatrix}\bigr)\in \GL_2(\RR)^0$ parametrize the lattice
$\Lambda=\ZZ(a + ib) + \ZZ(c+id)$; the quotient by $\SL_2(\ZZ)$ removes the
dependence on the choice of an oriented basis. For future purposes, we will
actually view $H$ as $\GL_2(\RR)/\O(2)\RR^*_+$.

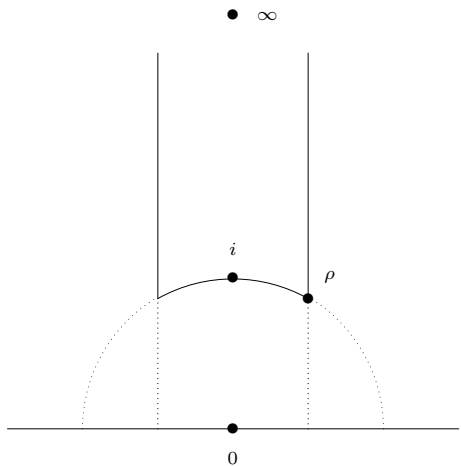
\begin{figure}[t]
\begin{equation*}
\begin{xy}
0;<2cm,0cm>:
(-1.5,0);(1.5,0)**\crv~l{},
(-.5,2.5);(-.5,.866025404)**\crv~l{},
(.5,.866025404);(.5,2.5)**\crv~l{},
(-.5,.866025404);(0,1.154700538),**{},(.5,.866025404)="rho",{\ellipse_{}},
(-1,0);(-1,1),**{},(-.5,.866025404),{\ellipse_{.}},
(1,0);(1,1),**{},(.5,.866025404),{\ellipse^{.}},
(-.5,.866025404);(-.5,0)**\crv{~*=<4mm>{}~**@{.}},
(.5,.866025404);(.5,0)**\crv{~*=<4mm>{}~**@{.}},
(0,1)*{\bullet}*+++!DC{\scriptstyle i},
(.5,.866025404)*{\bullet}*++!LD{\scriptstyle\rho},
(0,0)*{\bullet}*+++!UC{\scriptstyle 0},
(0,2.75)*{\bullet}*+++!CL{\scriptstyle \infty}
\end{xy}
\end{equation*}
\caption{The fundamental domain for the action of $\SL_2(\ZZ)$ on the
  upper half-plane}
\label{figSLTwoZ}
\end{figure}

A fundamental domain for the action of $\SL_2(\ZZ)$ is pictured in Figure
~\ref{figSLTwoZ}; the action of $\SL_2(\ZZ)$ further identifies $z
\leftrightarrow -\bar z$ for points in the boundary of the domain.  Note
the ``cusp'' at infinity which is not properly part of $X =\SL_2(\ZZ) \back
H$ but may be added to yield a compactification $\Xstar$.
One may also consider $X=\G\back H$ where $\G$ is a finite-index subgroup
of $\SL_2(\ZZ)$.  Then $\Xstar$ is formed by adding possibly several cusp
points.  Any representation $(\sigma,E)$ of $\GL_2(\RR)$ defines a locally
constant sheaf $\EE$ on $X$, namely $\G\back (H\times E) = H\times E /
{\sim}$, where $(x,v)\sim (\g x,\sigma(\g)v)$ for all $\g\in \G$.  The
classical Eichler-Shimura theorem \cite{refnShimura} states

\begin{thm}
\label{thmEichlerShimura}
Let $E_k$ be the representation of $\GL_2(\RR)$ with highest weight $k$.
Then $H_P^1(\Xstar;\EE_k) \cong \mathscr S_{k+2}(\G) \oplus
\overline{\mathscr S_{k+2}(\G)}$.
\end{thm}

The group appearing on the left-hand side is the \emph{parabolic
cohomology}:
\begin{equation*}
   H_P^1(\Xstar;\EE_k) \equiv \Ker \Bigl(H^1(X;\EE_k) \to
   \bigoplus_{x\in\text{cusps}} H^1(U_x;\EE_k)\Bigr),
\end{equation*}
where $U_x$ is a punctured neighborhood of the cusp point $x$.  It may also
be defined by a modification of the usual complex computing the group
cohomology $H(\G;E_k)$.  One should think of parabolic cohomology as a
modification of ordinary cohomology in which the local cohomology near a
cusp in degree $1$ has been killed or truncated.

The space $\mathscr S_{k+2}(\G)$ on the right-hand side is the space of
\emph{modular cusp forms of weight $k+2$ for $\G$}:
\begin{equation*}
f \in \mathscr S_{k+2}(\G) \Longleftrightarrow \left\{ \begin{array}{ll}
& f\colon H \to \CC \text{ holomorphic},\\
& f\left(\frac{az+b}{cz+d}\right) = (cz+d)^{k+2} f(z) \text{ for }
        \Bigl(\begin{matrix} a & b \\ c & d \end{matrix}\Bigr)\in \G,\\
&\text{$f$ vanishes at the cusps}.
\end{array}\right.
\end{equation*}

Modular forms (and their generalization, automorphic forms) play a central
role in number theory.  For example, note that the space $\Xstar$ is
actually an algebraic variety defined over a number field.  An important
application of Theorem ~\ref{thmEichlerShimura} is that, for $\G$ a
congruence subgroup, it allows one to relate the Hasse-Weil zeta function
of $\Xstar$ (which encodes the number of points of $\Xstar$ defined over
all finite fields) to the $L$-functions associated to modular forms.  We
will not discuss this further here; for details see Shimura's book
\cite{refnShimura}.

\section{$L^2$-cohomology}
We now interpolate a group within the Eichler-Shimura isomorphism:
\begin{equation}
\label{eqnEichlerShimuraLtwo}
H_P^1(\Xstar;E_k) \cong H_{(2)}^1(X;\EE_k) \cong \mathscr S_{k+2}(\G)
\oplus \overline{\mathscr S_{k+2}(\G)}
\end{equation}
The new group $H_{(2)}(X;\EE_k)$ is the \emph{$L^2$-cohomology}.

The $L^2$-cohomology is a variant of de Rham cohomology formed by imposing
an $L^2$-growth condition on the differential forms.  We will briefly
review the theory in a general context; good references are the papers of
Cheeger \cite{refnCheeger}, Cheeger, Goresky, and MacPherson
\cite{refnCheegerGoreskyMacPherson}, and Zucker \cite{refnZuckerWarped}.
Suppose that $M$ is an oriented Riemannian manifold (a $V$-manifold
\cite{refnSatakeVManifold} or orbifold is sufficient) and that $\EE$ is a
metrized locally constant sheaf.  In other words, $\EE$ is the sheaf of
sections of a flat vector bundle on $M$ equipped with a fiber-wise
Hermitian metric (which is not necessarily flat).  Let $d$ denote the
exterior derivative on $A(M;\EE)$, the smooth $\EE$-valued differential
forms on $M$.  For every $x\in M$, the norm on the tangent space $T_xM$ and
the norm on the fiber $\EE_x$ induce a norm on $\bigwedge T^*_xM\otimes
\EE_x$.  Thus for $\phi\in A(M;\EE)$, we have a function $x\mapsto
|\phi(x)|$.  On the other hand, the Riemannian metric also induces a volume
form $dV$.  We define the $L^2$-norm of $\phi$ by
\begin{equation}
\|\phi\|^2 \equiv \int_M |\phi|^2 dV
\end{equation}
and set
\begin{equation}
A_{(2)}(M;\EE) \equiv \{\, \phi\in A(M;\EE) \mid \|\phi\|<\infty
\text{ and } \|d\phi\|<\infty\,\}.
\end{equation}
The $L^2$-cohomology $H_{(2)}(M;\EE)$ is defined to be the cohomology of
the complex $(A_{(2)}(M;\EE),d)$; it has a topology (which may not be
Hausdorff) given by the semi-norm
\begin{equation*}
\bigl\|\,[\phi]\,\bigr\| \equiv \inf_{\eta\in A_{(2)}(M;\EE)} \|\phi+d\eta\|
\qquad \text{for $[\phi]\in H_{(2)}(M;\EE)$.}
\end{equation*}
The $L^2$-cohomology is a \emph{quasi-isometry} invariant of the metrics on
$M$ and $\EE$: two metrics $h$ and $h'$ are said to be quasi-isometric if
$c^{-1} h \le h'\le ch$ uniformly for some $c>0$.

Sometimes it is useful to consider forms which are not necessarily smooth,
but only measurable.  The measurable forms $\phi$ which are $L^2$ and which
have exterior derivatives in the distribution sense represented by $L^2$
forms yield a complex which we denote $(\dom \dbar,\dbar)$; this is indeed
the closure in the sense of unbounded operators on Hilbert space of the $d$
operator acting on the domain $A_{(2)}(M,\EE)$.  A smoothing argument
\cite{refnCheeger}, \cite{refndeRham} shows that this extended complex
computes the same cohomology,
\begin{equation*}
H_{(2)}(M;\EE) = H(\dom \dbar) = \Ker \overline d\left/ \Im \dbar\right..
\end{equation*}
It follows that
\begin{align*}
H_{(2)}(M;\EE) &= \Bigl({\Ker \dbar}\left/{\overline{\Im
  \dbar}}\right.\Bigr) \oplus \Bigl( \left.{\overline{\Im
  \dbar}}\right/{\Im \dbar}\Bigr) \\ &\cong \Bigl(\Ker \dbar\cap \Ker
  \dbar^*\Bigr) \oplus \Bigl( \left.{\overline{\Im \dbar}}\right/{\Im
  \dbar}\Bigr).
\end{align*}
This is Hodge theory: the first term is a certain space of $L^2$ harmonic
forms (in fact all $L^2$ harmonic forms if the Riemannian metric is
complete) and the second term is either $0$ or infinite dimensional.  Note
that $H_{(2)}(M;\EE)$ is Hausdorff if and only if the second term is $0$.

In the example of \S\ref{sectExample}, we give $H$ the hyperbolic metric
$(dx^2 + dy^2)/y^2$ under which $H$ becomes a symmetric space.  Since the
metric is invariant under the action of $\SL_2(\ZZ)$, it descends to a
metric on $X=\G\back H$, a locally symmetric space.  Furthermore, any
finite dimensional representation $(\sigma,E)$ of $\GL_2(\RR)$ admits an
\emph{admissible inner product}, that is, one for which $\sigma(g)^* =
\sigma(\theta g)^{-1}$, where $\theta$ denotes the Cartan involution
$g\mapsto \lsp{t}g^{-1}$.  Such an inner product induces a Hermitian metric
on $\EE = \G\back (H\times E)$, namely
\begin{equation}
\| (g\SO(2),v) \| = \|\sigma(g^{-1})v\| \qquad \text{for $g\in
  \SL_2(\RR)$.}
\end{equation}
A variant of this definition applies to $g\in\GL_2(\RR)$.  These metrics on
$X$ and $\EE$ are unique up to quasi-isometry (in fact up to scalar
multiples if $E$ is irreducible) and thus the $L^2$-cohomology
$H_{(2)}(X;\EE_k)$ is well-defined.

The first isomorphism of \eqref{eqnEichlerShimuraLtwo} is heuristically
true since the metric near a cusp, say the cusp at $\infty$, may be
expressed as
\begin{equation}
\label{eqnCuspMetric}
dr^2 + e^{-2r}d\theta^2, \qquad\text{where $r=\log y\in[b,\infty)$ and
    $\theta = x$ modulo $\ZZ$,}
\end{equation}
and hence $d\theta$ is not $L^2$.  The second isomorphism of
\eqref{eqnEichlerShimuraLtwo} follows from Hodge theory and the Hodge
decomposition: the coefficients of harmonic forms representing type $(1,0)$
cohomology classes can easily be seen to be modular forms of the
appropriate weight; the $L^2$ condition implies the form is cuspidal.

Other moduli problems yield locally symmetric spaces and in view of the
application mentioned in \S\ref{sectExample} and others, one would like a
generalization of \eqref{eqnEichlerShimuraLtwo}.  In this paper, we will
focus on the first isomorphism of \eqref{eqnEichlerShimuraLtwo}.

\section{Zucker's conjecture}
\label{sectZuckersConjecture}
In general a \emph{locally symmetric space} has the form
\begin{equation}
\label{eqnLocallySymmetricSpace}
X =  \G \back G / K A_G = \G \back D
\end{equation}
where
\newcommand{\defnwd}{.6\textwidth}
\begin{align*}
G &\equiv \parbox[t]{\defnwd}{the real points of a connected
  reductive algebraic group defined over $\QQ$,} \\
K &\equiv \parbox[t]{\defnwd}{a maximal compact subgroup of $G$,} \\
A_G &\equiv \parbox[t]{\defnwd}{the identity component  $(\RR^*_+)^s$
  of a maximal $\QQ$-split torus in the center of $G$,} \\
D &\equiv \parbox[t]{\defnwd}{$G / K A_G$, the corresponding
  symmetric space,} \\
\G &\equiv \parbox[t]{\defnwd}{an arithmetic subgroup of $G$.}
\end{align*}
We always give $D$ a $G$-invariant Riemannian metric; this induces a
Riemannian metric on $X$.  Although the locally symmetric space $X$ may
have finite quotient singularities, for simplicity we will treat it and
related spaces as if they were smooth; the correct treatment involves the
language of $V$-manifolds or orbifolds.  As in the case of $\GL_2(\RR)$,
any representation $E$ of $G$ induces a metrized locally constant sheaf
$\EE$ on $X$.

\begin{notation}
In this paper we will indicate Lie groups, often defined as the real points
of an algebraic group, by uppercase roman letters; their Lie algebras will
be denoted by the corresponding lowercase fraktur letter.  In order to
lighten the exposition and the notation, we will take some liberties in
terminology in this paper.  Specifically if $H$ is the Lie group of the
real points (or the identity component of the real points) of an algebraic
group defined over $\QQ$, we may refer to $H$ when properly we should refer
to the underlying algebraic group.  For example we might speak of a
rational character of $H$ defined over $\QQ$.  Or we will speak of a
parabolic $\QQ$-subgroup $P$ of $H$, when we properly should be referring
to the real points of a parabolic $\QQ$-subgroup of the algebraic group
underlying $H$.
\end{notation}

An important special case is a \emph{Hermitian locally symmetric space},
where $D$ has a $G$-invariant complex structure.  The generalization of the
first isomorphism in \eqref{eqnEichlerShimuraLtwo} is

\begin{thm}[Zucker's Conjecture]
\label{thmZuckersConjecture}
Let $X$ be a Hermitian locally symmetric space.  There is a natural
isomorphism $I_pH(\Xstar;\EE)\cong H_{(2)}(X;\EE)$, where $I_pH$ denotes
middle-perversity intersection cohomology.
\end{thm}

This result was first conjectured by Zucker in \cite{refnZuckerWarped}
where some simple cases were proven.  Other special cases of this result
were proven by Borel \cite{refnBorelQrankOne}, Borel and Casselman
\cite{refnBorelCasselmanQrankTwo}, and Zucker \cite{refnZuckerLtwoIHTwo}.
The theorem in general was proved by Looijenga \cite{refnLooijenga} and
independently by the author and Stern \cite{refnSaperSternTwo}.

There are two new ingredients in this theorem we need to explain: the
middle-perversity intersection cohomology $I_pH(\Xstar;\EE)$, which plays
the role of parabolic cohomology, and the Baily-Borel-Satake
compactification $\Xstar$, which generalizes adjoining cusp points to
$\G\back H$.  Along the way we will also introduce the Borel-Serre
compactification $\Xbar$ and the reductive Borel-Serre compactification
$\Xhat$ which will be important for the conjecture of Rapoport and
Goresky-MacPherson.

\section{Intersection cohomology}
\label{sectIntersectionCohomology}
Intersection cohomology, developed by Goresky and MacPherson
\cite{refnGoreskyMacPhersonIHOne}, \cite{refnGoreskyMacPhersonIHTwo}, is a
cohomology theory for singular spaces $Z$ which shares many of the
properties of ordinary cohomology theory applied to smooth spaces
(e.g. Poincar\'e duality).  It depends on a locally constant sheaf $\EE$ on
the smooth locus and a \emph{perversity}:
\begin{equation*}
p\colon \{2,\dots,d\} \to \ZZ, \quad p(2)= 0, \quad p(k)\le p(k+1)\le p(k)
+ 1.
\end{equation*}
We will mainly be interested in $p$ being one of the \emph{middle
perversities}:
\begin{equation}
\label{eqnMiddlePerversities}
n(k)= \left\lfloor\frac{(k-1)}2\right\rfloor \qquad \text{or}
\qquad m(k)= \left\lfloor\frac{(k-2)}2\right\rfloor;
\end{equation}
this is what is required in Zucker's conjecture and later in the conjecture
of Rapoport and Goresky-MacPherson.  Besides the original papers, a useful
reference is Borel's book \cite{refnBorelIntersectionCohomology}.

\subsection{Simplicial intersection cohomology}
\label{ssectSimplicialIC}
The type of singular space for which $I_pH(Z;\EE)$ is usually defined is a
pseudomanifold (see \cite{refnHabeggerSaper}, \cite{refnKing}, and
\cite{refnQuinn} for more general contexts).  A \emph{$d$-dimensional
stratified pseudomanifold} $Z$ is a topological space with a filtration by
closed subsets
\begin{equation*}
Z = Z_d \supset Z_{d-1} =Z_{d-2} \supset Z_{d-3} \supset \dotsb \supset Z_0
\supset Z_{-1} = \emptyset
\end{equation*}
such that
\begin{enumerate}
\item  For all $k\ge0$, the subspace $S^k = Z_{d-k} \setminus Z_{d-k-1}$
is a topological $(d-k)$-manifold, called the \emph{codimension-$k$
stratum}.
\item The stratum $S^0$ is dense.
\item Any point $x\in S^k$ has an open neighborhood $U$ in $Z$ which is
  homeomorphic (as a stratified pseudomanifold) to $B_{d-k} \times
  c(L_{k-1})$, where $B_{d-k}$ is a ball neighborhood of $x$ in $S^k$ and
  $c(L_{k-1})$ is a cone on a $(k-1)$-dimensional stratified pseudomanifold
  called the \emph{link}.
\end{enumerate}
The coefficient system $\EE$ is a locally constant sheaf defined on $S^0$,
the nonsingular stratum of $Z$.  Examples of pseudomanifolds include
complex analytic varieties as well as the reductive Borel-Serre
compactification and the Satake compactifications of locally symmetric
spaces to be discussed later.

The first definition of intersection cohomology
\cite{refnGoreskyMacPhersonIHOne}, however, used a piecewise-linear
structure on $Z$.  Briefly, assume that $Z$ is a piecewise-linear
stratified pseudomanifold (that is, all spaces and maps in the above
definition are taken in the piecewise-linear category).  A locally finite
\emph{geometric chain} $\xi$ is a locally finite simplicial chain for some
triangulation of $Z$, however we identify two chains if they agree after
passing to a common refinement of the triangulations; let $|\xi|$ denote
the support of the chain. Let $\mathcal O$ denote the orientation sheaf of
$S^0$ and let $I_pC^j(Z;\EE)$ denote the locally finite geometric
$(d-j)$-chains $\xi$ with values in $\EE\otimes \mathcal O$ that satisfy
the \emph{allowability} conditions
\begin{align*}
\dim(|\xi|\cap Z_{d-k}) & \le d-j -k +p(k),\\
\dim(|\partial\xi|\cap Z_{d-k}) & \le d-j -1 -k +p(k),
\end{align*}
for all $k\ge2$.  Even though $\EE$ is not defined on all of $Z$, the
equations above for $k=2$ ensure that the interiors of the
$(d-j)$-simplices and the $(d-j-1)$-simplices in $\xi$ lie in $S^0$, so it
makes sense to say that $\xi$ has ``values in $\EE \otimes \mathcal O$''.
The role of the perversity $p(k)$ in these equations is to bound the amount
of intersection allowed beyond that which would occur in the case of
transversal intersection.

Under the usual boundary operator, $I_pC(Z;\EE)$ becomes a complex whose
cohomology is the intersection cohomology $I_pH(Z;\EE)$.  The key thing to
keep in mind about intersection cohomology is its local characterization.
If $U = B_{d-k} \times c(L_{k-1})$ is a local neighborhood of $x\in S^k$ as
above, then
\begin{equation}
\label{eqnLocalCalculation}
I_pH^j(B_{d-k} \times c(L_{k-1});\EE) \cong \begin{cases}
	I_pH^j(L_{k-1};\EE) & \text{for $j \le p(k)$,} \\
	0 		& \text{for $j > p(k)$.}
			\end{cases}
\end{equation}
Indeed this calculation and the facts that the Mayer-Vietoris sequence
holds for intersection cohomology and that intersection cohomology agrees
with ordinary cohomology on a manifold, suffice to inductively determine
$I_pH(Z;\EE)$.  Since $m(k)= n(k)$ for $k$ even, it follows that both $p=m$
and $p=n$ yield the same theory for $Z$ with only even codimension strata
(for example, a complex analytic variety); certain spaces with possibly odd
codimension strata will become important later.  In our example of a
quotient of the upper half-plane, a cusp $x$ has codimension $2$ and so the
local intersection cohomology at $x$ vanishes in degrees $j>p(2)=0$ for
$p=m$ or $n$; it is easy to see then that $I_pH^1(\Xstar;\EE_k) \cong
H_P^1(\Xstar;\EE_k)$.

The constraints on a ordinary perversity and the condition in the
definition of a stratified pseudomanifold that $S^1=\emptyset$ ensure that
$I_pH(Z;\EE)$ is independent of the stratification; it is not apparent
however with this combinatorial definition that intersection cohomology is
independent of the piecewise-linear structure.

\subsection{The derived category}
It has proven useful to take a sheaf-theoretic approach to intersection
cohomology \cite{refnGoreskyMacPhersonIHTwo}; in particular, this allows an
easy proof of topological invariance.  In this approach, intersection
cohomology is the hypercohomology of a certain object $\IpC(Z;\EE)$ in the
constructible derived category of sheaves.  We begin by briefly describing
this category; a good reference is the book of Kashiwara and Schapira
\cite{refnKashiwaraSchapira}.

The (bounded) \emph{derived category of sheaves} $\Der^b(Z)$ has as objects
bounded complexes of sheaves
\begin{equation*}
\xymatrix @M+4pt { {\Sheaf \equiv} & {\cdots} \ar[r]^-{d_{n-2}} &
{\Sheaf^{n-1}} \ar[r]^-{d_{n-1}} & {\Sheaf^n} \ar[r]^-{d_n} & {\cdots} }
\end{equation*}
on $Z$.  A single sheaf will always be treated as a complex which is zero
except in degree $0$.  The definition of morphism is more subtle.  Before
giving it, recall that the \emph{cohomology sheaf} $\Hsheaf(\Sheaf)$ of
$\Sheaf$ is the complex of sheaves (with zero differential) defined by
$\Hsheaf^n(\Sheaf) \equiv \Ker d_n /\Im d_{n-1}$.  Any morphism $\psi\colon
\Sheaf \to \Sheaf'$ between complexes of sheaves induces a morphism
$\Hsheaf(\psi)\colon \Hsheaf(\Sheaf) \to \Hsheaf(\Sheaf')$ between the
cohomology sheaves; $\psi$ is called a \emph{quasi-isomorphism} if
$\Hsheaf(\psi)$ is an isomorphism.  Finally a morphism $\phi \in
\Mor_{\Der^b(Z)}$ is a certain equivalence class of
diagrams $\xymatrix@1{{\Sheaf} & {\Sheaf''} \ar[l]_{\phi_1} \ar[r]^{\phi_2}
& {\Sheaf'}}$ of sheaf complex morphisms such that $\phi_1$ is a
quasi-isomorphism.  (We will not define the equivalence relation here; to
begin with, $\phi_1$ and $\phi_2$ should be taken to be homotopy classes of
sheaf complex morphisms.)  One defines a notion of composition of such
morphisms.

It is clear that any sheaf complex morphism $\psi$ induces a morphism of
the corresponding objects in the derived category of sheaves; one sets
$\Sheaf'' = \Sheaf$ with $\phi_1$ the identity and $\phi_2=\psi$.  The
definition is arranged so that a quasi-isomorphism becomes an isomorphism
in the derived category and in fact we say that two complexes of sheaves
are quasi-isomorphic if they are isomorphic in the derived category.
Loosely put, the derived category is a localization of the ordinary
category of complexes of sheaves obtained by ``inverting'' the
quasi-isomorphisms.  The cohomology sheaf $\Hsheaf(\Sheaf)$ is well-defined
for an object in the derived category and the \emph{constructible derived
category} is the full subcategory consisting of objects $\Sheaf$ such that
$\Hsheaf(\Sheaf)|_{S^k}$ is a locally constant sheaf with finitely
generated stalks for every stratum $S^k$.

\subsection{Functors on the derived category}
Any bounded complex of sheaves $\Sheaf$ has an \emph{injective resolution}:
a quasi-isomorphism $\Sheaf \to \mathcal I$ where $\mathcal I$ is a complex
of injective sheaves; this is unique up to isomorphism in the derived
category.  Thus for any left exact functor $F$ on sheaves one may define a
corresponding functor $RF$ (the \emph{right derived functor}) on the
derived category of sheaves:
\begin{equation*}
RF(\Sheaf) \equiv F(\mathcal I).
\end{equation*}
We will write $RF$ simply as $F$ since we will never consider the original
functor $F$ on complexes of sheaves.  In this way, for any continuous map
$k:Z\to W$ we obtain functors $k_*$, $k_!$, and $k^*$ on the derived
category.  (The functor $k^*$ is actually exact so one does not even need
to pass to an injective resolution.)  If $k$ is a stratified map
\cite[\S1.2]{refnGoreskyMacPhersonIHTwo} these are in fact functors on the
constructible derived category.

The \emph{hypercohomology} of an object $\Sheaf$ of the derived category is
the right derived functor of the global sections functor $\G(Z;\cdot)$.
Thus $H(Z;\Sheaf) \equiv H(\Gamma(Z;\mathcal I))$.

The functor $k^*$ is a left adjoint to $k_*$:
\begin{equation*}
\Mor_{\Der^b(Z)}(k^*\Sheaf,\Sheaf') = \Mor_{\Der^b(W)}(\Sheaf, k_*\Sheaf').
\end{equation*}
It follows that there is a natural \emph{adjunction morphism}
\begin{equation*}
\Sheaf \longrightarrow  k_*k^*\Sheaf
\end{equation*}
characterized by the condition that $k^*\Sheaf \longrightarrow
k^*k_*k^*\Sheaf = k^* \Sheaf$ is the identity morphism.

With slightly more work, one may also define a functor $k^!$ on the derived
category which is a right adjoint to $k_!$.  We will only need $k^!$ when
$k$ is the inclusion of a locally closed subset; in this case $k^!=
k^*\circ R\G_Z$, the derived ``sections supported by $Z$'' functor.
Concretely, $k^!(\Sheaf)$ is the sheaf defined by
\begin{equation*}
k^!(\Sheaf)(U) \equiv \Ker \bigl(\G(\widetilde U;\mathcal I) \to
\G(\widetilde U \setminus (Z \cap \widetilde U);\mathcal I)\bigr),
\end{equation*}
where $\widetilde U\subseteq W$ is an open subset such that $U= Z \cap
\widetilde U$ is closed in $\widetilde U$.

It is not always necessary to use injective resolutions.  For example, for
hypercohomology it suffices to use a resolution by a complex of soft or
fine sheaves provided the space is paracompact (as we will assume).

\subsection{Sheaf-theoretic intersection cohomology}
\label{ssectSheafTheoreticIC}
We need one more functor on the derived category.  This is the truncation
functor $\t^{\leqslant n}$ defined by
\begin{equation*}
\xymatrix @M+4pt { {\t^{\leqslant n}\Sheaf \equiv} & {\cdots}
\ar[r]^-{d_{n-2}} & {\Sheaf^{n-1}} \ar[r]^-{d_{n-1}} & {\ker d_n}
\ar[r]^-{d_n} & {0} \ar[r] & {\cdots} }\ ;
\end{equation*}
it is characterized by the fact there is a natural morphism $\t^{\leqslant
n}\Sheaf \to \Sheaf$ which induces an isomorphism on cohomology sheaves in
degrees $\le n$ and the zero map on cohomology sheaves in degrees $>n$.

For all $k\ge 2$, let $U_k \equiv Z\setminus Z_{d-k}$, and let $j_k\colon
U_k\hookrightarrow U_{k+1}= U_k \cup S^k$ and $\i_k\colon S^k
\hookrightarrow U_{k+1}$ be the inclusions.  (We will find it convenient in
situations like this to use the same symbols to denote the inclusions of
the same sets into larger subspaces.)  The local calculation
\eqref{eqnLocalCalculation} is encoded into an object of the constructible
derived category by Deligne's formula \cite{refnBeilinsonBernsteinDeligne},
\cite{refnGoreskyMacPhersonIHTwo}:
\begin{equation}
\label{eqnIntersectionCohomologySheaf}
\IpC(Z;\EE) \equiv  \t^{\leqslant p(d)} j_{d*}\cdots\t^{\leqslant
  p(3)}j_{3*}\t^{\leqslant p(2)} j_{2*}\EE.
\end{equation}
This object is characterized up to isomorphism in the derived category by
the following conditions:
\begin{enumerate}
\item\label{itemCoefficientsOnSmooth} $\i_2^*\Sheaf \cong \EE$.
\item\label{itemLocalVanishing} $\Hsheaf^j(\i_k^*\Sheaf)=0$ for $j>p(k)$.
\item\label{itemLocalAttaching} The \emph{attaching morphism} $\i_k^*\Sheaf
  \to \i_k^*j_{k*}^{\vphantom{*}}j_k^*\Sheaf$ induces an isomorphism on
  cohomology sheaves in degrees $\le p(k)$.
\end{enumerate}
By definition, the intersection cohomology $I_pH(Z;\EE)$ is the
hypercohomology of $\IpC(Z;\EE)$.  It is possible to recast these
conditions in an
equivalent form that does not make explicit use of the stratification
\cite{refnBorelIntersectionCohomology}, \cite{refnGoreskyMacPhersonIHTwo};
this proves that intersection cohomology is a topological invariant.

In the piecewise-linear case, the allowable locally finite geometric chains
define a complex of sheaves $U \mapsto I_pC(U;\EE)$ which satisfies the
above conditions (essentially by \eqref{eqnLocalCalculation}) and hence is
quasi-isomorphic to $\IpC(Z;\EE)$.  These sheaves are also soft, so their
hypercohomology may be computed as the cohomology of global sections; it
follows that the two definitions of intersection cohomology agree.

As in \cite{refnBeilinsonBernsteinDeligne}, we will sometimes consider a
more general notion of perversity which is simply an integer-valued
function $p\colon \mathcal S \to \ZZ$ on the set of strata $\mathcal S$.
The formula \eqref{eqnIntersectionCohomologySheaf} may be generalized to
such $p$ and the corresponding intersection cohomology may be defined; in
general such a cohomology group depends upon the stratification.

\subsection{Sheaf-theoretic \boldmath $L^2$-cohomology}
If the  nonsingular stratum $S^0$ of $Z$ has a Riemannian metric, one may
localize the $L^2$-cohomology of $S^0$ at points of $Z$ to define the
\emph{$L^2$-cohomology sheaf} $\mathcal L_{(2)}(Z;\EE)$.  This is the
complex of sheaves defined by the presheaf
\begin{equation*}
U \mapsto A_{(2)}(U\cap S^0;\EE).
\end{equation*}
In general this sheaf may not be fine; the difficulty is finding a
partition of unity for a cover of $Z$ with functions $\eta$ such that
$\eta|_{S^0}$ has bounded differential.  However when this sheaf is fine,
the $L^2$-cohomology of $S^0$ may be expressed as a hypercohomology group
on $Z$: $H(Z;\mathcal L_{(2)}(Z;\EE)) = H_{(2)}(S^0;\EE)$.

Clearly the $L^2$-cohomology sheaf on the compactification $\Xstar$ of
$\SL_2(\ZZ)\back H$ considered in \S\ref{sectExample} is fine.  If instead
one considers the compactification $\Xbar$ obtained by adjoining
$S^1=\ZZ\back\RR$ (an example of the Borel-Serre compactification to be
considered later), the sheaf would not be fine since $|d\theta|\to \infty$
near the cusp by \eqref{eqnCuspMetric}.  It is not difficult to show that
$\mathcal L_{(2)}(Z;\EE)$ is fine \cite[(4.4)]{refnZuckerWarped},
\cite{refnZuckerLtwoIHTwo} when $Z$ is the reductive Borel-Serre
compactification $\Xhat$ or one of the Satake compactifications $\Xstar$ of
a locally symmetric space to be introduced later.

The more precise form of Zucker's conjecture and Theorem
~\ref{thmZuckersConjecture} is that there is an isomorphism
$\IpC(\Xstar;\EE) \cong \mathcal L_{(2)}(\Xstar;\EE)$ in the derived
category.

\section{Compactifications of locally symmetric spaces}
We need to consider three compactifications of a locally symmetric space
$X$ which fit into a diagram like this:
\begin{equation*}
\xymatrix @M+4pt @R-.4in{
{\Xbar}  \ar[r] & {\Xhat} \ar[r]^{\pi} & {\Xstar} \\
{\txt{Borel-Serre}} & {\txt{reductive\\Borel-Serre}} & {\txt{Satake}}
}
\end{equation*}

\begin{itemize}
\item The Borel-Serre compactification $\Xbar$ \cite{refnBorelSerre} is
``topological'': it corresponds to adding a boundary (or more generally
corners) at infinity and does not change the homotopy type of $X$.
\item The reductive Borel-Serre compactification $\Xhat$ (first constructed
  by Zucker in \cite[(4.1)]{refnZuckerWarped}) is ``geometric'': the
  boundary strata of $\Xbar$ become collapsed so that the Riemannian metric
  on $X$ extends to a nondegenerate metric on each stratum at infinity.
\item The compactification $\Xstar$ is one of the various Satake
compactifications \cite{refnSatakeCompact},
\cite{refnSatakeQuotientCompact}.  Zucker
\cite{refnZuckerSatakeCompactifications} has shown that it may be viewed as
a quotient of $\Xhat$.  In the special case where $X$ is Hermitian locally
symmetric, one of these is the Baily-Borel-Satake compactification.  This
is a complex-analytic compactification: the strata at infinity are further
collapsed so that they also become Hermitian locally symmetric.  In fact it
is also an algebraic compactification: Baily and Borel
\cite{refnBailyBorel} show in this case that $\Xstar$ is a projective
algebraic variety defined over a number field.
\end{itemize}

We will give the constructions of the three compactifications carefully in
\S\S\ref{sectBorelSerre}--\ref{sectSatake} together with examples.  In view
of their importance to intersection cohomology, we will highlight for each
compactification the boundary strata and their links.  We will not however
discuss the \emph{toroidal compactifications} $\widetilde X$
\cite{refnAshMumfordRapoportTai}; these are defined in the Hermitian case
and map to the Baily-Borel-Satake compactification $\Xstar$.

In the remainder of this section, we will describe the three
compactifications qualitatively.  First we give two simple examples.

\subsection{Example: \boldmath $X = \SL_2(\ZZ)\back H$}
\begin{figure}
\begin{equation*}
\begin{xy}
0;<2cm,0cm>:
(-1.5,0);(1.5,0)**\crv~l{},
(-.5,2.5);(-.5,.866025404)**\crv~l{},
(.5,.866025404);(.5,2.5)**\crv~l{},
(-.5,.866025404);(0,1.154700538),**{},(.5,.866025404)="rho",{\ellipse_{}},
(-1,0);(-1,1),**{},(-.5,.866025404),{\ellipse_{.}},
(1,0);(1,1),**{},(.5,.866025404),{\ellipse^{.}},
(-.5,.866025404);(-.5,0)**\crv{~*=<4mm>{}~**@{.}},
(.5,.866025404);(.5,0)**\crv{~*=<4mm>{}~**@{.}},
(0,1)*{\bullet}*+++!DC{\scriptstyle i},
(.5,.866025404)*{\bullet}*++!LD{\scriptstyle\rho},
(0,0)*{\bullet}*+++!UC{\scriptstyle 0},
(-1.5,2.75);(1.5,2.75)**\crv~l{},
(1.25,2.75)*+!DC{\txt{\scriptsize  $y= \infty$}},
(-1.5,2);(1.5,2)**\crv~l{},
(1.25,2)*+!DC{\txt{\scriptsize  $y=b$}},
\end{xy}
\end{equation*}
\caption{The   upper half-plane with an additional line at $y=\infty$}
\label{figHbar}
\end{figure}
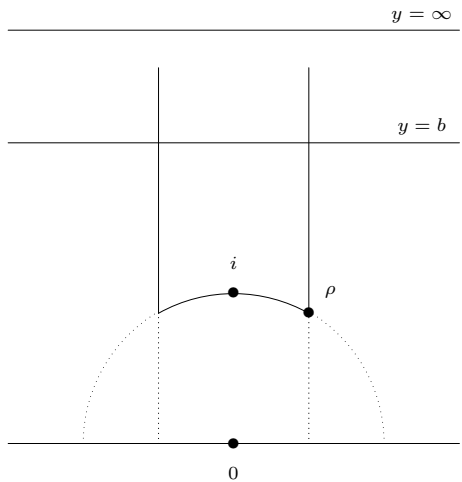

In this case there is only one non-proper stratum for all three
compactifications and the situation is quite simple.  A basic result of
reduction theory shows that for $y \gg 0$, the identifications induced by
$\SL_2(\ZZ)$ are already induced by its parabolic subgroup
$\left\{\,\left.\left(\begin{smallmatrix} 1&
n\\0&1\end{smallmatrix}\right)\right| n\in \ZZ\,\right\}$.  This group
preserves the lines $y=b$ and acts on them by translation by elements of
$\ZZ$.  Thus in $X$ for sufficiently large $y$ these lines wrap up into a
family of circles; for the compactification $\Xbar$ we adjoin yet another
circle at $y=\infty$; see Figure~\ref{figHbar}.
However the norm of the tangent vector to the circles in this family decays
as $y\to \infty$: $\left|\frac{\partial}{\partial x}\right| = \frac1{y} \to
0$.  Thus the circle at infinity is collapsed to a point in $\Xhat$ and
$\Xstar$.  In summary:
\begin{center}
\begin{tabular}{c|c|c}
& Boundary stratum & Link\\
\hline
$\Xbar^{\vphantom{A}}$ & $S^1$ & point\\
$\Xhat$ & point & $S^1$\\
$\Xstar$ & point & $S^1$
\end{tabular}
\end{center}

\subsection{Example: Hilbert modular surface}
The simplest Hermitian example where all three compactifications are
different is the Hilbert modular surface.  (See \cite[\S
I.5]{refnAshMumfordRapoportTai} for a nice short exposition and
\cite{refnHirzebruch} for more details.)  This is the moduli space of
complex tori $\CC^2/\Lambda$ with ``real multiplication'' by $\mathcal
O_k$, the ring of integers of a real quadratic extension $k=\QQ(\sqrt{d})$,
$d>0$, of $\QQ$.  The symmetric space is simply $H\times H$; the arithmetic
group $\G= \SL_2(\mathcal O_k)$ acts on $H\times H$ via the embedding
$\SL_2(\mathcal O_k) \hookrightarrow \SL_2(\RR)\times \SL_2(\RR)$ given by
the two real embeddings of $k$ in $\RR$.  Again, reduction theory shows
that ``near infinity'' corresponding to $y_1y_2 \gg 0$, the group $\G=
\SL_2(\mathcal O_k)$ acts via its parabolic subgroup
\begin{equation}
\left\{\,\left.\left(\begin{smallmatrix} u&
0\\0&u^{-1}\end{smallmatrix}\right)\right| u\in \mathcal
O_k^\times\,\right\} \ltimes
\left\{\,\left.\left(\begin{smallmatrix} 1&
a\\0&1\end{smallmatrix}\right)\right| a\in \mathcal O_k\,\right\}
\label{eqnHilbertParabolic}
\end{equation}
and there is only one non-proper stratum for all three compactifications.

\begin{figure}
\begin{minipage}{2.4in}
\def\latticebody{%
\drop{\cdot}}
\xy +<2.5cm,2.5cm>="o",0*\xybox{%
0;<.5cm,.5cm>:<.707106664cm,-.707106664cm>::,{"o"
\croplattice{-4}{4}{-4}{4}{-4.5}{4.5}{3.5}{-3.5}}
,"o"+(1,0)="a"*+!D{\scriptstyle(1,1)}
,"o"+(0,1)="b"*+!L{\scriptstyle(\delta,\bar\delta)}
,"o"+(1,1)="c"
,{"o" \ar "a"}
,{"o" \ar "b"}
,{"a" \ar @{-} "c"}
,{"b" \ar @{-} "c"}
}="L","o"."L"="L"
,{"L"+L \ar "L"+R*+!L{\scriptstyle x_1}}
,{"L"+D \ar "L"+U*+!D{\scriptstyle x_2}}
\endxy
\end{minipage}
\hfill
\begin{minipage}{2.4in}
\begin{xy}
0="o";<.5cm,0cm>:
,{0 \ar (7.071067812,7.071067812)}
,(1.414213562,1.414213562)*{\scriptstyle\bullet}
*+!L(1.2){\scriptstyle (1,1)}
,{0 \ar (9.995670057, .294245055)}
,(8.242640684, .2426406881)*{\scriptstyle\bullet}
*+!D{\scriptstyle (u_0^2,u_0^{-2})}
,{0 \ar (.294245055, 9.995670057)}
,(.2426406881, 8.242640684)*{\scriptstyle\bullet}
*+!L{\scriptstyle (u_0^{-2},u_0^2)}
,{(4.25,1.5)*+!LD{\txt{\scriptsize $y_1y_2=b$}} \ar (4,.7)}
,{(.2, 10.00000000);(10.0, .2000000000)
**\crv~c{
(.4, 5.000000000)&
(.6, 3.333333334)&
(.8, 2.500000000)&
(1.0, 2.000000000)&
(1.2, 1.666666667)&
(1.4, 1.428571429)&
(1.6, 1.250000000)&
(1.8, 1.111111111)&
(2.0, 1.000000000)&
(2.2, .9090909090)&
(2.4, .8333333334)&
(2.6, .7692307692)&
(2.8, .7142857142)&
(3.0, .6666666666)&
(3.2, .6250000000)&
(3.4, .5882352942)&
(3.6, .5555555556)&
(3.8, .5263157894)&
(4.0, .5000000000)&
(4.2, .4761904762)&
(4.4, .4545454546)&
(4.6, .4347826086)&
(4.8, .4166666666)&
(5.0, .4000000000)&
(5.2, .3846153846)&
(5.4, .3703703704)&
(5.6, .3571428572)&
(5.8, .3448275862)&
(6.0, .3333333334)&
(6.2, .3225806452)&
(6.4, .3125000000)&
(6.6, .3030303030)&
(6.8, .2941176470)&
(7.0, .2857142858)&
(7.2, .2777777778)&
(7.4, .2702702702)&
(7.6, .2631578948)&
(7.8, .2564102564)&
(8.0, .2500000000)&
(8.2, .2439024390)&
(8.4, .2380952380)&
(8.6, .2325581396)&
(8.8, .2272727272)&
(9.0, .2222222222)&
(9.2, .2173913044)&
(9.4, .2127659574)&
(9.6, .2083333334)&
(9.8, .2040816326)}}="L","o"."L"="L"
,{"L"+L \ar "L"+R*+!L{\scriptstyle y_1}}
,{(0,0) \ar (0,10)*+!D{\scriptstyle y_2}}
\end{xy}
\end{minipage}
\caption{$\SL_2(\mathcal O_k)$ acting on $H\times H$}
\label{figHilbert}
\end{figure}
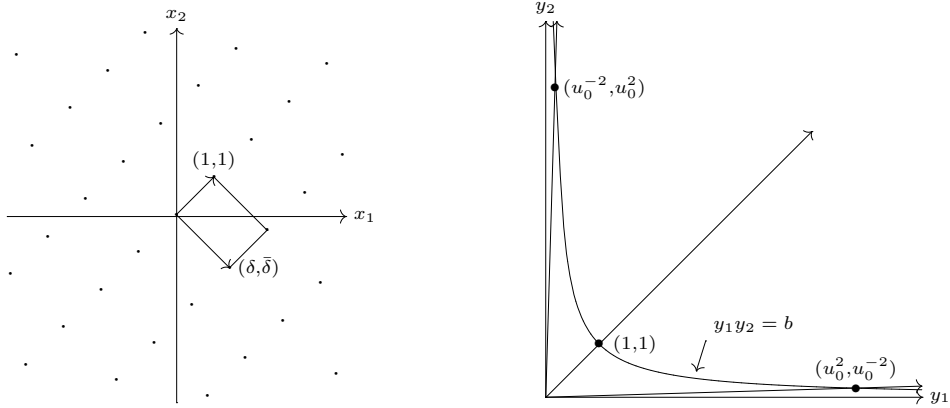

Write $\mathcal O_k = \ZZ + \ZZ\delta$; the image of the map $\mathcal O_k
\to \RR^2$ given by the two real embeddings of $k$ in $\RR$ is a lattice
$\ZZ(1,1) + \ZZ(\delta,\bar\delta)$.  The second factor of
\eqref{eqnHilbertParabolic} acts in the $x_1x_2$-plane via translation by
this lattice.  Let $u_0$ be a generator of the group of positive units
$\mathcal O_{k+}^\times$.  The first factor of \eqref{eqnHilbertParabolic}
preserves the hyperbolas $y_1y_2=b$ and acts on them via $(y_1,y_2)\mapsto
(u_0^{2k},u_0^{-2k})$ for all $k\in \ZZ$.  Furthermore it effects an
automorphism of the lattice in the $x_1x_2$-plane.  Thus in $X$ for
sufficiently large $y_1y_2$, we obtain a family of flat $T^2$-bundles (from
the $x_1x_2$-plane) over $S^1$'s (from the hyperbolas); for the
compactification $\Xbar$ we adjoin another copy of this bundle at
$y_1y_2=\infty$; see Figure ~\ref{figHilbert}.

To understand the non-proper stratum in the other two compactifications
first note that the product metric on $H\times H$ may be rewritten as
\begin{equation*}
dr^2 + dt^2 + e^{-r}(e^{-t}dx_1^2 + e^t dx_2^2) \qquad\text{where $r=\log
    y_1y_2$ and $t = \log y_1/y_2$;}
\end{equation*}
up to quasi-isometry this induces the metric
\begin{equation*}
dr^2 + ds_{S^1}^2 + e^{-r}ds_{T^2}^2
\end{equation*}
on $X$ where $r \gg 0$.  We see that the norm of a tangent vector to the
$T^2$-fibers decays to $0$ as $r\to \infty$ while the norm of a tangent
vector to the base $S^1$ does not.  Thus the $T^2$-fibers are collapsed in
$\Xhat$ while the $S^1$ remains; since $S^1$ cannot have a complex
structure, however, it is collapsed in $\Xstar$.  In summary:
\begin{center}
\begin{tabular}{c|c|c}
& Boundary stratum & Link\\
\hline
$\Xbar^{\vphantom{A}}$ & flat $T^2$-bundle over $S^1$ & point\\
$\Xhat$ & $S^1$ & $T^2$\\
$\Xstar$ & point &  flat $T^2$-bundle over $S^1$
\end{tabular}
\end{center}

\subsection{The general case}
\label{ssectGeneralCase}
In the general case the strata are no longer compact.  For both the
Borel-Serre compactification $\Xbar$ and the reductive Borel-Serre
compactification $\Xhat$ they are indexed by the set $\Pl$ of
$\G$-conjugacy classes of parabolic $\QQ$-subgroups $P$ of $G$; for every
$P\in\Pl$, we let $Y_P$ (respectively $X_P$) denote the corresponding
stratum of $\Xbar$ (respectively $\Xhat$).  (Note that we will use the same
symbol, say $P$, to denote both a parabolic $\QQ$-subgroup and its
$\G$-conjugacy class.)

Associated to $G\in \Pl$ we have $X_G=Y_G=X$, the top-dimensional open
stratum, while associated to a minimal parabolic $\QQ$-subgroup $P_0\in\Pl$
is a minimal-dimensional stratum.  In fact the correspondence between the
strata and $\Pl$ is order-preserving if we define the partial order on the
set of strata by $S<T$ if and only if $S\subseteq \cl{T}$ (topological
closure) and the partial order on $\Pl$ by $P<Q$ if and only if there
exists $\g\in\G$ so that $P \subset \lsp\g Q \equiv \g Q \g^{-1}$.

For $P\in\Pl$, let $r(P)$ denote the \emph{parabolic rank}; this is the
number of proper parabolic $\QQ$-subgroups in a maximal chain
$P=P_r<P_{r-1}< \cdots P_1<G$.  Let $|\D_P|^\circ$ denote an open simplex
of dimension $r(P)-1$.  We can now crudely describe the strata and their
links.  Each stratum $X_P$ of $\Xhat$ will be a locally symmetric space; in
the above example this was $S^1$.  Each stratum $Y_P$ of $\Xbar$ is a flat
bundle over $X_P$ with general fiber a compact nilmanifold $\nil_P$; in the
above example the nilmanifold was $T^2$.  The links are now themselves
stratified.

The corresponding data for a Satake compactification of $X=\G\back D$ is in
general more complicated; we will give Zucker's realization of $\Xstar$ as
a quotient of $\Xhat$ \cite{refnZuckerSatakeCompactifications}.  The
construction depends on the choice of a certain representation $U$ of $G$
which must satisfy a condition called ``geometric rationality''
\cite{refnCasselmanGeometricRationality}.  This data induces a
decomposition $X_P = X_{P,\l}\times X_{P,h}$ (modulo a finite quotient) for
each stratum $X_P$ of $\Xhat$.  The strata of $\Xstar$ will be the various
$X_{P,h}$, and the quotient map $\Xhat\xrightarrow{\pi} \Xstar$ is induced
by projection, stratum by stratum.  The situation is complicated by the
fact that several $P\in \Pl$ may yield the same stratum in $\Xstar$; we
index the boundary stratum $X_{P,h}$ of $\Xstar$ by the maximal such $P$,
denoted $P^\dag$.  Let $\Pl^*= \{\, P^\dag \mid P\in \Pl\,\}$.

We can summarize the above discussion in the following table;  for
simplicity we indicate just the smooth locus of the links.

\begin{center}
\begin{tabular}{c|c|c}
& Boundary stratum associated to $P$ & Link$^\circ$\\
\hline
$\Xbar^{\vphantom{A}}$ & $Y_P =$ flat $\nil_P$-bundle over $X_P$ &
$|\D_P|^\circ$\\
$\Xhat$ & $X_P$ & $\nil_P\times|\D_P|^\circ$\\
$\Xstar$ & $X_{P^\dag,h}$ & 
$(\nil_{P^\dag}\times|\D_{P^\dag}|^\circ)$-bundle over $X_{P^\dag,\l}$
\end{tabular}
\end{center}

\section{The Borel-Serre compactification}
\label{sectBorelSerre}
As in \S\ref{sectZuckersConjecture}, we consider a locally symmetric space
$X=\G\back D$ where $D= G/KA_G$.  In this section we outline the actual
construction of $\Xbar$ due to Borel and Serre \cite{refnBorelSerre}.  The
construction uses the language of algebraic groups, in particular their
parabolic subgroups, for which good references are the texts of Borel
\cite{refnBorelLAG} and Humphreys
\cite{refnHumphreysLinearAlgebraicGroups}, as well as the paper by Borel
and Tits \cite{refnBorelTits}.  For the reduction theory of arithmetic
groups, one may consult Borel's book \cite{refnBorelArithmetiques}.  We
illustrate the discussion with the example of $G=\GL_n(\RR)$ and
$\G=\GL_n(\ZZ)$.  In this case $A_G$ is the subgroup of positive scalar
matrices and we may pick $K=\O(n)$.  The symmetric space
$D=\GL_n(\RR)/\O(n)A_G$ is the space of positive definite $n\times n$
symmetric matrices modulo homothety and $X=\G\back D$ is the moduli of
positive definite quadratic forms on $\ZZ^n$ modulo homothety and isometry.

Let $P$ be a parabolic $\QQ$-subgroup, let $N_P$ be its unipotent radical
and $L_P = P/N_P$ its Levi quotient.  We can when necessary lift $L_P$ to a
subgroup $\widetilde L_P$ of $P$ such that $P= \widetilde L_P N_P$.  In the
case $G=\GL_n(\RR)$, any such $P$ is conjugate under $\GL_n(\QQ)$ to the
group of block upper-triangular matrices associated to a partition
$n=n_0+\dots + n_r$, $n_i\ge0$.  A lift of the Levi quotient is the group
of block diagonal matrices, $\GL_{n_0}(\RR)\times\dots \times
\GL_{n_r}(\RR)$, while the unipotent radical $N_P$ consists of those block
upper-triangular matrices with identity matrices along the diagonal; see
Figure ~\ref{figParabolic}.  Conjugation by elements of $N_P$ would take
one lift $\widetilde L_P$ to all possible such lifts.  The stabilizer of a
basepoint $x_0\in D$ has the form $KA_G$ where $K$ is a maximal compact
subgroup.  Thus the choice of a basepoint determines a maximal compact
subgroup $K$ with a corresponding Cartan involution $\theta$; there exists
a unique lift $\widetilde L_P$ of $L_P$ which is stable under $\theta$.
Such a lift however is not necessarily defined over $\QQ$.

\begin{figure}
\begin{equation*}
\def\objectstyle{\scriptstyle}
\left(
\vcenter{\xymatrix @M=0pt @!0{
{\GL_{n_0}(\RR)} & {}  & {} &
*!L(2){\txt{\Large $*$}} \\
{}  & *+++{\GL_{n_1}(\RR)}  \POS[];[rrdd]**\dir{} ?<<<;?>>>**\dir*{\cdot} & {} & {} \\
{} & {} & {} & {} \\
*!R(2){\txt{\Large $0$}} & {}  & {} & {\GL_{n_r}(\RR)} 
}}
\right)
= 
\left(
\vcenter{\xymatrix @M=0pt @!0{
{\GL_{n_0}(\RR)} & {}  & {} & *!L(2){\txt{\Large $0$}} \\
{}  & {\GL_{n_1}(\RR)} \POS[];[rrdd]**\dir{} ?<<<;?>>>**\dir*{\cdot} & {} & {} \\
{} & {} & {} & {} \\
*!R(2){\txt{\Large $0$}} & {}  & {} & {\GL_{n_r}(\RR)} 
}}
\right)
\left(
\vcenter{\xymatrix @M=0pt @!0{
{I_{n_0}} & {} & {} & {\txt{\Large $*$}} \\
{} & {I_{n_1}}
\POS[];[rrdd]**\dir{} ?<<<;?>>>**\dir*{\cdot} & {} & {} \\
{} & {} & {} & {} \\
{\txt{\Large $0$}} & {} & {} & {I_{n_r}} 
}}
\right)
\end{equation*}
\caption{A parabolic subgroup $P=\widetilde L_P N_P$ of $\GL_n(\RR)$}
\label{figParabolic}
\end{figure}

Let $A_P$ be the identity component of the maximal $\QQ$-split torus in the
center of $L_P$; it has a natural complement $M_P=\lsp0L_P$ defined as
\begin{equation*}
\lsp0L_P \equiv \bigcap_{\al\in X(L_P)} \ker \al^2,
\end{equation*}
where $X(L_P)$ denotes the (rational) characters of $L_P$ defined over
$\QQ$.  Thus $L_P = M_P A_P$ and this decomposition lifts to $\widetilde
L_P = \widetilde M_P \widetilde A_P$.  In the example, $\widetilde A_P$
consists of diagonal matrices which are positive scalar matrices in each
block and $\widetilde M_P$ consists of block diagonal matrices
$\SL^\pm_{n_0}(\RR)\times\dots \times \SL^\pm_{n_r}(\RR)$, where
$\SL^\pm_{n_i}(\RR)$ denotes matrices with determinant $\pm 1$.

There is a natural inclusion $A_G\subseteq A_P$ and a complement
$A_P^G\subseteq A_P$ given as the joint kernel of the characters in $X(G)$.
In our standard example, $a =\diag(\lambda_0\cdot I_{n_0},\dots,
\lambda_r\cdot I_{n_r})\in A_P$ lies in $A_G$ if and only if
$\lambda_0=\dots=\lambda_r$ and $a$ lies in $A_P^G$ if and only if $\det a
=1$.  Since $A_P^G$ is $\QQ$-split, it has coordinates given by characters
of $P$ defined over $\QQ$.  In the next paragraph we will define an
important and canonical set of characters $\D_P$ so that we have an
isomorphism
\begin{equation}
\label{eqnSplitCoordinates}
A_P^G \longtildearrow (\RR^*_+)^{\D_P},\qquad a \mapsto
(a^{\al_1},\dots,a^{\al_r}).
\end{equation}

There is an adjoint action of $z\in L_P$ on the Lie algebra $\n_P$ of $N_P$
given by lifting $z$ to $\widetilde L_P$ and then acting by the adjoint
action.  Although this action depends on the lift, its $A_P^G$-weights
(which are induced from characters of $P$) do not.  For $P$ a minimal
parabolic $\QQ$-subgroup these are the positive roots of a root system,
namely the $\QQ$-root system of $G$, but this does not hold true in
general.  However it is still true that if we let
$\D_P=\{\al_1,\dots,\al_r\}$ be the indecomposable elements among the
$A_P^G$-weights of $\n_P$, then all $A_P^G$-weights of $\n_P$ can be
expressed uniquely as non-negative integral linear combinations of elements
of $\D_P$.  We call the weights in $\D_P$ the \emph{simple roots} even
though they may not generate a root system.  There is also a ``dual'' basis
$\Dhat_P=\{\b_1,\dots,\b_r\}$ which in the case of a minimal parabolic
$\QQ$-subgroup is simply the set of fundamental weights.  In our example,
we have
\begin{equation*}
a^{\al_i} = \lambda_i/\lambda_{i-1} \quad\text{and}\quad
a^{\b_i} = 1/(\lambda_0\cdots\lambda_{i-1})
\end{equation*}
for $a=\diag(\lambda_0\cdot I_{n_0},\dots, \lambda_r\cdot I_{n_r})\in
A_P^G$.

For later use we single out two subsets of $A_P^G$, the \emph{strictly
dominant cone}
\begin{equation*}
A_P^{G+} \equiv \{ a\in A_P^G \mid a^\al > 1 \text{ for all
  $\al\in\D_P$}\},
\end{equation*}
and the \emph{negative codominant cone},
\begin{equation*}
\lsp-A_P^G \equiv \{ a\in A_P^G \mid a^\b \le 1 \text{ for all
  $\b\in\Dhat_P$}\}.
\end{equation*}
The Lie algebra $\sa_P^G$ of $A_P^G$, endowed with a Weyl group invariant
inner product, is illustrated in Figure ~\ref{figSimpleRoots}(a) in the
case of a minimal parabolic $\QQ$-subgroup of $\GL_3(\RR)$ together with
the subsets $\sa_P^{G+}\equiv \log A_P^{G+}$ and $\lsp-\sa_P^G \equiv \log
\lsp-A_P^G$.  The figure also includes $\D_P$ and $\Dhat_P$, identified
with elements of $\sa_P^{G*}$ and then transfered to $\sa_P^G$ via the
inner product.
\begin{figure}
\begin{equation*}
\begin{xy}
0;<.5in,0cm>:
,(-.866025,-1.5)="o"
%
,{"o";"o"+(1.73205,1)**[lightgrey]@{.}?(.1)@+c?(.2)@+c?(.3)@+c?(.4)@+c?(.5)@+c?(.6)@+c?(.7)@+c?(.8)@+c?(.9)@+c?(1)@+c}
,{"o";"o"+(0,2)**[lightgrey]@{.}?(.1)@+c?(.2)@+c?(.3)@+c?(.4)@+c?(.5)@+c?(.6)@+c?(.7)@+c?(.8)@+c?(.9)@+c?(1)@+c}
,s1;s{11}**[lightgrey]@{.}
,s2;s{12}**[lightgrey]@{.}
,s3;s{13}**[lightgrey]@{.}
,s4;s{14}**[lightgrey]@{.}
,s5;s{15}**[lightgrey]@{.}
,s6;s{16}**[lightgrey]@{.}
,s7;s{17}**[lightgrey]@{.}
,s8;s{18}**[lightgrey]@{.}
,s9;s{19}**[lightgrey]@{.}
,"o"*+!C!/va(240)2.45cm/{\sa_P^{G+}}
,{"o";"o"+(.866025,.5)**@{-}?>*\dir{>}?(1)/4mm/*{\b_1}}
,{"o";"o"+(0,1)**@{-}?>*\dir{>}?(1)/4mm/*{\b_2}}
,{"o";"o"+(1.73205,0)**@{-}?>*\dir{>}?(1)/4mm/*{\al_1}}
,{"o";"o"+(-.866025,1.5)**@{-}?>*\dir{>}?(1)/4mm/*{\al_2}}
@i,{"o";"o"+(-1.73205,0)**[lightgrey]@{-}?(.1)@+c?(.3)@+c?(.5)@+c?(.7)@+c?(.9)@+c}
,{"o";"o"+(.866025,-1.5)**[lightgrey]@{-}?(.1)@+c?(.3)@+c?(.5)@+c?(.7)@+c?(.9)@+c}
,s0;s{5}**[lightgrey]@{.}
,s1;s{6}**[lightgrey]@{.}
,s2;s{7}**[lightgrey]@{.}
,s3;s{8}**[lightgrey]@{.}
,s4;s{9}**[lightgrey]@{.}
,"o"*+!C!/va(60)1.5cm/{\lsp-\sa_P^{G}}
,"o"-(0,1.75)*++!C{\txt{(a)}}
\end{xy}
\quad\qquad
\begin{xy}
0;<.5in,0cm>:
,(-.866025,-1.5)="o"
%
,{"o"-(1.73205,1);"o"+(1.73205,1) *+{\sa_{Q_2}^G} **@{-}}
,{"o"-(0,1.75)*++!C{\txt{(b)}};"o"+(0,1.75) *+{\sa_{Q_1}^G} **@{-}}
,{"o"-(1.73205,0);"o"+(1.73205,0) *+{\sa_P^{Q_1}} **@{-}}
,{"o"-(-.866025,1.5);"o"+(-.866025,1.5) *+{\sa_P^{Q_2}} **@{-}}
\end{xy}
\end{equation*}
\caption{$\sa_P^G$, its subcones, and its subspaces for $P$ a minimal
parabolic subgroup of $\GL_3(\RR)$.  (Here $\D_P=\{\al_1,\al_2\}$ and
$\D_P^{Q_i}=\{\al_i\}$.)}
\label{figSimpleRoots}
\end{figure}

The parabolic $\QQ$-subgroups $Q$ containing $P$ are in one-to-one order
preserving correspondence with subsets $\D_P^Q$ of $\D_P$.  Namely,
consider $Q = \widetilde L_Q N_Q \supseteq P=\widetilde L_P N_P$; we can
decompose $N_P = (N_P\cap\widetilde L_Q)N_Q$.  Define $\D_P^Q$ to consist
of those $\al\in \D_P$ whose weight space lies in $\n_P\cap \widetilde
\levi_Q$ as opposed to $\n_Q$.  For such a $Q$, there is a natural
identification
\begin{equation*}
A_Q^G \equiv \bigcap_{\al\in\D_P^Q} \ker \al\subseteq A_P^G.
\end{equation*}
The set of simple roots $\D_Q$ for $Q$ is the set of restrictions
$\al|_{A_Q^G}$ for $\al\in \D_P \setminus \D_P^Q$.  Furthermore there is a
complementary subgroup
\begin{equation*}
A_P^Q \equiv \bigcap_{\al\in\D_P\setminus \D_P^Q}
\ker \b_\al\subseteq A_P^G
\end{equation*}
to $A_Q^G$; here $\b_\al\in\Dhat_P$ is the dual element to $\al\in\D_P$.
The group $A_P^Q$ may also be viewed as the analogue of $A_P^G$ for the
parabolic $\QQ$-subgroup $P/N_Q$ of $L_Q$.  We have a direct decomposition
\begin{equation*}
A_P^G = A_Q^G \times A_P^Q.
\end{equation*}
The corresponding decomposition of the Lie algebra in the example of
$\GL_3(\RR)$ is illustrated in Figure ~\ref{figSimpleRoots}(b).  Note that
the closure of the strictly dominant cone $A_P^{G+}$ is stratified with
strata $\{A_Q^{G+}\}_{Q\supseteq P}$, while $\{\lsp-A_P^Q\}_{Q\supseteq P}$
are the closures of the strata of the negative codominant cone
$\lsp-A_P^G$.

\medskip
The compactification $\Xbar$ will be the quotient by $\G$ of a
bordification $\Dbar$ which will be constructed in three steps.

\subsection{}
The first step is to form a partial bordification $\Abar_P^G$ of $A_P^G$, a
``corner'', by allowing the coordinates in \eqref{eqnSplitCoordinates} to
attain infinity:
\begin{equation*}
\Abar_P^G \longtildearrow (\RR^*_+ \cup \{\infty\})^{\D_P},\qquad a \mapsto
(a^{\al_1},\dots,a^{\al_r}).
\end{equation*}
The corner $\Abar_P^G$ is given a real-analytic structure by means of the
inverses of these coordinates: $\Abar_P^G \to (\RR_{\ge0})^{\D_P}$,
$a\mapsto (a^{-\al_1},\dots,a^{-\al_r})$.  For $Q\supseteq P$, the
inclusion $A_Q^G \hookrightarrow A_P^G$ extends to the partial
bordifications, $\Abar_Q^G \hookrightarrow \Abar_P^G$.  The action of
$A_P^G$ on itself by multiplication extends naturally to an action of
$A_P^G$ on $\Abar_P^G$.  In fact $\Abar_P^G$ is a semigroup.  Set
$\Abar_P^{G+}$ to be the subset where all coordinates are greater than $1$.

Let $o_Q \in \Abar_P^G$ be the point with coordinates
$\{\infty\}^{\D_P\setminus \D_P^Q}\times \{1\}^{\D_P^Q}$; this point has
stabilizer $A_Q^G$ so the orbit it generates is $A_P^Q\cdot o_Q$.  The
space $\Abar_P^G$ is stratified by these $A_P^G$-orbits:
\begin{equation*}
\Abar_P^G = \coprod_{Q\supseteq P} A_P^Q\cdot o_Q.
\end{equation*}
Topologically $\Abar_P^G$ is homeomorphic to a cone $c(|\D_P|)$ on the
simplex $|\D_P|$.  (We denote by $|\D_P|$ the $(r-1)$-dimensional closed
simplex whose vertices are identified with $\al\in \D_P$.)  The vertex of
the cone is $o_P = \{\infty\}^{\D_P}$ and the sub-cone $c(|\D_P^Q|)$
corresponds to the closure of the orbit $A_P^Q\cdot o_Q$.

All of the above illustrated in Figure ~\ref{figSplitComponent} for $P$ a
minimal parabolic subgroup in $G=\GL_3(\RR)$.  The figure is drawn using
analytic coordinates based at $o_P$ and chosen so that at $o_G$ the diagram
accurately reflects the geometry of a Weyl-group invariant inner product on
the tangent space $\sa_P^G$.  It is important to note that whereas the two
$1$-parameter subgroups $A^G_{Q_1}$ and $A^G_{Q_2}$ each intersect the
boundary at infinity in $1$ point, neither of the subgroups $A_P^{Q_1}$ or
$A_P^{Q_2}$ intersect the boundary at all.  This is due to the fact that
the two simple roots are not orthogonal.

\begin{figure}
\begin{equation*}\begin{xy}
0="o";<1in,0cm>:
,{(-.866025,-1.5); (-2.508,-2.448) **[lightgrey]@{-}}
,{(-.866025,-1.5); (-.866025,-3) **[lightgrey]@{-}}
,{(-.866025,-1.5);
(-.1003,  -2.9956)
**[lightgrey]\crv~c{
(-.641349641230121392, -1.53231485280888058)&
(-.488298956031021303, -1.61366865210952337)&
(-.384427196404046434, -1.72287163012768940)&
(-.308697502240292465, -1.85316467241659798)&
(-.259106279177597208, -1.97780620323190814)&
(-.216341843634927338, -2.12566528165960511)&
(-.183233298864616617, -2.27980769442398934)&
(-.158027924900405958, -2.43222120498072724)&
(-.139387492802535318, -2.57308012409121067)&
(-.121217765341058770, -2.74288183585022516)}}
,{(-.866025,-1.5);
(-2.49767, -1.562256)
**[lightgrey]\crv~c{
(-.864229646762526671, -1.50312327712821348)&
(-1.00634876851364274, -1.32158250841776104)&
(-1.15332856800193029, -1.22971362661904670)&
(-1.29983700094806286, -1.19435953305537866)&
(-1.45053893258849497, -1.19392221523319584)&
(-1.58327727617248203, -1.21329572166381627)&
(-1.73270921204235817, -1.25019017331921045)&
(-1.88275472968209722, -1.29858853884797920)&
(-2.02735138868630526, -1.35296679996145474)&
(-2.15865900703453661, -1.40725317178242149)&
(-2.31479646675466410, -1.47641858210045029)}}
,{(-.866025,-1.5) *{\scriptscriptstyle\bullet}
             *+!C!/va(60)3mm/{\scriptstyle o_G};
             (0,-1) **@{-} *{\scriptscriptstyle\bullet}
             ?<(.5)*+!U{\scriptscriptstyle A_{Q_2}^{G+}}
             ?(1)/4mm/*{\scriptstyle o_{Q_2}}}
,{(-.866025,-1.5); (-.866025,-.5) **@{-} *{\scriptscriptstyle\bullet}
             ?(.5)*+!R{\scriptscriptstyle A_{Q_1}^{G+}}
             ?(1)/4mm/*{\scriptstyle o_{Q_1}}}
,{(-2.483, -2.0242);
(-.866025,-1.5)
**\crv~c{
(-1.93792145972396912, -1.78735303963998615)&
(-1.22900262672383498, -1.54900354857376321)&
(-.869628119631486228, -1.50000646736410670)}?>>(.5)*+!D{\scriptscriptstyle
  \lsp-A_P^{Q_1}}}
,{(-.52597, -3.01473);
(-.866025,-1.5)
**\crv~c{
(-.726975110255204115, -1.83884927034752632)}?>>(.5)*+!L{\scriptscriptstyle
  \lsp-A_P^{Q_2}}}
,{0 *{\scriptscriptstyle\bullet}; (-2.4312557,-1.4036861) **@{-}
             ?(0)/-2mm/*!LD{\scriptstyle o_P= A_P^P\cdot \,o_P}
             ?(.65)*+!RD{\scriptstyle A_P^{Q_1}\cdot \,o_Q}}
,{0; (0,-3) **@{-} ?(.65)*+!LD{\scriptstyle A_P^{Q_2}\cdot \,o_Q}}
\end{xy}
\end{equation*}
\caption{$\Abar_P^G$ for $P$ a minimal parabolic subgroup of $\GL_3(\RR)$}
\label{figSplitComponent}
\end{figure}

\subsection{}
The second step in the construction of $\Dbar$ is to extend the partial
bordification of $A_P^G$ to one of $D$.  To do this, note that even though
$A_P$ is not naturally a subgroup of $G$, Borel and Serre define an action
of $A_P$ on $D$ called the \emph{geodesic action}
\cite[\S3]{refnBorelSerre}.  To describe it, fix a basepoint in $D$ which
determines a maximal compact subgroup $K\subseteq G$, a Cartan involution
$\theta$, and a $\theta$-stable lift $\widetilde L_P$ of $L_P$.  Note that
$G=PK$, that is, $P$ acts transitively on $D$; this is a generalization of
Gram-Schmidt orthonormalization.  Every point $x\in D$ can thus be written
as $x=p K A_G$ for some $p\in P$, and the geodesic action of $a\in A_P$ on
$x$ is given by
\begin{equation*}
a \geo x \equiv p \tilde a K A_G,
\end{equation*}
where $\tilde a\in \tilde A_P$ is the lift of $a$.  The geodesic action is
independent of the choice of basepoint and commutes with the usual action
of $P$.

Let $e_P\equiv A_P^G\back D$ be obtained by collapsing the orbits of the
geodesic action of $A_P^G$ to points.  We will further discuss the
structure of $e_P$ below.  The quotient map $r_P\colon D \to e_P$ is a
$A_P^G$-principal bundle; natural trivializing sections are given by the
orbits of $\lsp 0P \equiv \bigcap_{\al\in X(P)} \ker \al^2$.  Thus if we
let $\mathscr A_P^G \equiv \lsp 0P\back D$ be the $A_P^G$-homogeneous space
obtained by collapsing these sections to points, there is a natural product
decomposition
\begin{equation*}
D = \mathscr A_P^G \times e_P.
\end{equation*}
This generalizes the product decomposition of the upper half-plane into
geodesics tending to $\infty$ (rays parallel to the positive $y$-axis)
times horocycles at $\infty$ (lines parallel to the $x$-axis).  For $P$ a
minimal parabolic $\RR$-subgroup this is the decomposition induced from the
Iwasawa decomposition of $G$.

The space $\mathscr A_P^G$ is an affine version of $A_P^G$; given a choice
of basepoint in $D$, there is a canonical isomorphism of
$A_P^G$-homogeneous spaces $\mathscr A_P^G \cong A_P^G$ such that the class
of the basepoint maps to the identity.  We can then define a bordification
$\bar{\mathscr A}_P^G$ such that this isomorphism extends to an isomorphism
$\bar{\mathscr A}_P^G \cong \Abar_P^G$.  The resulting $A_P^G$-space (in
fact $\Abar_P^G$-space) is independent of the choice of basepoint.  Note
that although the points $o_Q\in \Abar_P^G$ do not uniquely determine
points in $\bar{\mathscr A}_P^G$ (except for $Q=P$), the orbits $A_P^Q\cdot
o_Q\subseteq \Abar_P^G$ do indeed define well-defined orbits in
$\bar{\mathscr A}_P^G$.  We will again use the notation $A_P^Q\cdot o_Q$ to
denote the corresponding orbit of $\bar{\mathscr A}_P^G$.

Borel and Serre construct a partial bordification $D(P)$ by replacing
$\mathscr A_P^G$ in the above decomposition with $\bar{\mathscr A}_P^G$:
\begin{equation}
\label{eqnPartialBordification}
D(P) = \bar{\mathscr A}_P^G \times e_P.
\end{equation}
We consider $e_P$ as a stratum of $D(P)$ by identifying it with $o_P \times
e_P= \{\infty\}^{\D_P}\times e_P$; the link of any point in $e_P$ is
$|\D_P|$.  The projection map
\begin{equation}
\label{eqnGeodesicRetraction}
D(P) = \bar{\mathscr A}_P^G \times e_P \longrightarrow
\{\infty\}^{\D_P}\times e_P \equiv e_P
\end{equation}
is called \emph{geodesic retraction}.

\subsection{}
\label{ssectBorelSerreThirdStep}
The third and final step in the construction of $\Dbar$ is to note that the
construction above of the partial bordification is functorial in the sense
that for $P \subseteq Q$ we can canonically identify $D(Q)$ with an open
subset of $D(P)$.  Indeed $e_Q \subseteq D(Q)$ corresponds to $(A_P^Q\cdot
o_Q) \times e_P \subseteq D(P)$ and so $D(P)$ may be re-expressed as
$\coprod_{Q\supseteq P} e_Q$.  Borel and Serre then define a full
bordification of $D$ by setting
\begin{equation*}
\Dbar = \bigcup_P D(P) = \coprod_P e_P
\end{equation*}
where the union is over all parabolic $\QQ$-subgroups $P$ and the above
identifications are taken into account.  Let $\ebar_P$ denote the closure
of $e_P$ in $\Dbar$.

\medskip
The action of $G$ on $D$ extends to an action of $G(\QQ)$ on $\Dbar$ and
the \emph{Borel-Serre compactification} is defined by
\begin{equation*}
\Xbar \equiv \G \back \Dbar.
\end{equation*}
Given a parabolic $\QQ$-subgroup $P$, the various strata $e_{\lsp \g P}$ of
$\Dbar$ for $\gamma \in \G$ become identified to one stratum $Y_P$ in
$\Xbar$.  Furthermore reduction theory shows that the identifications
induced by $\G$ on $e_P$, and indeed even on a small neighborhood of
$\ebar_P$, are already induced by $\G_P \equiv \G \cap P$.  Thus there is a
stratification
\begin{equation*}
\Xbar = \coprod_{P\in \Pl} Y_P,\qquad \text{ where $Y_P \equiv \G_P \back
e_P$.}
\end{equation*}

The geodesic retraction \eqref{eqnGeodesicRetraction} extends to a
projection, also called  geodesic retraction,
\begin{equation}
\label{eqnQuotientGeodesicRetraction}
V\longrightarrow Y_P,
\end{equation}
where $V$ is a sufficiently small neighborhood of $Y_P$ in $\Xbar$.  In
fact, if $O_P$ is a relatively compact open subset of $Y_P$ and $\widetilde
O_P$ is its inverse image in $e_P$, reduction theory shows that there
exists $s_P\in \mathscr A_P^G$ such that the identifications induced by
$\G$ on
\begin{equation*}
\widetilde W_P \equiv (\Abar_P^{G+}\cdot s_P) \times \widetilde
O_P\subseteq D(P)
\end{equation*}
agree with those induced by $\G_P$.  However $\G_P$ acts only on the second
factor.  Thus
\begin{equation}
\label{eqnCylindricalSet}
W_P \equiv \G_P\back \widetilde W_P = (\Abar_P^{G+}\cdot s_P) \times
O_P\subseteq \G_P\back D(P)
\end{equation}
may be identified with a subset of $\Xbar$; the geodesic retraction is
projection on the second factor.  The open subset $W_P$ is called a
\emph{cylindrical set}; for appropriate choices of $O_P$ and $s_P$ (for all
$P\in\Pl$) one obtains a \emph{cylindrical cover}
\begin{equation}
\label{eqnCylindricalCover}
\Xbar = \bigcup_{P\in\Pl} W_P.
\end{equation}

Let $x\in Y_P$.  If $O_P$ is a ball neighborhood $B$ of $x$ in $Y_P$ the
construction above yields a neighborhood in $\Xbar$ of $x$ which is
homeomorphic (as a stratified space) to $B\times c(|\D_P|)$.  (Here we have
interchanged the factors in \eqref{eqnCylindricalSet}.)  As $B$ shrinks and
$s_P$ tends to $o_P$ we obtain a basis of local neighborhoods of $x$.  The
link $|\D_P|$ is itself stratified by its intersections with the strata
$Y_Q>Y_P$; this stratification is
\begin{equation*}
|\D_P| = \coprod_{Q > P} |\D_P^Q|^\circ,
\end{equation*}
where $|\D_P^Q|^\circ$ denotes the open face with vertices $\D_P^Q$.

\section{The reductive Borel-Serre compactification}
\label{sectReductiveBorelSerre}
We outline the construction of $\Xhat$ which was introduced by Zucker
\cite{refnZuckerWarped}; see also \cite{refnGoreskyHarderMacPherson}.

Recall that $e_P$ is a homogeneous space for $\lsp 0P$.  Let $D_P =
N_P\back e_P$ be obtained by collapsing the orbits of the unipotent
radical; it a symmetric space for $L_P$, namely $D_P = L_P/K_PA_P$ where
$K_P = K\cap P$.  The quotient map $e_P \to D_P$ is a principal
$N_P$-bundle.  We can define a geodesic action of $L_P$ on $D$
\cite[\S1.1]{refnSaperLModules}%
\footnote{Precise references to numbered results in
\cite{refnSaperLModules} refer to version 2 as posted on the
\texttt{arXiv}.  The numbering in \cite{refnSaperLModules} may change in
the published version due to revision or a possible splitting of the
paper.}
which extends that of $A_P$; the orbits of the induced geodesic action of
$\lsp0L_P$ on $e_P$ yield natural trivializing sections of $e_P \to D_P$.
We let $\mathscr N_P$ be obtained from $e_P$ be collapsing these sections;
thus there is a natural decomposition
\begin{equation}
e_P=  D_P \times \mathscr N_P.
\label{eqnePDecomposition}
\end{equation}
More concretely, the decomposition $\lsp0 P = \widetilde M_P \ltimes N_P$
(which depends on the choice of basepoint $x_0\in D$) induces a
decomposition
\begin{equation}
e_P \tildearrow D_P\times N_P, \qquad nmK_P \mapsto (mK_P,n),\qquad
  \text{($n\in N_P$, $m\in \widetilde M_P$).}
\label{eqnePBasepointDecomposition}
\end{equation}
This agrees with \eqref{eqnePDecomposition} under the canonical isomorphism
of $N_P$-homogeneous spaces $\mathscr N_P \cong N_P$ which sends the class
of $x_0$ to the identity.  In the example considered in
\S\ref{sectBorelSerre} (see Figure ~\ref{figParabolic}), $D_P$ is the
product of the spaces of positive definite $n_i\times n_i$ symmetric
matrices modulo homothety, for $i=0,\dots,r$.

The action of $\g\in \G_{N_P} \equiv \G \cap N_P \subseteq \G_P$ is solely
on the second factor of \eqref{eqnePDecomposition}; thus we obtain
\begin{equation*}
\G_{N_P}\back e_P = D_P\times \nil_P, \qquad
\text{where $\nil_P \equiv \G_{N_P}\back \mathscr N_P$.}
\end{equation*}
To achieve the full quotient $Y_P = \G_P\back e_P$ we need to now take the
quotient by $\G_{L_P} \equiv \G_P / \G_{N_P}$.  An element $\g\in \G_{L_P}$
acts in the usual way on the first factor $D_P$ but in view of
\eqref{eqnePBasepointDecomposition} it also conjugates the second factor.
Thus the trivial $N_P$-principal bundle $e_P\to D_P$ descends to a flat
bundle of the arithmetic quotients, called the \emph{nilmanifold
fibration}:
\begin{equation*}
\xymatrix{
{Y_P} \ar[d]^{\txt{flat $\nil_P$-bundle}}\\
{X_P \equiv \G_{L_P}\back D_P}
}
\end{equation*}

The \emph{reductive Borel-Serre compactification} $\Xhat$ is obtained by
collapsing the fibers of the nilmanifold fibration on each stratum $Y_P$ of
$\Xbar$.  Thus there is a quotient map $\Xbar \to \Xhat$ which is
stratified with respect to the stratification
\begin{equation*}
\Xhat = \coprod_{P\in \Pl} X_P.
\end{equation*}
The closure of a stratum is $\Xhat_P$, the reductive Borel-Serre
compactification of $X_P$.

In order to understand the structure of $\Xhat$ near a point $x\in X_P$, we
first look at the inverse image of $x$ in $\Xbar$.  This is a full
nilmanifold fiber $\nil_P$ in $Y_P$ which has a basis of
neighborhoods of the form
\begin{equation*}
\widetilde U = B \times \nil_P \times c(|\D_P|),
\end{equation*}
where $B$ maps diffeomorphically to a ball neighborhood of $x$ in $X_P$.
(Again we have interchanged the factors in \eqref{eqnCylindricalSet}.)  The
strata of $\widetilde U$ are
\begin{equation}
\label{eqnTildeNeighborhoodStratification}
\widetilde U\cap Y_Q = \begin{cases}
  B\times \nil_P \times c(|\D_P|)^\circ & Q=G, \\
  B\times \nil_P \times c(|\D_P^Q|)^\circ & P<Q<G, \\
  B\times \nil_P \times c(\emptyset) & Q= P,
		       \end{cases}
\end{equation}
where $c(\emptyset)$ denotes the vertex of the cone and $c(|\D_P^Q|)^\circ
\equiv c(|\D_P^Q|^\circ)\setminus c(\emptyset)$.

A neighborhood $U$ of $x$ is obtained by projecting $\widetilde U$ to
$\Xhat$.  Since the $\nil_Q$-fibers in $Y_Q$  need to be collapsed
for each $Q\ge P$, \eqref{eqnTildeNeighborhoodStratification} is replaced
by
\begin{equation}
\label{eqnNeighborhoodStratification}
U\cap X_Q = \begin{cases}
  B\times \nil_P \times c(|\D_P|)^\circ & Q=G, \\
  B\times \mathscr N_P^{\prime Q} \times c(|\D_P^Q|)^\circ & P<Q<G, \\
  B\times  c(\emptyset) & Q= P,
		       \end{cases}
\end{equation}
where $\mathscr N_P^{\prime Q} \equiv (\G_{N_P}/\G_{N_Q})\back (N_Q\back
\mathscr N_P)$.  In particular, over $X=X_G$ no collapsing occurs, while
over  $X_P$ the entire $\nil_P$ factor is
collapsed to a point.

The link of the stratum $X_P\subseteq \Xhat$ is thus
\begin{equation}
\label{eqnLinkReductiveBorelSerre}
(\nil_P \times |\D_P|)/{\sim} = \coprod_{Q > P}  \mathscr
  N_P^{\prime Q} \times |\D_P^Q|^\circ,
\end{equation}
where $(\G_{N_P}n,a) \sim (\G_{N_P}n',a)$ if $a\in |\D_P^Q|^\circ$ and
$n'=un$ for some $u\in N_Q$.

\section{Satake compactifications}
\label{sectSatake}
The Satake compactification $\Xstar$ depends on an irreducible
representation $(\sigma,U)$ of $G$.  The construction proceeds in several
stages.  Satake \cite{refnSatakeCompact} first constructs a
compactification $\lsb\RR\Dstar$ of $D$ depending on $(\sigma,U)$.  The
compactification $\lsb\RR\Dstar$ is a disjoint union of subsets called
\emph{real boundary components}.  The action of $G$ on $D$ extends to an
action on $\lsb\RR\Dstar$ which permutes the boundary components.  In
\cite{refnSatakeQuotientCompact}, Satake defines $\Dstar\subseteq
\lsb\RR\Dstar$ to be the union of certain of these boundary components,
called \emph{rational boundary components}, and shows that $G(\QQ)$ acts on
$\Dstar$.  He defines the \emph{Satake topology} on $\Dstar$ and sets
$\Xstar \equiv \G\back \Dstar$.  These steps, and certain assumptions that
have to be made on $(\sigma,U)$, are described in more detail below.  In
addition to the above papers of Satake, references are the papers of Borel
\cite{refnBorelEnsemblesFondamentaux}, Casselman
\cite{refnCasselmanGeometricRationality}, and Zucker
\cite{refnZuckerSatakeCompactifications}.

\subsection{}
Assume
\begin{enumerate}
\item $(\sigma,U)$  is nontrivial on each $\RR$-simple component;
\item the maximal $\QQ$-split torus and the maximal $\RR$-split torus in
  the center of $G$ coincide.
\end{enumerate}
(This second assumption is not essential but it simplifies the exposition.)
Give $U$ an \emph{admissible inner product}, that is, one for which
\begin{equation*}
\sigma(g)^* = \sigma(\theta g)^{-1},\qquad\text{for all $g\in G$;}
\end{equation*}
here $\theta$ is the Cartan involution associated to a choice of maximal
compact subgroup $K\subseteq G$.  There is an action of $G$ on the real
vector space $\Herm(U)$ of Hermitian endomorphisms of $U$ given by
\begin{equation*}
S\mapsto\sigma(g)S\sigma(g)^*,\qquad\text{for $S\in \Herm(U)$, $g\in G$.}
\end{equation*}
In particular $g\in K$ acts by conjugation by $\sigma(g)$, while $g\in A_G$
acts by a real scalar (since the center of $G$ acts by a rational character
on an irreducible representation).  Thus if we consider the induced
projective representation, the class of the identity $[I]\in \PP\Herm(U)$
is fixed by $KA_G$ and there is a well-defined map
\begin{equation*}
D = G/KA_G \longrightarrow \PP\Herm(U),\qquad x=gKA_G \mapsto g\cdot [I] =
\sigma(g)\sigma(g)^*.
\end{equation*}
This map is an embedding by the assumptions above and $\lsb\RR\Dstar$ is
defined to be the closure of the image of $D$ in $\PP\Herm(U)$.

An important class of examples are the Hermitian symmetric spaces, in which
$D$ has a $G$-invariant complex structure.  In this case the $\RR$-root
system is always type $C_r$ or $BC_r$ and there is a unique simple
$\RR$-root $\al_r$ at the end of the diagram with a double bond (if the
type is $C_2$ (resp. $BC_2$), $\al_2$ is the long (resp. short) simple
root).  If the highest restricted weight of $(\sigma,U)$ is a multiple of
the fundamental weight dual to $\al_r$, the compactification
$\lsb\RR\Dstar$ is homeomorphic to the closure of $D$ in $\CC^N$ under
Harish-Chandra's embedding as a bounded symmetric domain
\cite{refnHarishChandraRepresentationsSSLieGroupsVI}.  This is called the
\emph{natural compactification}.

\subsection{}
A \emph{real boundary component} is by definition the set of points of
$\lsb\RR\Dstar$ fixed by the unipotent radical of a parabolic
$\RR$-subgroup of $G$; if $P$ is a parabolic $\RR$-subgroup, we let
$D_{P,h}\equiv \Fix(N_P)\cap \lsb\RR\Dstar$ denote the associated real
boundary component.  The compactification $\lsb\RR\Dstar$ is the disjoint
union of its real boundary components.  The action of $P$ on $D_{P,h}$
descends to an action of the Levi quotient $L_P$ and there is an almost
direct product decomposition (direct product of subgroups modulo finite
subgroup) of reductive $\RR$-subgroups
\begin{equation}
\label{eqnLeviFactorization}
L_P = \widetilde L_{P,h} \times L_{P,\l},
\end{equation}
where $L_{P,\l}$ is the subgroup fixing $D_{P,h}$ pointwise and the induced
action of $L_{P,h}\equiv L_P/L_{P,\l}$ on $D_{P,h}$ is that of a real
reductive group on its associated symmetric space.  This decomposition
yields a product decomposition of the associated symmetric spaces,
\begin{equation}
\label{eqnBoundaryStratumFactorization}
D_P =  D_{P,h} \times D_{P,\l}.
\end{equation}

In the case of the natural compactification of a bounded symmetric domain
in $\CC^N$, the real boundary components, which are again Hermitian
symmetric spaces, can also be described as the \emph{holomorphic arc
components}: two points are in the same holomorphic arc component if there
are a finite number of holomorphic maps of $\{z\in\CC\mid |z|<1\}$ into the
closure of $D$ for which the union of the images is connected and contains
the two points.

\subsection{}
We can be more specific about the subgroups $L_{P,\l}$ and $\widetilde
L_{P,h}$, in particular, their root systems.  First some notation. Choose a
maximal $\RR$-split torus $\lsb\RR A$ of $L_P$; via an appropriate lift we
may identify it with a maximal $\RR$-split torus of $G$.  Let $\lsb \RR\D$
be the simple $\RR$-roots of $G$ with respect to a fixed positive system of
$\RR$-roots which contains the roots in $\n_P$.  To any set of
$\RR$-weights $\psi$ is associated a graph with vertex set $\psi$ and an
edge between $\al$, $\b\in\psi$ if $(\al,\b)\neq 0$. (Here we view weights
as belonging to $\lsb\RR\sa^*$ endowed with a Weyl group invariant inner
product.)  The set $\psi$ is \emph{connected} if the associated graph is
connected.

Let $\lsb\RR\u$ be the highest weight of $U$ restricted to $\lsb\RR A$.
The simple $\RR$-roots $\lsb\RR\D^P\subseteq \lsb\RR\D$ of $L_P$ may be
orthogonally decomposed
\begin{equation*}
\lsb\RR\D^P = \kap(\lsb\RR\D^P) \sqcup \z(\lsb\RR\D^P),
\end{equation*}
where $\kap(\lsb\RR\D^P)$ is the largest subset of $\lsb\RR\D^P$ such that
$\kap(\lsb\RR\D^P) \cup \{\lsb\RR\u\}$ is connected and $\z(\lsb\RR\D^P)$
is the complement of $\kap(\lsb\RR\D^P)$ in $\lsb\RR\D^P$.  Alternatively,
a root $\al\in \lsb\RR\D^P$ is in $\kap(\lsb\RR\D^P)$
(resp. $\z(\lsb\RR\D^P)$) if it can (resp. cannot) be connected to
$\lsb\RR\u$ via a chain of weights $\chi_0\equiv \al$, $\chi_1$, \dots,
$\chi_{k-1}$, $\chi_k\equiv\lsb\RR\u$ satisfying $\chi_i\in\lsb\RR\D^P$ for
$i<k$ and $(\chi_{i-1},\chi_i)\neq 0$ for $i\le k$.

We can now describe the identity component of the centralizer $L_{P,\l}$ of
$D_{P,h}$ as the maximal normal connected $\RR$-subgroup of $L_P$ with
simple $\RR$-roots $\z(\lsb\RR\D^P)$.  Similarly the identity component of
$\widetilde L_{P,h}$ is the minimal normal connected $\RR$-subgroup of
$L_P$ with simple $\RR$-roots $\kap(\lsb\RR\D^P)$.

Given a real boundary component $D_{P,h}$, the parabolic subgroup $P$ is
not uniquely determined: any parabolic $\RR$-subgroup $P'\supseteq
P$ such that $\kap(\lsb\RR\D^{P'}) = \kap(\lsb\RR\D^P)$ yields the same
boundary component.  The largest such parabolic $\RR$-subgroup is denoted
$P^\dag$ and we set
\begin{equation*}
\o(\lsb\RR\D^P) \equiv \lsb\RR\D^{P^\dag}.
\end{equation*}
The parabolic $\RR$-subgroups $R$ for which $\lsb\RR\D^R = \o(\psi)$ for
some $\psi\subseteq \lsb\RR\D$ are called \emph{saturated}; they are the
normalizers of the real boundary components.  An example for
$G=\Sympl_{18}(\RR)$ is pictured in Figure ~\ref{figDynkinDiagrams}; here
$D$ is the Siegel generalized upper half-plane (a Hermitian symmetric
space) and $\lsb\RR\Dstar$ is the natural compactification.

\begin{figure}
\entrymodifiers={}
\begin{align*}
&\xymatrix @!0 @M=0pt {{\txt{\llap{$\lsb\RR\D\colon$}}} &
{\circ}\save[]+<0in,-.2in>*{\scriptstyle\al_1}\restore \ar@{-}[r] &
{\circ}\save[]+<0in,-.2in>*{\scriptstyle\al_2}\restore \ar@{-}[r] &
{\circ}\save[]+<0in,-.2in>*{\scriptstyle\al_3}\restore \ar@{-}[r] &
{\circ}\save[]+<0in,-.2in>*{\scriptstyle\al_4}\restore \ar@{-}[r] &
{\circ}\save[]+<0in,-.2in>*{\scriptstyle\al_5}\restore \ar@{-}[r] &
{\circ}\save[]+<0in,-.2in>*{\scriptstyle\al_6}\restore \ar@{-}[r] &
{\circ}\save[]+<0in,-.2in>*{\scriptstyle\al_7}\restore
\ar@{-}@<.15ex>[r] \ar@{-}@<-.15ex>[r] \ar@{}[r]|{<}&
{\circ}\save[]+<0in,-.2in>*{\scriptstyle\al_8}\restore\save []="aaa",\POS "aaa"+/u1ex/+/r4ex/*++++++!L{\scriptstyle \lsb\RR\u}
\ar @{.} "aaa"\restore} \\[2.5ex]
&\xymatrix @!0 @M=0pt {{\txt{\llap{$\lsb\RR\D^P\colon$}}} &
{\circ}\save[]+<0in,-.2in>*{\scriptstyle\al_1}="a"\restore \ar@{-}[r] &
{\circ}\save[]+<0in,-.2in>*{\scriptstyle\al_2}\restore &
{} &
{\circ}\save[]+<0in,-.2in>*{\scriptstyle\al_4}="b"\restore
\save "a";"b" **{}?."a"."b" *!/u1ex/\frm{_\}} *!/u3.25ex/{\scriptstyle
\z(\lsb\RR\D^P)}\restore &
{} &
{\circ}\save[]+<0in,-.2in>*{\scriptstyle\al_6}="a"\restore \ar@{-}[r] &
{\circ}\save[]+<0in,-.2in>*{\scriptstyle\al_7}\restore
\ar@{-}@<.15ex>[r] \ar@{-}@<-.15ex>[r] \ar@{}[r]|{<}&
{\circ}\save[]+<0in,-.2in>*{\scriptstyle\al_8}="b"\restore
\save []="aaa",\POS "aaa"+/u1ex/+/r4ex/*++++++!L{\scriptstyle \lsb\RR\u}
\ar @{.} "aaa"\restore
\save "a";"b" **{}?."a"."b" *!/u1ex/\frm{_\}} *!/u3.25ex/{\scriptstyle
\kap(\lsb\RR\D^P)}\restore} \\[1.5ex]
&\xymatrix @!0 @M=0pt {{\txt{\llap{$\lsb\RR\D^{P^\dag}\colon$}}} &
{\circ}\save[]+<0in,-.2in>*{\scriptstyle\al_1}="a"\restore \ar@{-}[r] &
{\circ}\save[]+<0in,-.2in>*{\scriptstyle\al_2}\restore \ar@{-}[r] &
{\circ}\save[]+<0in,-.2in>*{\scriptstyle\al_3}\restore \ar@{-}[r] &
{\circ}\save[]+<0in,-.2in>*{\scriptstyle\al_4}="b"\restore
\save "a";"b" **{}?."a"."b" *!/u1ex/\frm{_\}} *!/u3.25ex/{\scriptstyle
\z(\lsb\RR\D^{P^\dag})}\restore &
{} &
{\circ}\save[]+<0in,-.2in>*{\scriptstyle\al_6}="a"\restore \ar@{-}[r] &
{\circ}\save[]+<0in,-.2in>*{\scriptstyle\al_7}\restore
\ar@{-}@<.15ex>[r] \ar@{-}@<-.15ex>[r] \ar@{}[r]|{<}&
{\circ}\save[]+<0in,-.2in>*{\scriptstyle\al_8}="b"\restore
\save []="aaa",\POS "aaa"+/u1ex/+/r4ex/*++++++!L{\scriptstyle \lsb\RR\u}
\ar @{.} "aaa"\restore
\save "a";"b" **{}?."a"."b" *!/u1ex/\frm{_\}} *!/u3.25ex/{\scriptstyle
\kap(\lsb\RR\D^{P^\dag})}\restore}
\end{align*}
\caption{An example of the Dynkin diagrams of the $\RR$-root systems of
 $G$, $L_P$, and $L_{P^\dag}$, where $G=\Sympl_{18}(\RR)$ and $(\sigma,U)$
 has highest weight $\u=\b_8$.}
\label{figDynkinDiagrams}
\end{figure}

\subsection{}
Let $A \subseteq \lsb\RR A$ be the identity components of a maximal
$\QQ$-split torus and a maximal $\RR$-split torus respectively.  The
restriction of characters on $\lsb\RR A$ to characters on $A$ yields a map
on simple roots
\begin{equation*}
\r_{\RR/\QQ}\colon \lsb\RR\D \to \D \cup \{0\}
\end{equation*}
and we set 
\begin{equation*}
\epsilon_{\RR/\QQ}(\varUpsilon) \equiv
\r_{\RR/\QQ}^{-1}(\varUpsilon\cup\{0\}),\qquad \text{for any
  $\varUpsilon\subseteq \D$.}
\end{equation*}

We wish now to define a $\G$-invariant union of real boundary components
$\Dstar \subseteq \lsb\RR\Dstar$ so that $\G$ acts discontinuously on
$\Dstar$ and $\G\back \Dstar$ is compact.  It is clear that not all real
boundary components need to be included in $\Dstar$ nor are desirable to be
included.  The following two points motivate the exact definition of
$\Dstar$.

\begin{itemize}
\item \emph{$\Dstar$ should be big enough}: By reduction theory, a
fundamental domain $\Omega$ for the action of $\G$ on $D$ may be formed
from a finite union of \emph{Siegel sets} $\mathfrak S \equiv C \widetilde
A_{P_0}^+b \cdot x_0$.  Here $P_0$ is a minimal parabolic $\QQ$-subgroup,
$C\subseteq P_0$ is compact, $b\in \widetilde A_{P_0}$, and $x_0\in D$ is a
fixed basepoint.  Call a real boundary component \emph{rationally visible}
if it occurs in the closure $\overline\Omega\subseteq \lsb\RR\Dstar$ for
some fundamental domain as above.  For the quotient $\G\back \Dstar$ to be
compact, we wish to include in $\Dstar$ at least all rationally visible
boundary components.  One may check
\cite[Lemma~8.1]{refnCasselmanGeometricRationality} that the rationally
visible boundary components are those whose normalizers $R$ contain a
minimal parabolic $\QQ$-subgroup and satisfy
\begin{equation*}
\lsb\RR\D^R = \o(\epsilon_{\RR/\QQ}(\varUpsilon)) \qquad \text{for some
$\varUpsilon\subseteq \D$.}
\end{equation*}

\item
\emph{$\Dstar$ should not be too big}: On the other hand, call a real
boundary component \emph{geometrically rational} if its normalizer $R$ is
defined over $\QQ$ and if there exists a normal subgroup
$L'_{R,\l}\subseteq L_{R,\l}$ defined over $\QQ$ such that
$L_{R,\l}/L'_{R,\l}$ is compact.  In this case, the decomposition in
\eqref{eqnLeviFactorization} can be modified to be defined over $\QQ$,
essentially by shifting certain compact factors from $L_{R,\l}$ to
$\widetilde L_{R,h}$; we denote these modified terms again by $L_{R,\l}$
and $\widetilde L_{R,h}$.  The geometrical rationality condition on a
boundary component ensures that $\G_{N_R}$, $\G_{L_{R,\l}}\equiv
\G_{L_R}\cap L_{R,\l}$, and $\G_{L_{R,h}} \equiv \G_{L_R}/\G_{L_{R,\l}}$
are all arithmetic subgroups and hence act discontinuously on $N_R$,
$D_{R,\l}$, and $D_{R,h}$.  To ensure that $\G$ acts discontinuously on
$\Dstar$, we do not want to include in $\Dstar$ boundary components that
are not geometrically rational.
\end{itemize}

Any geometrically rational boundary component (in fact, any boundary
component whose normalizer is defined over $\QQ$) is rationally visible.
The two points above suggest that we must assume conversely that
\begin{enumerate}
\setcounter{enumi}{2}
\item  every rationally visible boundary component is geometrically
  rational.
\end{enumerate}
A compactification $\lsb\RR\Dstar$ of $D$ satisfying this assumption is
called \emph{geometrically rational}.  Under this assumption, the
geometrically rational boundary components will simply be called
\emph{rational boundary components} and are just those with normalizers
defined over $\QQ$.  We define
\begin{equation*}
\Dstar \equiv \coprod_{
  \substack{\text{rational}\\\text{boundary components}}}
D_{R,h}.
\end{equation*}
Note that under the assumption of geometric rationality, the Levi quotient
$L_P$ of any parabolic $\QQ$-subgroup $P$ will have a decomposition
(possibly different from \eqref{eqnLeviFactorization} due to shifts of
compact factors) into an almost direct product of reductive $\QQ$-subgroups
\begin{equation}
\label{eqnRationalLeviDecomposition}
L_P = \widetilde L_{P,h} \times L_{P,\l},
\end{equation}
where the $\QQ$-root systems are $\kap(\D^P)$ and $\z(\D^P)$ respectively.
(The notation $\kap(\D^P)$, $\z(\D^P)$, and $\o(\D_P)$ is extended to
subsets of $\D$ by using $\lsb\QQ\u \equiv \u|_A$ instead of $\lsb\RR\u$.)
The normalizers $R$ of the rational boundary components will be the saturated
parabolic $\QQ$-subgroups; these are those with $\o(\D^R) = \D^R$.

The Satake compactification is known to be geometrically rational if
$(\sigma,U)$ is strongly $\QQ$-rational
\cite{refnBorelEnsemblesFondamentaux} or if $\lsb\RR\Dstar$ is the natural
compactification of a Hermitian symmetric space \cite{refnBailyBorel}.  In
the latter case, as we have mentioned, all real boundary components are
again Hermitian symmetric.  Borel has suggested a natural generalization to
consider.  A real boundary component $D_{R,h}$ is \emph{equal-rank} if
$\CCrank L_{R,h} =\rank K_{R,h}$, where $K_{R,h}$ is a maximal compact
subgroup of $L_{R,h}$.  (A Hermitian symmetric space is automatically
equal-rank.)  Call $\lsb\RR\Dstar$ an \emph{real equal-rank Satake
compactification} if $D$ and all its real boundary components are
equal-rank.  Zucker \cite[(A.2)]{refnZuckerLtwoIHTwo} has enumerated the
possible examples when $G$ is $\RR$-simple; there are two new infinite
families, $G=\SO(p,q)$ with $p+q$ odd and $G=\Sympl(p,q)$, with $\lsb\RR\u$
again connected to a simple $\RR$-root at the end of the Dynkin diagram.
In \cite{refnSaperGeometricRationality} we prove that a real equal-rank
Satake compactification is geometrically rational, aside from some
explicitly described $\QQrank$ $1$ and $2$ examples obtained by restriction
of scalars.

\subsection{}
The subspace topology on $\Dstar\subseteq \lsb\RR\Dstar$ is not the right
one to use.  Instead $\Dstar$ is given the \emph{Satake topology}.  To
describe this topology, recall that $\overline \Omega$ is the closure in
$\lsb\RR\Dstar$ of a fundamental domain $\Omega$ constructed from Siegel
sets.  We have $\Dstar = \G \cdot \overline \Omega$.  Then a basis of local
neighborhoods of $x\in\Dstar$ in the Satake topology is given by subsets
$U$ such that
\begin{enumerate}
\item If $x\in \g \overline \Omega$, then $U\cap \g\overline\Omega$ is open
  in the subspace topology on $\g \overline \Omega\subseteq \lsb\RR\Dstar$.
\item $U$ is $\G_x$-invariant, where $\G_x=\{\g\in\G\mid \g\cdot x = x\}$.
\end{enumerate}
The \emph{Satake compactification} $\Xstar$ is defined by
\begin{equation*}
\Xstar \equiv \G\back \Dstar;
\end{equation*}
it has a stratification
\begin{equation*}
\Xstar = \coprod_{R\in \Pl^*} X_{R,h}
\end{equation*}
where $\Pl^*\subseteq \Pl$ is the set of $\G$-equivalence classes of
saturated parabolic $\QQ$-subgroups and
\begin{equation*}
X_{R,h} \equiv \G_{L_{R,h}} \back D_{R,h}.
\end{equation*}

The Satake compactification $\Xstar$ associated to the natural
compactification of a Hermitian symmetric space is known as the
\emph{Baily-Borel-Satake compactification}; Baily and Borel
\cite{refnBailyBorel} show that it admits the structure of a normal
projective variety defined over a number field.  In this case $\Pl^*$
simply consists of $\G$-equivalence classes of maximal parabolic
$\QQ$-subgroups.  More generally, the Satake compactification $\Xstar$
associated to a geometrically rational real equal-rank Satake
compactification $\Dstar$ is again called a \emph{real equal-rank Satake
compactification}.

\subsection{}
\label{ssectProjectSatakeToRBS}
There is an alternative description due to Zucker
\cite{refnZuckerSatakeCompactifications} of the topology on $\Xstar$ that
is quite useful.  For $P\in\Pl$, we have the projection
\begin{equation*}
D_P =  D_{P,h} \times D_{P,\l} \longrightarrow D_{P,h} = D_{P^\dag,h}
\end{equation*}
onto the first factor of 
\eqref{eqnBoundaryStratumFactorization}.  This induces a map on the
arithmetic quotients
\begin{equation*}
\pi_P\colon X_P \longrightarrow X_{P^\dag,h}.
\end{equation*}
In general this is a flat bundle with typical fiber $X_{P,\l}\equiv
\G_{L_{P,\l}}\back D_{P,\l}$, however if we replace $\G_{L_P}$ by a certain
finite index subgroup the corresponding bundle becomes trivial.

The maps $\pi_P$ for the various $P\in\Pl$ fit together to define a
surjective map
\begin{equation*}
\pi\colon \Xhat \to \Xstar
\end{equation*}
and Zucker proves that the topology on $\Xstar$ is the quotient topology
induced by $\pi$.  (See Figure ~\ref{figProjection} in
\S\ref{sectFunctorialityMicroSupport} for a picture of the map $\pi$.)  One
can also define a continuous map $\widehat{D} \equiv \coprod_P D_P \to
\Dstar$ before taking the arithmetic quotient, however this map is not
necessarily open.

\section{The Rapoport/Goresky-MacPherson conjecture}
The preceding sections, and especially the table at the end of
\S\ref{ssectGeneralCase}, explicitly illustrate how the links of strata
become more complicated as we pass from $\Xbar$ to $\Xhat$ to $\Xstar$.
The intersection cohomology sheaf will likewise become more complicated.  A
hope therefore is that one might be able to transfer the study of the
cohomology of various sheaves, in particular the intersection cohomology
sheaf, from $\Xstar$ to the simpler space $\Xhat$.  (From now on in this
paper, a ``sheaf'' means a complex of sheaves, viewed as an element of the
constructible derived category.)  In particular, the replacement of
$\Xstar$ by $\Xhat$ in Theorem ~\ref{thmZuckersConjecture}, Zucker's
conjecture, would be of considerable use in applications.

For example, if $\Xstar$ is the Baily-Borel-Satake compactification of a
Hermitian locally symmetric space, part of Langlands's program involves
relating the Hasse-Weil zeta function of $\Xstar$ to automorphic
$L$-functions; see for example \cite{refnLanglandsRamakrishnan} where
Picard modular surfaces are considered.  One approach to establishing this
relationship would be to compare the Arthur-Selberg trace formula
\cite{refnArthurLTwoLefschetzNumbers} for a Hecke operator on
$H_{(2)}(X;\EE)$ with a topological Lefschetz fixed point formula for
$I_pH(\Xstar;\EE)$ such as Goresky and MacPherson give in
\cite{refnGoreskyMacPersonLefschetzFixedPointFormula}; this is reasonable
since Zucker's conjecture shows that the groups are isomorphic.  However
Arthur's formula involves a sum over $P\in\Pl$ whereas a fixed point
formula on $\Xstar$ will involve a sum over $R\in\Pl^*$; since $\Pl$
indexes the strata of $\Xhat$, a fixed point formula on $\Xhat$ might be
more amenable to a comparison with Arthur's formula.

These considerations were one motivation for the following conjecture.  It
was first posed by Rapoport in 1986 within a letter to Borel
\cite{refnRapoportLetterBorel}; the conjecture was later independently
posed by Goresky and MacPherson and described in an unpublished preprint
\cite{refnGoreskyMacPhersonWeighted}.

\begin{conj}[Rapoport/Goresky-MacPherson]
\label{conjRapoportGoreskyMacPherson}
Let $X$ be a Hermitian locally symmetric space and let $p$ be a
middle-perversity.  There is a natural quasi-isomor\-phism
$\IpC(\Xstar;\EE)\cong \pi_*\IpC(\Xhat;\EE)$ and hence
$I_pH(\Xstar;\EE)\cong I_pH(\Xhat;\EE)$.
\end{conj}

Previously Zucker had noticed that the conjecture held for
$G=\Sympl_4(\RR)$ with $\EE=\CC$.  In addition Goresky and MacPherson
\cite{refnGoreskyMacPhersonWeighted} announced the conjecture held for
$G=\Sympl_4(\RR)$, $\Sympl_6(\RR)$, and (for $\EE=\CC$) $\Sympl_8(\RR)$.
In the case where $\QQrank G=1$ a proof was given by the author and Stern
in an appendix to a paper of Rapoport's \cite{refnRapoport}.  The full
conjecture (in fact a generalization to the equal-rank setting) was finally
proved by the author in 2001 \cite{refnSaperLModules} (see Theorem
~\ref{thmRapoportConjecture} below).

\section{Weighted cohomology}
With similar motivation, but in a different direction from the above
conjecture, Goresky, Harder, and MacPherson have defined the weighted
cohomology $W^\eta H(\Xhat;\EE)$ and have proven
\cite{refnGoreskyHarderMacPherson}:
 
\begin{thm}[Goresky-Harder-MacPherson]
\label{thmGoreskyHarderMacPherson}
Let $X$ be a Hermitian locally symmetric space, let $p$ be a
middle-perversity, and let $\eta$ be a middle-weight profile.  There is a
natural quasi-isomorphism $\IpC(\Xstar;\EE) \cong \pi_*\WnC(\Xhat;\EE)$ and
hence an isomorphism $I_pH(\Xstar;\EE)\cong W^\eta H(\Xhat;\EE)$.
\end{thm}

Weighted cohomology depends on a locally constant sheaf $\EE$ on $X$
arising from a representation $E$ of $G$, as well as on a \emph{weight
profile} $\eta \in \sa_{P_0}^{G*}$ where $P_0$ is a minimal parabolic
$\QQ$-subgroup.  \pagebreak The \emph{middle-weight profiles} are $\v = -\r$ and $\u =
-\r+\epsilon\r$,  where $\epsilon>0$ is very small and $\r$ is one-half the
sum of the positive $\QQ$-roots of $G$ (counted with multiplicity).  Like
$L^2$-cohomology and intersection cohomology, there is a sheaf
$\WnC(\Xhat;\EE)$ whose hypercohomology is weighted cohomology.

We will give a local calculation of $W^\eta H(\Xhat;\EE)$ analogous to the
one for intersection cohomology in \eqref{eqnLocalCalculation} of
\S\ref{ssectSimplicialIC}.  Recall that the smooth part of the link of a
stratum $X_P\subseteq \Xhat$ is $\nil_P\times|\D_P|^\circ$ (see
\eqref{eqnLinkReductiveBorelSerre} in \S\ref{sectReductiveBorelSerre}).
Van Est's theorem \cite{refnvanEst} shows that the cohomology
$H(\nil_P\times|\D_P|^\circ;\EE)\cong H(\nil_P;\EE)$ is naturally an
$A_P^G$-module (in fact, an $L_P$-module); we will explain this in more
detail in \S\ref{ssectLModulesOverview}.  Thus there is a decomposition
$H(\nil_P;\EE) = \bigoplus_\chi H(\nil_P;\EE)_\chi$ into $A_P^G$-weight
spaces.  If $U$ is a small neighborhood of $x\in X_P$, the local weighted
cohomology $W^\eta H^j(U;\EE)$ similarly has a decomposition into
$A_P^G$-weight spaces and the local characterization is
\begin{equation}
\label{eqnLocalWeightedCohomology}
W^\eta H(U;\EE)_\chi \cong
\begin{cases}
  H(\nil_P;\EE)_\chi & \text{for $\chi \ge \eta_P$,} \\
  0                  & \text{for $\chi \not\ge \eta_P$;}
\end{cases}
\end{equation}
here $\eta_P$ is the restriction of $\eta$ to $\sa_P^G$ and the inequality
$\chi \ge \eta_P$ means that $\chi-\eta_P$ is a nonnegative linear
combination of elements of $\D_P$.

\begin{rems}
(i) Despite the similarity of \eqref{eqnLocalWeightedCohomology} with the
local calculation of intersection cohomology \eqref{eqnLocalCalculation},
an important difference is that \eqref{eqnLocalWeightedCohomology} is not
inductive.  Thus in general local calculations involving weighted
cohomology are easier than the corresponding local calculations for
intersection cohomology.

(ii) There is a natural projection $U\cap X\to \nil_P$.  If we represent a
class in $H(\nil_P;\EE)_\chi\cong H(U\cap X;\EE)_\chi$ by a differential
form $\phi$ which is the pullback of a form on $\nil_P$, then the condition
$\chi \ge \u_P$ is precisely the condition that $\phi$ is $L^2$.  Thus
middle-weight profile weighted cohomology may be viewed as an algebraic
analogue of $L^2$-cohomology.
\end{rems}

Goresky and MacPherson \cite{refnGoreskyMacPhersonTopologicalTraceFormula}
have applied their fixed point formula
\cite{refnGoreskyMacPersonLefschetzFixedPointFormula} to calculate the
Lefschetz number of a Hecke correspondence on $W^\eta H(\Xhat;\EE)$.  They
prove in collaboration with Kottwitz
\cite{refnGoreskyKottwitzMacPhersonDiscreteSeries} that the result agrees
with Arthur's trace formula \cite{refnArthurLTwoLefschetzNumbers} on
$H_{(2)}(X;\EE)$, in line with the motivation from the previous section.

A proof of Theorem ~\ref{thmGoreskyHarderMacPherson} (different from the
original one) can easily be sketched using Zucker's conjecture and later
work of Nair \cite{refnNair}.  First of all, Zucker's conjecture yields a
quasi-isomorphism
\begin{equation}
\label{eqnZuckerConjecture}
\IpC(\Xstar;\EE) \cong \mathcal L_{(2)}(\Xstar;\EE).
\end{equation}
However it is fairly immediate that
\begin{equation}
\label{eqnLtwoInvariant}
\mathcal L_{(2)}(\Xstar;\EE) \cong \pi_* \mathcal L_{(2)}(\Xhat;\EE).
\end{equation}
Finally Nair's result is that there is a quasi-isomorphism
\begin{equation*}
\mathcal L_{(2)}(\Xhat;\EE)\cong \WnC(\Xhat;\EE)
\end{equation*}
as suggested by Remark (ii) above and hence
\begin{equation}
\label{eqnPushForwardNair}
\pi_*\mathcal L_{(2)}(\Xhat;\EE)\cong \pi_*\WnC(\Xhat;\EE).
\end{equation}
Equations \eqref{eqnZuckerConjecture}--\eqref{eqnPushForwardNair} yield
Theorem ~\ref{thmGoreskyHarderMacPherson}.

\section{Main theorem}
Despite the success of weighted cohomology, Conjecture
~\ref{conjRapoportGoreskyMacPherson} is still of interest.  The short
argument given above for Theorem ~\ref{thmGoreskyHarderMacPherson} does not
apply to prove the conjecture since in general the analogue of Nair's
theorem does not hold: $\IpC(\Xhat;\EE)$ is not quasi-isomorphic to
$\mathcal L_{(2)}(\Xhat;\EE)$ or $\WnC(\Xhat;\EE)$ on $\Xhat$.  The only
hope for the conjecture is that a quasi-isomorphism exists downstairs on
$\Xstar$ and indeed this is the case:
\begin{thm}
\label{thmRapoportConjecture}
Let $\Xstar$ be a real equal-rank Satake compactification and let $p=m$ or
$n$ be a middle-perversity.  There is a natural quasi-isomorphism
$\IpC(\Xstar;\EE) \cong \pi_*\IpC(\Xhat;\EE)$ and hence an isomorphism
$I_pH(\Xstar;\EE)\cong I_pH(\Xhat;\EE)$.
\end{thm}
This result is proved in \cite{refnSaperLModules} as a consequence of a
number of results in the theory of $\L$-modules which we will outline in
the remaining sections.

Note that the condition on $\Xstar$, that all real boundary components of
$\lsb\RR\Dstar$ are equal-rank, always holds if $D$ is Hermitian and
$\Xstar$ is the Baily-Borel-Satake compactification.  Thus the theorem
implies the conjecture of Rapoport and Goresky-MacPherson.  Note that
although the equal-rank condition implies that the strata of $\Xstar$ have
even dimension, the strata of $\Xhat$, on the other hand, may have odd
codimension, even in the Hermitian case.

\section{$\L$-modules}
In the course of work on the Rapoport/Goresky-MacPherson conjecture, it was
realized that many of the techniques could be applied to study other
sheaves.  Furthermore, many parts of the proof, including a crucial and
difficult local property of $\IpC(\Xhat;\EE)$, apparently held even when
$D$ was non-Hermitian or non-equal-rank.  So it seemed desirable in view of
other applications to separate the components of the proof and to work in a
more general context.  This was one motivation for us to develop the theory
of $\L$-modules \cite{refnSaperLModules}.

\subsection{Overview}
\label{ssectLModulesOverview}
First recall the notation.  The reductive Borel-Serre compactification has
a stratification
\begin{equation*}
\Xhat = \coprod_{P\in\Pl} X_P,
\end{equation*}
where $\Pl$ is the partially ordered set of $\G$-conjugacy classes of
parabolic $\QQ$-subgroups of $G$.  For a parabolic $\QQ$-subgroup $P$ of
$G$, let $N_P$ denote its unipotent radical and let $L_P=P/N_P$ denote its
Levi quotient.  The reductive Borel-Serre boundary stratum $X_P$ is the
locally symmetric space associated to $L_P$, namely
\begin{equation*}
X_P =  \G_{L_P}\back  L_P / K_P A_P =  \G_{L_P}\back D_P,
\end{equation*}
\enlargethispage*{\baselineskip}
where $A_P$ is the identity component of a maximal $\QQ$-split torus in the
center of $L_P$, $\G_P = \G\cap P$, $\G_{N_P}= \G\cap N_P$ and $\G_{L_P} =
\G_P/\G_{N_P}$.  The link of $X_P\subseteq \Xhat$ has as its smooth open
stratum $\nil_P \times |\D_P|^\circ$, where $\nil_P$ is a compact
nilmanifold which is canonically diffeomorphic to $\G_{N_P}\back N_P$ once
a basepoint is chosen.  The rest of the link is described in
\S\ref{sectReductiveBorelSerre}.  In general, if $P\subseteq Q$, we extend
all this notation to the parabolic $\QQ$-subgroup $P/N_Q$ of $L_Q$ simply
by adding a superscript $Q$; thus we have $N_P^Q = N_P/N_Q$, etc.

An $\L$-module
\begin{equation*}
\M = {(E_P,f_{PQ})}_{P \le Q\in \Pl}
\end{equation*}
on $\Xhat$ is a combinatorial analogue of a constructible complex of
sheaves.  (The precise definition will be given below.)  Like the
constructible derived category of sheaves on $\Xhat$, the category of
$\L$-modules has the functors $k_*$, $k_!$, $k^*$, and $k^!$ as well as
truncation functors.  Furthermore there is a realization functor $\M
\mapsto \Sheaf(\M)$ from $\L$-modules to the constructible derived category
of sheaves, commuting with these functors.  Many familiar sheaves can be
realized as $\Sheaf(\M)$ (or ``lifted to an $\L$-module''), for example,
the intersection cohomology sheaf $\IpC(\Xhat;\EE)$, the weighted
cohomology sheaf $\WnC(\Xhat;\EE)$, or the sheaf $\i_{G*}\EE$ whose
hypercohomology is the ordinary cohomology $H(X;\EE)$.  The global
cohomology of an $\L$-module is defined to be the hypercohomology of its
realization:
\begin{equation*}
H(\Xhat;\M) \equiv H(\Xhat;\Sheaf(\M)).
\end{equation*}
In fact the realization functor yields an explicit complex of fine sheaves
so that the cohomology may be computed from the complex of their global
sections; we will see this complex is easy to work with.

Let $\i_P\colon X_P\hookrightarrow \Xhat$ and $\ihat_P\colon
\Xhat_P\hookrightarrow \Xhat$ denote the inclusion of a stratum and of its
closure.  Recall that for a sheaf $\Sheaf$ to be in the constructible
derived category, the local cohomology sheaf $H(\i_P^*\Sheaf)$ along $X_P$
must be a finitely generated locally constant sheaf for all $P$.  It is
equivalent%
\footnote{This follows by the methods of
  \cite[V.3]{refnBorelIntersectionCohomology} though it is not explicitly
  stated there.}
to require that the local cohomology $H(\i_P^!\Sheaf)$ supported on $X_P$
is a finitely generated locally constant sheaf for all $P$; we prefer this
point of view for reasons that will become clear later.  Consequently for a
sheaf, $H(\i_P^!\Sheaf)$ may be identified with a graded representation of
the arithmetic subgroup $\G_{L_P}\subseteq L_P$.  An important difference
between a sheaf and an $\L$-module is that for an $\L$-module $\M$,
$H(\i_P^!\M) = H(E_P,f_{PP})$ is actually a graded representation of the
full reductive group $L_P$.

The reason this is possible, namely that the usual functorial operations
applied to representations of Levi quotients yield representations of Levi
quotients (as opposed to merely representations of fundamental groups) is
van Est's theorem \cite{refnvanEst}.  We will illustrate this with the
simplest example of such a functorial operation.  Suppose we start with a
representation $E$ of $\G \subseteq G$ and hence a locally constant sheaf
$\EE$ on $X=X_G$.  Pushforward $\EE$ to $\i_{G*}\EE$ on $\Xhat$, and then
restrict it to $\i_P^*\i_{G*}\EE$ on $X_P$.  The local cohomology
$H(\i_P^*\i_{G*}\EE)$ is simply the cohomology of the smooth open stratum
of the link $\nil_P\times |\D_P|^\circ$, which is homotopic to $\nil_P$.
Thus the local cohomology is $H(\nil_P;\EE)$ and since this is the
cohomology of the fiber of the flat bundle $Y_P\to X_P$, it has an action
on it of $\G_{L_P}\subseteq L_P$.  However van Est's theorem says that if
$E$ actually arose from a representation of $G$, there is an isomorphism of
$\G_{L_P}$-modules
\begin{equation}
\label{eqnVanEst}
H(\nil_P;\EE) \cong H(\n_P;E).
\end{equation}
The $\G_{L_P}$-action on the complex $\bigwedge \n_P^*\otimes E$ is
obtained by lifting $\G_{L_P}$ to $G$; the induced action on the cohomology
$H(\n_P;E)$ is independent of the lift.  This $\G_{L_P}$-action, however,
extends to an $L_P$-action defined by the same procedure.

The idea of representing a sheaf by combinatorial and linear algebraic data
is not new.  In other contexts there is the theory of $S$-constructible
sheaves \cite{refnKashiwaraRiemannHilbert}, \cite[Exercise
~VIII.1]{refnKashiwaraSchapira} and MacPherson's theory of cellular sheaves
\cite{refnShepardThesis}, \cite{refnVybornovConstructible},
\cite{refnVybornovSheaves}.  These two examples are based on a
decomposition into contractible cells, unlike the decomposition of $\Xhat$
into strata.  For examples with non-contractible cells, there is the theory
of sheaves on fans \cite{refnBarthelEtAlEquivariant},
\cite{refnBarthelEtAlCombinatorial},
\cite{refnBresslerLuntsToricVarieties},
\cite{refnBresslerLuntsIntersectionCohomology},
\cite{refnMcConnellIntersectionCohomology} and Braden and MacPherson's
moment graphs \cite{refnBradenMacPhersonMomentGraphs}.  There is also
forthcoming work of Braden \cite{refnBradenKoszulDuality} on
Koszul duality, in which he constructs a combinatorial category of mixed
sheaves on affine toric varieties.  The forthcoming work of Pink and
Wildeshaus on automorphic sheaves \cite{refnPinkJAMITalk} applies in a
context closest to that of $\L$-modules.  The basic linear algebraic data
in all these theories is the local cohomology along a cell or stratum,
$\i_P^*\Sheaf$.  An $\L$-module uses instead the dual notion of the local
cohomology supported on a stratum $E_P = \i_P^!\Sheaf$ as the building
block.  This approach turns out to make the construction of a useful
realization technically easier.

\subsection{\boldmath $\L$-modules}
An \emph{$\L$-module} $\M = (E_\cdot,f_{\cdot \cdot})$ on $\Xhat$ consists
of
\begin{itemize}
\item for all $P\in\Pl$, a graded algebraic representation $E_P$ of $L_P$;
\item  for all $P \le  Q$, morphisms
  $f_{PQ}\colon H(\n_P^Q;E_Q)\xrightarrow{[1]} E_P$.
\end{itemize}
(The superscript $[1]$ indicates that $f_{PQ}$ is graded of degree $1$,
that is, it maps $H^i(\n_P^Q;E_Q)$ into $E_P^{i+1}$.)  This data must
satisfy the condition that for all $P \le R$,
\begin{equation}
\label{eqnLmoduleDifferential}
\sum_{P\le Q\le R} f_{PQ}\circ H(\n_P^Q;f_{QR}) = 0.
\end{equation}
In particular, $(E_P,f_{PP})$ is a complex of $L_P$-modules.  Here
$H(\n_P^Q;f_{QR})$ is viewed as a map $H(\n_P^R;E_R) \to H(\n_P^Q;E_Q)$ via
the canonical isomorphism
\begin{equation}
\label{eqnNilpotentCanonicalIsomorphism}
H(\n_P^R;E_R) \cong H(\n_P^Q;H(\n_Q^R;E_R)).
\end{equation}

This is a good place to recall the standard notation regarding shifts of
degree: if $C$ is a complex, then $C[k]$ denotes the shifted complex
defined by
\begin{equation*}
C[k]^i \equiv C^{i+k}, \qquad d_{C[k]} \equiv (-1)^k d_C.
\end{equation*}
This applies to graded objects (viewed as complexes with zero differential)
and ordinary objects (viewed as complexes with nonzero entries only in
degree $0$) as well.  Thus if $P<Q$ and there are no elements of $\Pl$
strictly between $P$ and $Q$, the condition \eqref{eqnLmoduleDifferential}
shows that $f_{PQ}$ is a morphism of complexes,
\begin{equation*}
(H(\n_P^Q;E_Q),H(\n_P^Q;f_{QQ})) \longrightarrow (E_P,f_{PP})[1].
\end{equation*}

A \emph{morphism} $\phi =(\phi_{\cdot\cdot})$ of $\L$-modules $\M \to \M'$
consists of graded maps
\begin{equation*}
\phi_{PQ}\colon H(\n_P^Q;E_Q) \longrightarrow E'_P \qquad\text{for all
  $P\le Q$}
\end{equation*}
satisfying
\begin{equation*}
  \sum_{P\le Q\le R} \phi_{PQ} \circ H(\n_P^Q;f_{QR}) =  \sum_{P\le Q\le R}
  f'_{PQ} \circ H(\n_P^Q;\phi_{QR}) \qquad\text{for all
  $P\le R$.}
\end{equation*}

A subspace $W\subseteq \Xhat$ is \emph{admissible} if it is a locally
closed union of strata; set $\Pl(W) = \{\, P\in \Pl\mid X_P\subseteq
W\,\}$.  In this case we define an $\L$-module on $W$ as above but with
$P$, $Q$, and $R$ restricted to belong to $\Pl(W)$.

\begin{figure}
\begin{equation*}
\xymatrix @R=.75pc @C=1.5pc @!0 { {} & {} & {P} \ar@{-}[dddd] \\ {} & {} & {}
\\ {Q_1} \ar@{-}[rruu] \ar@{-}[dddd] \\ {} & {} & {} \\ {} & {} & {Q_2} \\
{} & {} & {} \\ {G} \ar@{-}[rruu] }
\end{equation*}
\caption{The (inverted) Hasse diagram for  $\Pl= \{P<Q_1,Q_2< G\}$}
\label{figHasseDiagram}
\end{figure}
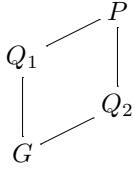

\begin{figure}
\begin{equation*}
\xymatrix 
@R=1.25cm @C=3cm @!0
{
{} & {} & {E_P} \ar@(r,u)[]_{f_{PP}} \\
{} & {H(\n_P^{Q_1};E_{Q_1})} \ar[ur]^-{f_{PQ_1}} & {} \\
{E_{Q_1}} \ar@(r,u)[]_{f_{Q_1Q_1}} & {} & {H(\n_P^{Q_2};E_{Q_2})}
\ar[uu]^-*-<1\jot>{\labelstyle f_{PQ_2}}\\
{} & {H(\n_P;E_G)} \ar[uuur]^-*-<1\jot>{\labelstyle f_{PG}}& {} \\
{H(\n_{Q_1};E_G)} \ar[uu]^{f_{Q_1G}} & {} & {E_{Q_2}}
\ar@(r,u)[]_{f_{Q_2Q_2}} \\
{} & {H(\n_{Q_2};E_G)}  \ar[ru]^-{f_{Q_2G}} & {} \\
{E_G} \ar@(r,u)[]_{f_{GG}}
}
\end{equation*}
\caption{The data of a typical $\L$-module $\M$ for the $\Pl$ from Figure
  ~\ref{figHasseDiagram}.  (All objects are graded and all morphisms are
  degree $1$.)}
\label{figExampleLModule}
\end{figure}
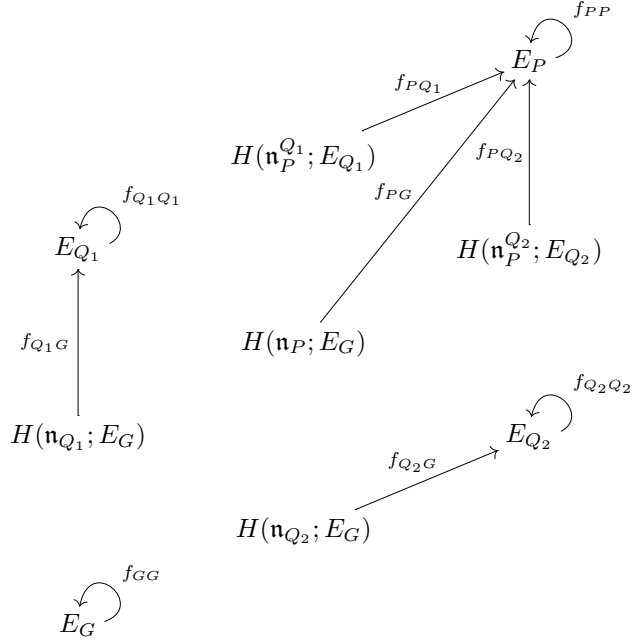

To illustrate the above definition, suppose that $\QQrank G = 2$ and that
there is only one $\G$-conjugacy class of parabolic $\QQ$-subgroups of each
type; an example where this occurs is $G=\Sympl_4(\RR)$ and
$\G=\Sympl_4(\ZZ)$.  Then $\Pl= \{P,Q_1,Q_2, G\}$ with the covering
relations $P < Q_1, Q_2 < G$.  In Figure ~\ref{figHasseDiagram} we
represent $\Pl$ by its Hasse diagram (which we invert so as to match with
Figure ~\ref{figSplitComponent}).
The data of an $\L$-module $\M$ for this example (but not condition
\eqref{eqnLmoduleDifferential}) is displayed in Figure
~\ref{figExampleLModule}.

\subsection{Functors on \boldmath$\L$-modules}
\label{ssectFunctorsLModules}
We will only define the usual functors in the cases we need here.  If
$k\colon Z\hookrightarrow W$ is an inclusion of admissible spaces, then
$k^!\M$ for an $\L$-module $\M$ on $W$ is defined by restricting the
subscripts of $E_P$ and $f_{PQ}$ to belong to $\Pl(Z)$.  In particular,
\begin{equation*}
\i_P^!\M = (E_P,f_{PP}).
\end{equation*}
for all $P\in \Pl(W)$.  As we have noted, this is a complex of
$L_P$-modules; its cohomology $H(\i_P^!\M)$ is the \emph{local cohomology
supported on $X_P$}.  The less trivial example of $\ihat_{Q_1}^!\M$ for our
example is represented in Figure ~\ref{figIHatShriekExample}.
\begin{figure}
\begin{equation*}
\xymatrix 
@R=1.25cm @C=3cm @!0
{
{} & {} & {E_P} \ar@(r,u)[]_{f_{PP}} \\
{} & {H(\n_P^{Q_1};E_{Q_1})} \ar[ur]^-{f_{PQ_1}} & {} \\
{E_{Q_1}} \ar@(r,u)[]_*-<1\jot>{\labelstyle f_{Q_1Q_1}}
}
\end{equation*}
\caption{The $\L$-module $\ihat_{Q_1}^!\M$ ($\M$ as in Figure
~\ref{figExampleLModule})}
\label{figIHatShriekExample}
\end{figure}

On the other hand, for $P\in \Pl(W)$ we define
\begin{equation*}
\i_P^*\M \equiv \left( \bigoplus_{P\le Q} H(\n_P^Q;E_Q) , \sum _{P\le Q\le
  R} H(\n_Q^R; f_{QR})\right).
\end{equation*}
Again by \eqref{eqnLmoduleDifferential} this is a complex of $L_P$-modules;
its cohomology $H(\i_P^*\M)$ is the \emph{local cohomology at $X_P$} (or
\emph{at $P$}).  The complexes $\i_{Q_1}^*\M$ and $\i_P^*\M$ for our usual
example are pictured in Figure~\ref{figUpperStarExample}.
\newbox\subfigbox
\newdimen\subfigwd
\begin{figure}
\begin{center}
  \mbox{\hfil%
    \vbox{%
      \setbox\subfigbox=\hbox{\labelmargin-{1pt}%
	$\vcenter{\xymatrix 
	  @R=1.25cm @C=3cm @!0
	  {
	    {E_{Q_1}} \ar@(r,u)[]_{f_{Q_1Q_1}} \\
	    {} \\
	    {H(\n_{Q_1};E_G)} \ar[uu]^{f_{Q_1G}}
	    \ar@(r,u)[]!CR;[]!/r1ex/_{H(\n_{Q_1};f_{GG})}
	}}$}%
      \dimen\subfigwd=\wd\subfigbox%
      \box\subfigbox%
      \vskip\baselineskip%
      \hbox to \dimen\subfigwd{\hfil\footnotesize(a) $\i_{Q_1}^*\M$\hfil}%
    }%
    \qquad%
    \vbox{%
      \setbox\subfigbox=\hbox{
	$\vcenter{\xymatrix 
	  @R=1.25cm @C=3cm @!0
	  {
	    {} & {E_P} \ar@(r,u)[]_{f_{PP}} \\
	    {H(\n_P^{Q_1};E_{Q_1})} \ar@(u,l)[];[]!CL_-*+h{\labelstyle H(\n_P^{Q_1};f_{Q_1Q_1})}
	    \ar[ur]^-{f_{PQ_1}} & {} \\
	       {} & {H(\n_P^{Q_2};E_{Q_2})}
	       \ar@(r,u)[]!CR;[]!/r1ex/_{H(\n_P^{Q_2};f_{Q_2Q_2})}
      \ar[uu]^*-<1\jot>{\labelstyle f_{PQ_2}}\\
		   {H(\n_P;E_G)} \ar[uu]^{H(\n_P^{Q_1};f_{Q_1G})}
      \ar[uuur]^-*-<1\jot>{\labelstyle f_{PG}}
		   \ar[ur]_{H(\n_P^{Q_2};f_{Q_2G})} \ar@(l,d)[]!CL;[]_{H(\n_P;f_{GG})} & {}
	}}$
      }%
      \dimen\subfigwd=\wd\subfigbox%
      \box\subfigbox%
      \vskip\baselineskip%
      \hbox to \dimen\subfigwd{\footnotesize\hfil(b) $\i_P^*\M$\hfil}%
    }%
  \hfil}%
\end{center}
\caption{The complexes computing local cohomology of $\M$ at $Q_1$ and at
$P$ ($\M$ as in Figure ~\ref{figExampleLModule})}
\label{figUpperStarExample}
\end{figure}
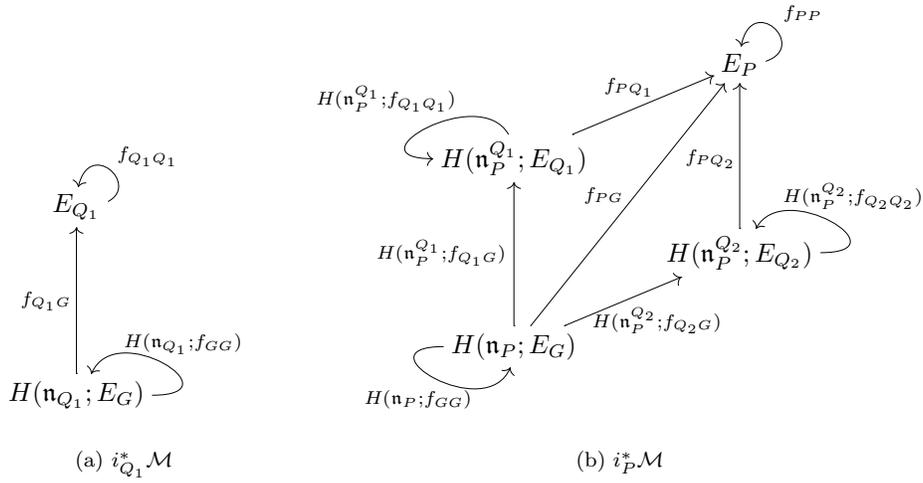

If $R\in\Pl(W)$, the pullback $\ihat_R^*\M = (E'_\cdot,f'_{\cdot\cdot})$ to
an entire closed face $\Xhat_R\cap W$ is defined by
\begin{equation*}
E'_P \equiv \bigoplus_{\substack{Q\ge P\\Q\cap R=P}} H(\n_P^Q;E_Q) \qquad
\text{for all $P\le R$},
\end{equation*}
with the obvious induced morphisms for $f'_{\cdot\cdot}$.

Finally $k_*\M'$ for an $\L$-module $\M'$ on $Z$ is defined by extending
$E_P$ and $f_{PQ}$ to be zero if $P$ or $Q$ belongs to $\Pl(W)\setminus
\Pl(Z)$.

In the cases we have defined, the usual adjoint relations hold among these
functors.  (For some of the other cases, it may be necessary to pass to the
``homotopy category of $\L$-modules'' which we do not consider here.)  In
particular, $\i_P^*$ is a left adjoint to $\i_{P*}$ on $\L$-modules and so
there is a natural \emph{adjunction morphism}
\begin{equation*}
\M \longrightarrow  \i_{P*}\i_P^*\M
\end{equation*}
which is characterized by the condition that $\i_P^*\M \longrightarrow
\i_P^*\i_{P*}\i_P^*\M = \i_P^* \M$ is the identity map.

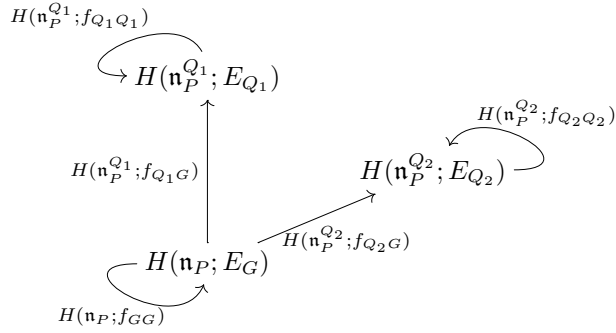
\begin{figure}
\begin{equation*}
\xymatrix 
@R=1.25cm @C=3cm @!0
{
{H(\n_P^{Q_1};E_{Q_1})} \ar@(u,l)[];[]!CL_-*+{\labelstyle H(\n_P^{Q_1};f_{Q_1Q_1})} \\
{} & {H(\n_P^{Q_2};E_{Q_2})} 
\ar@(r,u)[]!CR;[]!/r1ex/_{H(\n_P^{Q_2};f_{Q_2Q_2})} \\
{H(\n_P;E_G)} \ar[uu]^{H(\n_P^{Q_1};f_{Q_1G})}
\ar[ur]_-{H(\n_P^{Q_2};f_{Q_2G})} \ar@(l,d)[]!CL;[]_{H(\n_P;f_{GG})} & {}
}
\end{equation*}
\caption{The complex computing the link cohomology of $\M$ at $P$ ($\M$ as
in Figure ~\ref{figExampleLModule})}
\label{figLinkCohomologyExample}
\end{figure}

Besides $\i_{P*}\i_P^*\M$, certain other compositions of these functors are
important.  For example, let $\j_P\colon W\setminus X_P \hookrightarrow W$
denote the inclusion of the complement of $X_P$.  The functor
$\i_P^*\j_{P*}^{\vphantom{*}}\j_P^*\M$ is a complex of $L_P$-modules whose
cohomology is the \emph{link cohomology of $\M$ at $P$}.  An example is
pictured in Figure ~\ref{figLinkCohomologyExample}.

The composed functor $\i_P^*\ihat_Q^!\M$, where $P\le Q$, will be important
in \S\ref{sectMicroSupport} when we define the micro-support of an
$\L$-module.  Here are some examples.  Since $\ihat_G^!$ is the identity,
Figure ~\ref{figUpperStarExample} already gives two examples in the case
$Q=G$.  The other extreme is $Q=P$; in this case $\i_P^*\ihat_Q^!\M =
\ihat_P^!\M = (E_P,f_{PP})$.  For an example with $P<Q< G$, consider
\begin{equation}
\i_P^*\ihat_{Q_1}^!\M = \qquad\vcenter{\xymatrix 
@R=1.25cm @C=3cm @!0
{
{} & {E_P} \ar@(r,u)[]_{f_{PP}} \\
{H(\n_P^{Q_1};E_{Q_1})}  
\ar@(u,l)[];[]!CL_-*+{\labelstyle H(\n_P^{Q_1};f_{Q_1Q_1})} 
\ar[ur]^-{f_{PQ_1}} & {}
}}
\label{eqnExampleMicroSupportExpression}
\end{equation}
where $\M$ is as in Figure ~\ref{figExampleLModule}.

\subsection{Realization of \boldmath$\L$-modules}
\label{ssectRealization}
In order to define the realization sheaf of an $\L$-module, we need to
recall the notion of special differential forms.  The sheaf of
\emph{special differential forms} $\Asp(\Xhat; \EE)$ on $\Xhat$ for a
$G$-module $E$ was introduced by Goresky, Harder, and MacPherson
\cite{refnGoreskyHarderMacPherson}.  It is a fine resolution of $i_{G*}\EE$
whose sections on $U\subseteq \Xhat$ are smooth differential forms on the
inverse image $\widetilde U\subseteq \Xbar$ with coefficients in $\EE$
which satisfy certain boundary conditions.  Namely, in a neighborhood
$\widetilde U$ of a boundary point in $Y_P$, the forms are locally
pullbacks of forms on $Y_P$ under the geodesic retraction $\widetilde
U=B\times c(|\D_P|) \to B \subseteq Y_P$.  Furthermore these boundary
values are required to be constant along the nilmanifold fibers $\nil_P$.
In the case that $E = \bigoplus_i E^i[-i]$ is graded, the differential $d$
is the alternating sum of the corresponding de Rham exterior differentials.
One may show, essentially by applying harmonic projection along the
nilmanifold fibers and using van Est's theorem \eqref{eqnVanEst}, that
there is a natural morphism, ``restriction to a stratum'',
\begin{equation*}
k_P\colon \Asp(\Xhat; \EE) \longrightarrow \ihat_{P*}\Asp(\Xhat_P;
\HH(\n_P;E))
\end{equation*}
which restricts over $\Xhat_P$ to be a quasi-isomorphism.  Given a graded
$L_Q$-module $E_Q$ we may likewise define the sheaf of special differential
forms $\Asp(\Xhat_Q; \EE_Q)$, a differential $d_Q$, and a natural morphism
\begin{equation*}
k_{PQ}\colon \Asp(\Xhat_Q; \EE_Q) \longrightarrow \ihat_{P*}\Asp(\Xhat_P;
\HH(\n_P^Q;E_Q)).
\end{equation*}

The \emph{realization} $\Sheaf(\M)$ of an $\L$-module $\M$ on $\Xhat$ can
now be defined to be the complex of sheaves
\begin{equation*}
\left(\bigoplus_{P\in\Pl} \ihat_{P*}\Asp(\Xhat_P; \EE_P), \sum_{P\le
Q\in\Pl} d_{PQ}\right)
\end{equation*}
where
\begin{equation*}
  d_{PQ} =
  \begin{cases}
    d_P  +  \Asp(\Xhat_P;f_{PP})  &\text{if $P=Q$} \\ 
    \Asp(\Xhat_P;f_{PQ})\circ k_{PQ}   &\text{if $P<Q$}.
  \end{cases}
\end{equation*}
The realization sheaf $\Sheaf(\M)$ for $\M$ as in
Figure~\ref{figExampleLModule} is pictured in
Figure~\ref{figExampleRealization}.  (We do not show the factorization of
the arrows however.)
\begin{figure}
\begin{equation*}
\xymatrix 
@R=1.25cm @C=3cm @!0
{
{} & {} & {\Asp(\Xhat_P;\EE_P)} \ar@(r,u)[]!CR;[]!/r1ex/_{d_{PP}} \\
{} & {}  & {} \\
{\Asp(\Xhat_{Q_1};\EE_{Q_1})} \ar@(u,l)[];[]!CL_{d_{Q_1Q_1}}
\ar[uurr]^-{d_{PQ_1}} & {} & {}\\
{} & {} & {} \\
{} & {} & {\Asp(\Xhat_{Q_2};\EE_{Q_2})}
\ar@(r,u)[]!CR;[]!/r1ex/_{d_{Q_2Q_2}}  \ar[uuuu]^{d_{PQ_2}} \\
{} & {} & {} \\
{\Asp(\Xhat;\EE_G)} \ar@(u,l)[];[]!CL_{d_{GG}} \ar[rruu]^-{d_{Q_2G}}
\ar[uuuu]^{d_{Q_1G}} \ar[uuuuuurr]^{d_{PG}}
}
\end{equation*}
\caption{The realization sheaf $\Sheaf(\M)$ ($\M$ as in Figure
~\ref{figExampleLModule})}
\label{figExampleRealization}
\end{figure}

It is common to compute the cohomology of a sheaf by a complex consisting
of differential forms satisfying boundary conditions or growth conditions
depending on the particular sheaf.  In view of the definition of
$\Sheaf(\M)$, however, we have a uniform type of complex to compute the
cohomology of an $\L$-module: special differential forms with coefficients.
The complication is that these ``forms'' are actually collections of forms
on $\Xhat$ and on the various strata $\Xhat_P$, and there are interaction
terms obtained by restricting to boundary strata and applying $f_{PQ}$.
The boundary conditions mentioned above are translated into the choice of
the coefficients and the $f_{PQ}$; see the example of intersection
cohomology below.

It is easy to check that the realization functor commutes with the usual
functors in the cases we have defined them, that is, there are natural
quasi-isomorphisms $\Sheaf(k^!\M) \cong k^!\Sheaf(\M)$, etc.

\subsection{Reduced $\L$-modules}
\label{ssectReducedLmodules}
An $\L$-module $\M=(E_\cdot,f_{\cdot \cdot})$ is called \emph{reduced} if
$f_{PP}=0$ for all $P$.  One can show that for any $\L$-module $\M$, there
exists a reduced $\L$-module $\M_0$ which is quasi-isomorphic to $\M$;
quasi-isomorphic here means that there exists a third $\L$-module $\M'$ and
morphisms $\M \leftarrow \M' \to \M_0$ which both induce isomorphisms on
local cohomology.  It follows that $H(\Xhat;\M)\cong H(\Xhat;\M_0)$.  One
could develop the theory by requiring $f_{PP}=0$ in the definition of an
$\L$-modules to begin with, however in that case the definitions of the
usual functors, in particular the truncation functors defined below, would
be more complicated.

\subsection{Simplest example}
\label{ssectSimplestExample}
The simplest sheaf that can be lifted to an $\L$-module is $\i_{G*}\EE$,
where $E$ is a $G$-module in degree $0$; the cohomology is
$H(\Xhat;\i_{G*}\EE) \cong H(X;\EE)\cong H(\G;E)$.  The lifted $\L$-module
$i_{G*}E$ is defined by
\begin{equation*}
\begin{aligned}
  E_P & \equiv
  \begin{cases}
    E & \text{if $P=G$,}\\
    0 & \text{if $P\neq G$,}
  \end{cases} \\
  f_{PQ} & \equiv 0.
\end{aligned}
\end{equation*}
The local cohomology at $X_P$ is simply $\i_P^*\i_{G*}^{\vphantom{*}}E =
H(\n_P;E)$.

\subsection{Truncation functors}
\label{ssectTruncationFunctors}
Recall that the formula for the intersection cohomology sheaf
\eqref{eqnIntersectionCohomologySheaf} involves the functor $\t^{\leqslant
n}$ which truncates local cohomology in degrees greater than $n$:
\begin{equation}
\label{eqnTruncationAbove}
\xymatrix @M+4pt { {\t^{\leqslant n}\Sheaf \equiv} & {\cdots}
\ar[r]^-{d_{n-2}} & {\Sheaf^{n-1}} \ar[r]^-{d_{n-1}} & {\ker d_n}
\ar[r]^-{d_n} & {0} \ar[r] & {\cdots} }\ .
\end{equation}
In order to lift intersection cohomology to an $\L$-module, one needs an
analogue of this truncation functor for $\L$-modules.

Consider the simplest $\L$-module $\i_{G*} E$ on the simplest nontrivial
$\Xhat$ where there is only one singular stratum corresponding to the
parabolic subgroup $P$.  We assume $E$ is purely in degree $0$ and $n\ge
0$.  What should the $\L$-module $\t^{\leqslant n} \i_{G*} E$ be?  The data
$E_G=E$ cannot be changed; instead one must define $E_P$ and $f_{PG}$ so
that the local cohomology $H(\i_P^* \t^{\leqslant n} \i_{G*} E)$ at $X_P$
changes from $H(\n_P;E)$ to $\bigoplus_{i\le n} H^i(\n_P;E)[-i]$.  Since
the pullback $\i_P^* \t^{\leqslant n} \i_{G*} E$ is
\begin{equation*}
H(\n_P;E) \xrightarrow{f_{PG}} E_P
\end{equation*}
and there is a natural decomposition
\begin{equation}
\label{eqnDegreeDecomposition}
H(\n_P;E) = \Bigl(\bigoplus_{i\le n} H^i(\n_P;E)[-i]\Bigr) \oplus
\Bigl(\bigoplus_{i>n} H^i(\n_P;E)[-i]\Bigr),
\end{equation}
the desired truncation of local cohomology can be arranged by setting
\begin{equation}
\label{eqnEPDefinition}
E_P = \Bigl(\bigoplus_{i > n} H^i(\n_P;E)[-i]\Bigr)[-1]
\end{equation}
and letting $f_{PG}$ be the natural degree $1$ map (projection onto the
second factor of \eqref{eqnDegreeDecomposition}).

In essence the ``internal'' truncation functor $\t^{\leqslant n}$
\eqref{eqnTruncationAbove} has been realized as an ``external'' truncation
functor.  McConnell does something similar with his ``flabby truncation''
functor $\check\t_{\leqslant n}$
\cite{refnMcConnellIntersectionCohomology}.  That this is possible is not
surprising.  Recall that for a map of complexes $f\colon (C,d_C) \to
(D,d_D)$ (where the objects lie in any additive category), the
\emph{mapping cone} $M(f)$ is the complex defined by
\begin{equation*}
\left( C[1] \oplus D,
\begin{pmatrix}
  d_{C[1]} & 0 \\
  f        & d_D
\end{pmatrix}
\right)
\end{equation*}
The truncation $\t^{\leqslant n}\Sheaf$ of a sheaf is quasi-isomorphic to
the mapping cone
\begin{equation*}
M( \Sheaf \to \t^{>n}\Sheaf)[-1],
\end{equation*}
where
\begin{equation}
\label{eqnTruncationBelow}
\xymatrix @M+4pt { {\t^{> n}\Sheaf \equiv} &
{\cdots\ 0} \ar[r] & {\Sheaf^{n+1}/\Im d_n}
\ar[r]^-{d_{n+1}} & {\Sheaf^{n+2}} \ar[r]^-{d_{n+2}} & {\cdots} }
\end{equation}
and
\begin{equation}
\label{eqnNaturalMapToTruncation}
\Sheaf \to \t^{>n}\Sheaf
\end{equation}
is the natural map.  With this notation, equation \eqref{eqnEPDefinition}
may be re-expressed as
\begin{equation}
E_P = \t^{>n} H(\n_P;E)[-1]
\label{eqnEPDefinitionClean}
\end{equation}
(the degree shift is applied after the truncation).

In order to define $\t^{\leqslant n}\M$ for a general $\L$-module $\M$, we
proceed stratum by stratum.  Namely assume that for a fixed stratum $X_P$,
the local cohomology of $\M$ at all $X_Q$ with $Q>P$ already vanishes in
degrees greater than $n$.  There is an obvious notion of a mapping cone for
a morphism of $\L$-modules and to truncate the local cohomology at $P$ we
define
\begin{equation*}
\t^{\leqslant n}_P\M \equiv M(\M \to
\i_{P*}^{\vphantom{*}}\t^{>n}\i_P^*\M)[-1].
\end{equation*}
The morphism here is the composition
\begin{equation*}
\M \longrightarrow \i_{P*}^{\vphantom{*}}\i_P^*\M \longrightarrow
\i_{P*}^{\vphantom{*}}\t^{>n}\i_P^*\M
\end{equation*}
of the adjunction morphism with one induced from
\eqref{eqnNaturalMapToTruncation}.  An example of $\t^{\leqslant n}_P$
where $P$ has parabolic rank $2$ will be given in the next subsection.  The
functor $\t^{\leqslant n}$ is defined to be the composition of all the
functors $\t^{\leqslant n}_P$, where one applies $\t^{\leqslant n}_R$
before $\t^{\leqslant n}_Q$ if $Q<R$.

It is easy to check that there is a natural quasi-isomorphism
\begin{equation*}
\Sheaf(\t^{\leqslant n}\M) \cong \t^{\leqslant n}\Sheaf(\M).
\end{equation*}

\subsection{Intersection cohomology}
\label{ssectionIntersectionCohomology}
The formula \eqref{eqnIntersectionCohomologySheaf} for the intersection
cohomology sheaf $\IpC(\Xhat;\EE)$ involves pushforward and truncation;
since both of these functors have been defined for $\L$-modules, commuting
with the realization functor, we can lift intersection cohomology to an
$\L$-module $\IpC(\Xhat;E)$.  Specifically, enumerate the elements of $\Pl$
as $P_0=G, P_1,\dotsc, P_N$ such that $\dim X_{P_i} \ge \dim X_{P_j}$ for
$i<j$.  For $k<N$, let
\begin{equation*}
\j_{P_k}\colon \coprod_{i< k} X_{P_i} \hookrightarrow 
\coprod_{i\le k} X_{P_i}
\end{equation*}
be the inclusion map.  Define
\begin{equation}
\label{eqnStrataPerversity}
p(Q) \equiv p(\dim \n_Q + \#\D_Q) \qquad \text{for $Q\in \Pl$;}
\end{equation}
this is reasonable since $\dim \n_Q + \#\D_Q$ is the codimension of $X_Q$
in $\Xhat$ (namely $1$ more than the dimension of the link
\eqref{eqnLinkReductiveBorelSerre}).  The intersection cohomology
$\L$-module is defined as
\begin{equation*}
\IpC(\Xhat;E) \equiv \t^{\leqslant p(P_N)} \j_{P_N*} \dotsm \t^{\leqslant
p(P_2)} j_{P_2*} \t^{\leqslant p(P_1)} j_{P_1*}E,
\end{equation*}
where the $G$-module $E$ is viewed as usual as an $\L$-module on $X=X_G$.

Let us illustrate this procedure for our standard $\QQrank 2$ example with
$\Pl$ as in Figure ~\ref{figHasseDiagram}.  Obtaining the $\L$-module over
$X_G \sqcup X_{Q_1} \sqcup X_{Q_2}$ is straightforward; the result is
pictured in Figure ~\ref{figIntersectionCohomologyPrankOne}.
\begin{figure}
\begin{equation*}
\xymatrix 
@R=1.25cm @C=3cm @!0
{
{\t^{> p(Q_1)}H(\n_{Q_1};E)[-1]}\\
{}\\
{H(\n_{Q_1};E)} \ar[uu] & {} & {\t^{> p(Q_2)}H(\n_{Q_2};E)[-1]} \\
{} & {H(\n_{Q_2};E)}  \ar[ru] \\
{E}
}
\end{equation*}
\caption{The $\L$-module $\IpC(X \sqcup X_{Q_1} \sqcup X_{Q_2};E)$.  (The
  unlabeled morphisms are the natural degree $1$ maps.)}
\label{figIntersectionCohomologyPrankOne}
\end{figure}
The complex $\i_P^*j_{P*}^{\vphantom{*}}\IpC(X \sqcup X_{Q_1} \sqcup
X_{Q_2};E)$ computing the intersection cohomology of the link at $P$ is
then
\begin{equation}
(C,d_C) \equiv \quad \vcenter{
\xymatrix @R=1.25cm @C=3cm @!0
{
*+!L(.25){H(\n_P^{Q_1};\t^{> p(Q_1)}H(\n_{Q_1};E))[-1]} \\
{} & *+!L(.25){H(\n_P^{Q_2};\t^{> p(Q_2)}H(\n_{Q_2};E))[-1]} \\
{H(\n_P;E)} \ar[uu] \ar[ur]
}};
\label{eqnLinkIntersectionCohomology}
\end{equation}
the final $\L$-module $\IpC(\Xhat;E)$ is pictured in Figure
~\ref{figIntersectionCohomologyPrankTwo}.
\begin{figure}
\begin{equation*}
\xymatrix 
@R=1.25cm @C=3cm @!0
{
{} & {} & {\t^{> p(P)}(C,d_C)[-1]}  \\
{} & *+!R(.4){H(\n_P^{Q_1};\t^{> p(Q_1)}H(\n_{Q_1};E))[-1]} \ar[ur] & {} \\
{\t^{> p(Q_1)}H(\n_{Q_1};E)[-1]} & {} &
*+!L(.4){H(\n_P^{Q_2};\t^{> p(Q_2)}H(\n_{Q_2};E))[-1]} \ar[uu]\\
{} & {H(\n_P;E)} \ar[uuur] & {} \\
{H(\n_{Q_1};E)} \ar[uu] & {} & {\t^{> p(Q_2)}H(\n_{Q_2};E)[-1]} \\
{} & {H(\n_{Q_2};E)}  \ar[ru] \\
{E}
}
\end{equation*}
\caption{The $\L$-module $\IpC(\Xhat;E)$.  (The
  unlabeled morphisms are the natural degree $1$ maps.)}
\label{figIntersectionCohomologyPrankTwo}
\end{figure}

\section{Micro-support of an $\L$-module}
\label{sectMicroSupport}
A key local invariant of an $\L$-module defined in \cite{refnSaperLModules}
is its \emph{micro-support}; the terminology reflects the rough analogy
that exists between this concept and that introduced by Kashiwara and
Schapira \cite{refnKashiwaraSchapira} for sheaves.  The micro-support
$\mS(\M)$ of an $\L$-module $\M$ is a certain set of irreducible
representations $V$ of the various $L_P$ which we now define.

\subsection{Definition of micro-support}
Let $\IrrRep(L_P)$ denote the irreducible algebraic representations of
$L_P$.  Recall we have decomposed $L_P=M_PA_P$, where $A_P$ is the identity
component of the maximal $\QQ$-split torus in the center of $L_P$.  For any
$V\in\IrrRep(L_P)$, let $\xi_V$ be the character by which $A_P$ acts on $V$
and let $\CC_{\xi_V}$ be the corresponding $1$-dimensional representation.
We can then write $V$ as $V|_{M_P}\otimes \CC_{\xi_V}$.  Recall that $\D_P$
denotes the simple $A_P$-weights of the adjoint action of (a lift of) $L_P$
on $\n_P = \Lie(N_P)$, and that the parabolic $\QQ$-subgroups $Q\ge P$ are
parametrized by subsets $\D_P^Q$ of $\D_P$.  Consequently we can define
parabolic $\QQ$-subgroups $P\le Q_V\le Q'_V$ by setting
\begin{align*}
\D_P^{Q_V} &=\{\,\al\in \D_P\mid (\xi_V+\r,\al)<0\,\}, \\
\D_P^{Q_V'} &=\{\,\al\in \D_P\mid (\xi_V+\r,\al)\le 0\,\};
\end{align*}
here $\r$ is one-half the sum of the $A_P$-weights of $\n_P$ (counted with
multiplicity).

Let $\M$ be an $\L$-module on $\Xhat$.  Define $\mS(\M)$ to consist of
those $V\in \IrrRep(L_P)$ (for $P\in\Pl$) such that
\begin{align}
\tag{i}&(V|_{M_P})^* \cong \overline{V|_{M_P}},\text{ and}
\label{eqnMicroSupportOne}\\
\tag{ii}&H(i_P^* \ihat_Q^! \M)_V \neq 0 \text{ for some $Q_V\le Q \le
  Q_V'$.} \label{eqnMicroSupportTwo}
\end{align}
Here $H(i_P^* \ihat_Q^! \M)_V$ denotes the $V$-isotypical component of
$H(i_P^* \ihat_Q^! \M)$.  Let
\begin{equation*}
c(V;\M)\le d(V;\M)
\end{equation*}
be the least and greatest degrees for which \eqref{eqnMicroSupportTwo}
holds.

Condition \eqref{eqnMicroSupportOne} says that $V|_{M_P}$ is
\emph{conjugate self-contragredient}.  This condition is automatic if $D_P$
is equal-rank, $\CCrank M_P = \rank K_P$ (see the discussion later in
\S\ref{ssectExampleMicroSupport}).  Condition \eqref{eqnMicroSupportTwo}
may be interpreted geometrically as follows.  Consider the long exact
sequence of the pair $(U,U\setminus(U\cap \Xhat_Q))$ where $U$ is a small
neighborhood of a point of $X_P$.  Since $H(i_P^* \ihat_Q^! \M)_V \cong
H(U,U\setminus(U\cap \Xhat_Q);\M)_V$, equation \eqref{eqnMicroSupportTwo}
is equivalent to
\begin{equation*}
\xy\xycompile{ 0;<-1mm,0mm>:<0mm,1mm>::
(0,20)*+!DL{\scriptstyle X_P};(20,23) **\crv{(10,20)},
(10,22)*+!D{\scriptstyle X_Q},
(0,20)*-{\bullet};(14,3) **\crv{(10,10)},
(16,14)*+{\scriptstyle X},
(10,-5)*{H(U;\M)_V}
}\endxy\qquad\xy 0;<1mm,0mm>:
(5,-5)\ar @{>}(15,-5)
\endxy\qquad\xy\xycompile{ 0;<-1mm,0mm>:<0mm,1mm>::
(0,20)*+!DL{\scriptstyle X_P};(20,23) **\crv{~*=<4pt>{.} (10,20)},
(10,22)*+!D{\scriptstyle X_Q},
(0,20)*\cir<2pt>{};(14,3) **\crv{(10,10)},
(16,14)*+{\scriptstyle X},
(10,-5)*{H(U\setminus(U\cap \Xhat_Q);\M)_V}}
\endxy
\end{equation*}
failing to be an isomorphism in some degree.  Note that a functor like
$i_P^* \ihat_Q^!$ was used by Goresky and MacPherson
\cite{refnGoreskyMacPersonLefschetzFixedPointFormula} to describe the local
contributions to a Lefschetz fixed point formula; more recently it has been
used in Braden's work \cite{refnBradenHyperbolicLocalization} on hyperbolic
localization.

We say that an element $V\in \mS(\M)$ belongs to $\emS(\M)$, the
\emph{essential micro-support}, if furthermore
\begin{equation}
\label{eqnMicroSupportTwoPrime}
\tag*{(ii)$'$} H(i_P^* \ihat_{Q_V}^! \M)_V  \longrightarrow H(i_P^*
\ihat_{Q_V'}^! \M)_V \qquad \text{is nonzero.}
\end{equation}

All of the above was for an $\L$-module on $\Xhat$.  Consider more
generally an admissible subspace $W\subseteq \Xhat$ which has a unique
maximal stratum $X_R$.  We define the micro-support and the essential
micro-support of an $\L$-module on $W$ similarly except that we restrict
$V$ to belong to $\coprod_{P\in\Pl(W)} \IrrRep(L_P)$ and replace $Q_V$ and
$Q'_V$ by parabolic $\QQ$-subgroups $Q_V^R$ and $Q_V^{\prime R}$ defined by
\begin{align*}
\D_P^{Q_V^R} &=\{\,\al\in \D_P^R\mid (\xi_V+\r,\al)<0\,\}, \\
\D_P^{Q_V^{\prime R}} &=\{\,\al\in \D_P^R\mid (\xi_V+\r,\al)\le 0\,\}.
\end{align*}

\subsection{Example: Micro-support of \boldmath $\i_{G*} E$}
\label{ssectExampleMicroSupport}
We illustrate the definition by calculating the micro-support of $\i_{G*}E$
\cite{refnSaperIHP}.  This also gives us an opportunity to introduce
Kostant's theorem which will be important later.

For $\M=\i_{G*}E$, we have $E_P=0$ for all $P\neq G$ and thus $\ihat_Q^!
\M =0$ unless $Q=G$; if $Q=G$ we have $\ihat_G^! \M =\M$.  Thus for any
$P\in \Pl$, we only need to consider irreducible $L_P$-modules $V$ such
that $Q_V'=G$, that is, such that
\begin{equation}
(\xi_V+\r,\al)\le 0 \qquad \text{for all $\al\in\D_P$;}
\end{equation}
such a $V$ will actually be in $\mS(\i_{G*}E)$ if $(V|_{M_P})^* \cong
\overline{V|_{M_P}}$ and if $V$ occurs in $H(\i_P^*\i_{G*}E)= H(\n_P;E)$.

We therefore need to know the irreducible components of $H(\n_P;E)$; in
fact a theorem of Kostant \cite{refnKostant} completely describes
$H(\n_P;E)$ as an $L_P$-module.  Since this will be essential later as
well, we review it here.  First some notation.  Let $\h= \h_{M_P} + \sa_P$
be a Cartan subalgebra of $\levi_P$.  Via a lift we may view $\h$ as a
Cartan subalgebra of $\mathfrak g$ as well and we choose a set of simple
$\CC$-roots $\lsb\CC\D$ for $\mathfrak g_\CC$ such that the roots in
$\n_{P\CC}$ are positive; as usual, let $\r$ denote one-half the sum of the
positive roots (its restriction to $\sa_P$ agrees with the previous
definition of $\r$).  Let $W$ be the Weyl group of $\mathfrak g_\CC$ and
let $W^P$ denote the subgroup corresponding to the Weyl group of
$\levi_{P\CC}$.  The length of $w\in W$ is denoted $\l(w)$ and we let $W_P$
be the set of minimal length coset representatives of $W^P\back W$.  Assume
$E$ is irreducible with highest weight $\lambda$.  Kostant's theorem says
that
\begin{equation}
H(\n_P; E) = \bigoplus_{w\in W_P} V_{w(\lambda+\r)-\r}[-\l(w)],
\end{equation}
where $V_\u$ denotes the irreducible $L_P$-module with highest weight $\u$.

We also need to understand when $V|_{M_P}$ is conjugate
self-contragredient.  The involution $V|_{M_P}\mapsto
(\overline{V|_{M_P}})^*$ on representations induces an involution on highest weights,
$\u|_{h_{M_P}} \mapsto \t_P(\u|_{h_{M_P}})$.  It is always the case that
$\t_P(\r|_{\h_{M_P}}) = \r|_{\h_{M_P}}$, so $V_{w(\lambda+\r)-\r}|_{M_P}$
is conjugate self-contragredient if and only if
\begin{equation*}
\t_P(w(\lambda+\r)|_{\h_{M_P}}) = w(\lambda+\r)|_{\h_{M_P}}.
\end{equation*}
The involution $\t_P$ is described in \cite{refnBorelCasselman}; it is the
composition of the opposition involution and the ``$*$-action'' of complex
conjugation \cite{refnTits}.  Concretely, let $\lsb\CC\D^P \subseteq
\lsb\CC\D$ denote the simple $\CC$-roots of $\levi_{P\CC}$.  Then
\begin{equation*}
\t_P = (-w^P_0)\circ(w^P_c c)
\end{equation*}
where $w^P_0\in W^P$ is the longest element, $c$ is complex conjugation,
and $w^P_c\in W^P$ is such that $w^P_c(c\lsb\CC\D^P) = \lsb\CC\D^P$.  Even
more concretely, if $\h$ is a \emph{fundamental} (that is, maximally
compact) $\theta$-stable Cartan subalgebra of $\levi_P$ and the positive
system is chosen so that $\lsb\CC\D^P$ is $\theta$-stable, then
$\t_P=\theta$.

The micro-support of $\i_{G*}E$ is thus
\begin{equation*}
\begin{split}
\mS(\i_{G*}E) = \coprod_P \{ \, V_{w(\lambda+\r)-\r} \,\mid\,
	& (w(\lambda+\r),\al) \le 0 \text{ for all }\al\in\D_P,\\
	& \t_P(w(\lambda+\r)|_{\h_{M_P}}) = w(\lambda+\r)|_{\h_{M_P}}\,\};
\end{split}
\end{equation*}
to obtain an expression for the essential micro-support $\emS(\i_{G*}E)$,
we impose instead the strict inequality $(w(\lambda+\r),\al) < 0 $.  For
$V=V_{w(\lambda+\r)-\r} \in \mS(\i_{G*}E)$, we have
\begin{equation}
\label{eqnPushforwardSheafMicroSupportDegree}
c(V;\i_{G*}E) = d(V;\i_{G*}E) = \l(w).
\end{equation}

\section{A vanishing theorem for the cohomology of an $\L$-module}
The justification for the definition of $\mS(\M)$ is that it enters into a
general vanishing theorem for the cohomology of an $\L$-module (see Theorem
~\ref{thmVanishingTheorem} below).

First some definitions.  Let $\M$ be an $\L$-module on $\Xhat$ and let
$V\in \mS(\M)$ be an $L_P$-module.  Let $\u$ be the highest weight of $V$
with respect to a fundamental Cartan subalgebra $\h$ of $\levi_P$ and a
$\theta$-stable positive system; view $\u$ as an element of
$\levi_{P\CC}^*$ by extending it to be $0$ on all root spaces.  Let
$L_P(\u)\subseteq L_P$ be the stabilizer of $\u$ for the coadjoint action;
this is a reductive $\RR$-subgroup whose $\CC$-roots are those $\CC$-roots
of $L_P$ which are orthogonal to $\u$.  Let
\begin{equation*}
D_P(V) = L_P(\u)/(K_P\cap L_P(\u))A_P
\end{equation*}
be the associated symmetric space.  The space $D_P(V)$ is not well-defined,
even up to isomorphism, since it depends on the choice of a $\theta$-stable
positive system; we will assume however that the positive system has been
chosen so that $D_P(V)$ has the maximum possible dimension.

Define
\begin{align*}
c(\M) &= \inf_{V\in\mS(\M)}  \tfrac12(\dim D_P - \dim D_P(V)) + c(V;\M), \\
d(\M) &= \sup_{V\in\mS(\M)}  \tfrac12(\dim D_P + \dim D_P(V)) + d(V;\M).
\end{align*}
As we will see in Theorem ~\ref{thmVanishingLTwoCohomology}, the
significance of $(\dim D_P \pm \dim D_P(V))/2$ is that this is the range of
degrees in which $H_{(2)}(X_P;\VV)$ can be nonzero.  In these definitions,
one can also consider just $V\in \emS(\M)$ and redefine $c(V;\M)\le
d(V;\M)$ to be the least and greatest degrees in which
\ref{eqnMicroSupportTwoPrime} is nonzero; one can prove that the same
values of $c(\M)$ and $d(\M)$ are obtained.

The following vanishing theorem is proved in \cite[Theorem
10.4]{refnSaperLModules}:

\begin{thm}
\label{thmVanishingTheorem}
Let $\M$ be an $\L$-module on $\Xhat$. Then
\begin{equation*}
H^j(\Xhat;\M) = 0 \qquad \text{for $j\notin [c(\M),d(\M)]$.}
\end{equation*}
In particular, if $\mS(\M)=\emptyset$ then $H(\Xhat;\M)=0$ identically.
\end{thm}

We will sketch a proof of the theorem in
\S\ref{sectSketchProofVanishingTheorem}.  It has a number of ingredients,
including the author's past work on tilings of locally symmetric spaces
\cite{refnSaperTilings}, various analytic estimates, and a vanishing
theorem for $L^2$-cohomology which we will recall in
\S\ref{sectVanishingLTwoCohomology}.

\section{Application: a vanishing theorem for ordinary cohomology}
One can immediately apply Theorem ~\ref{thmVanishingTheorem} and our
calculation of $\emS(\i_{G*}E)$ to obtain a vanishing theorem for the
ordinary cohomology $H(X;\EE)$.  This applies to all $X$ and $E$, however
in special cases we can be more explicit regarding the degree bounds.  For
example, we have proved \cite[Theorem ~5]{refnSaperIHP}:

\begin{thm}
\label{thmOrdinaryCohomologyVanishing}
If $D$ is equal-rank and $E$  has regular highest weight, then
\begin{equation*}
H^j(X;E) = 0 \qquad \text{for $j<\tfrac12 \dim X$.}
\end{equation*}
\end{thm}

The same theorem has been announced by Li and Schwermer
\cite{refnLiSchwermer} in a slightly strengthened form.  They drop the
equal-rank condition and replace the vanishing range by
\begin{equation*}
j < \tfrac12(\dim X - (\CCrank \lsp0\,G - \rank K));
\end{equation*}
it is clear that this strengthened form also follows from our proof in
\cite{refnSaperIHP}.  Theorem ~\ref{thmOrdinaryCohomologyVanishing} answers
a question of Tilouine \cite[\S8.7]{refnMokraneTilouine} in the Hermitian
case; the case $G=R_{k/\QQ}\GSympl_4(\RR)$ where $k$ is a totally real
number field was proved in \cite{refnTilouineUrban} using results of Franke
\cite{refnFranke}.  For an application of the theorem see Mauger's thesis
\cite{refnMauger}.

In order to deduce Theorem ~\ref{thmOrdinaryCohomologyVanishing} from
Theorem ~\ref{thmVanishingTheorem}, equation
\eqref{eqnPushforwardSheafMicroSupportDegree} shows that we need an
estimate on $\l(w)$.  The following basic lemma provides the needed
estimate; it will also be crucial later in studying the micro-support of
$\IpC(\Xhat;E)$ and its behavior under functorial operations, Theorems
~\ref{thmMicroSupportIC} and \ref{thmFunctoriality}.

\begin{lem}
\label{lemBasicLemma}
Let $E$ be an irreducible $G$-module with highest weight $\lambda$.  Let
$P$ be a parabolic $\QQ$-subgroup and let $w\in W_P$.  Set $V=
V_{w(\lambda+\r)-\r}$ and assume that $V|_{M_P}$ conjugate
self-contragredient.
\begin{enumerate}
\item\label{itemNegativeBasic} If $(w(\lambda+\r),\al) \le0$ for all
  $\al\in \D_P$, then
\begin{equation*}
\l(w) \ge \tfrac12(\dim\n_P + \dim \n_P(V)).
\end{equation*}
\item\label{itemPositiveBasic} If $(w(\lambda+\r),\al) \ge0$ for all
  $\al\in \D_P$, then
\begin{equation*}
\l(w) \le \tfrac12(\dim\n_P - \dim \n_P(V)).
\end{equation*}
\end{enumerate}
\end{lem}

The quantity $\n_P(V)$ in the lemma is defined as follows; see
\cite[\S24.1]{refnSaperLModules} for more details.  Let $\u$ be the highest
weight of $V$ for a fundamental $\theta$-stable Cartan subalgebra of
$\levi_P$ and $\theta$-stable positive system.  Let $\n_P(\u)$ be the sum
of the irreducible $L_P(\u)$-submodules of $\n_{P\CC}$ whose weights are
stable under the action of $-\t_P$; observe that $-\t_P$ is complex
conjugation in this situation so that $\n_P(\u)$ contains in particular all
positive real root spaces.  Let $\n_P(V)$ be any of the $\n_P(\u)$ with
maximal dimension as we vary the $\theta$-stable positive system.

We can now deduce the theorem (in its strengthened form).  The
Kostant-Sugiura classification of Cartan subalgebras
\cite{refnKostantConjugacyCartan}, \cite{refnSugiura},
\cite{refnSugiuraCorrection} shows that $\rank K - \rank K_P$ is at most
the number of positive real roots.  Thus
\begin{equation}
\dim\n_P(V) \ge \rank K - \rank K_P.
\end{equation}
The hypothesis that $E$ has regular highest weight implies that $\dim
D_P(V) = \CCrank M_P - \rank K_P$ for any $V\in \mS(\i_{G*}E)$.  These facts
and Lemma ~\ref{lemBasicLemma}\itemref{itemNegativeBasic} allow one to
estimate
\begin{equation*}
\begin{split}
&\tfrac12(\dim D_P - \dim D_P(V)) +
   c(V;\i_{G*}E) \\
&\qquad \ge \tfrac12\bigl(\dim D_P + \dim \n_P - (\CCrank M_P
   - \rank K)\bigr) \\
&\qquad = \tfrac12\bigl(\dim X - (\CCrank \lsp0G - \rank K)\bigr).
\end{split}
\end{equation*}

\begin{proof}[Sketch of the proof of Lemma ~\ref{lemBasicLemma}]
We will sketch the argument with $\dim \n_P(V)$ replaced by $\rank K -
\rank K_P$; see \cite[Lemma ~24.2]{refnSaperLModules} for the complete
proof.  Recall that $\l(w)$ can be viewed as the number of positive
$\CC$-roots $\g$ for which $w^{-1}\g$ is negative; for $w\in W_P$ all these
roots are in $\n_{P\CC}$.  Since $\lambda+\r$ is strictly dominant,
$(w(\lambda+\r),\g) = (\lambda+\r,w^{-1}\g)$ is never zero and will be
negative precisely when $w^{-1}\g<0$.  Thus $\l(w)$ is the number of roots
$\g$ in $\n_{P\CC}$ for which
\begin{equation}
(w(\lambda+\r),\g)\le 0.
\label{eqnPhiwTest}
\end{equation}
However complex conjugation defines an involution $\g\mapsto \bar \g$ on
the roots in $\n_{P\CC}$; we will prove \itemref{itemNegativeBasic} by
showing that \eqref{eqnPhiwTest} holds for at least one of $\g$ and
$\bar\g$.

Assume that $\h= \h_{M_P} + \sa_P$ is a fundamental $\theta$-stable Cartan
subalgebra of $\levi_P$ and that the set of positive roots of
$\levi_{P\CC}$ are $\theta$-stable.  Decompose $\h = \h_{M_P,1} + (\h_{M_P,-1}
+ \sa_P)$ according to the $\pm1$-eigenspaces of the Cartan involution.  In
this case, the operator $\t_P$ is simply $\theta$
\cite{refnBorelCasselman}, so the condition that
$(V_{w(\lambda+\r)-\r}|_{M_P})^* \cong
\overline{V_{w(\lambda+\r)-\r}|_{M_P}}$ is equivalent to
$w(\lambda+\r)|_{\h_{M_P,-1}} = 0$.  It follows that
\begin{equation*}
(w(\lambda+\r),\g) = (w(\lambda+\r),\g|_{\h_{M_P,1}}) +
  (w(\lambda+\r),\g|_{\sa_P}).
\end{equation*}
Since $\g|_{\sa_P}$ is a linear combination with nonnegative coefficients
of $\al\in \D_P$, the second term is nonpositive by hypothesis; on the
other hand, the first term is negated when $\g$ is replaced by $\bar \g$.
Thus \eqref{eqnPhiwTest} holds for one of either $\g$ or $\bar\g$ as
desired.  In particular it holds for all roots such that $\g=\bar \g$, that
is, the positive real roots.  Since, as we have already noted, the number
of positive real roots is at least $\rank K - \rank K_P$, assertion
\itemref{itemNegativeBasic} follows with $\dim \n_P(V)$ substituted by
$\rank K - \rank K_P$.  Assertion \itemref{itemPositiveBasic} is proven
similarly.
\end{proof}

\section{A vanishing theorem for $L^2$-cohomology}
\label{sectVanishingLTwoCohomology}

As we will see in the next section, the proof of Theorem
~\ref{thmVanishingTheorem} will be reduced to the study of the
$L^2$-cohomology groups $H_{(2)}(X_P; \VV)$ for all $P\in \Pl$ and all
$L_P$-modules $V\in \mS(\M)$.  Thus we need a vanishing theorem for
$L^2$-cohomology.  The theorem below was proved by the author and Stern in
\cite{refnSaperSternTwo} though it is not stated explicitly in this form;
it is based on work of Raghunathan \cite{refnRaghunathan},
\cite{refnRaghunathanCorrection}.  The statement and proof appeared in
\cite[Theorem 14.4]{refnSaperLModules}.

\begin{thm}
\label{thmVanishingLTwoCohomology}
Let $X$ be a locally symmetric space and let $E$ be a representation of $G$
with corresponding metrized locally constant sheaf $\EE$.
\begin{enumerate}
\item If $(E|_{\lsp0\,G})^* \not\cong \overline {E|_{\lsp0\,G}}$, then
$H_{(2)}(X;\EE)=0$.
\item If $(E|_{\lsp0\,G})^* \cong \overline {E|_{\lsp0\,G}}$, then
$H_{(2)}^j(X;\EE)=0$ for
\begin{equation*}
j \notin [\tfrac12(\dim D - \dim D(E)),\tfrac12(\dim D + \dim D(E))]
\end{equation*}
and $H_{(2)}^{(\dim D - \dim D(E))/2}(X;\EE)$ is Hausdorff.
\end{enumerate}
\end{thm}

The proof of the theorem is based on the following criterion which will be
used in the proof of Theorem ~\ref{thmVanishingTheorem} later as well:
\begin{prop}
\label{propEstimate}
Let $M$ be a complete Riemannian manifold with a metrized locally constant
sheaf $\EE$.  Then for every $j$, the following two conditions are
equivalent.
\begin{enumerate}
\item $H^j_{(2)}(M;\EE)=0$ and $H^{j+1}_{(2)}(M;\EE)$ is Hausdorff\textup;
\item there exists $c>0$ such that 
\begin{equation*}
\|d\phi\|^2 + \|d^*\phi\|^2 \ge c \|\phi\|^2
\end{equation*}
for all compactly supported smooth forms $\phi$ of degree $j$.
\end{enumerate}
\end{prop}

The complex of compactly supported forms in $A_{(2)}(X;\EE)$ may be
identified as in \cite[VII.2.7]{refnBorelWallach},
\cite{refnMatsushimaMurakami} with a Koszul complex $C(\lsp0\mathfrak g, K;
A_c^0(\G\back \lsp 0\,G)\otimes E)$.  Thus the differential may be
decomposed $d=d_0+d_1$, where $d_0$ corresponds to the algebraic action of
$\lsp0\,\mathfrak g$ on $E$ and $d_1$ corresponds to the action of
$\lsp0\,\mathfrak g$ on $A_c^0(\G\back \lsp 0\,G)$ by differentiation.  One
calculates using integration by parts that
\begin{equation*}
\|d\phi\|^2 + \|d^*\phi\|^2 \ge (\D_0\phi,\phi),
\end{equation*}
where $\D_0$ is the nonnegative algebraic Laplacian corresponding to $d_0$;
it is the Laplacian for the complex of invariant forms $C(\lsp0\,\mathfrak g,
K; E)$.  A careful analysis of the possible zero eigenvectors of $\D_0$
yields Theorem ~\ref{thmVanishingLTwoCohomology}.

Alternatively, note that a vanishing theorem for $(\mathfrak
g,K)$-cohomology of unitary representations was proved by Vogan and
Zuckerman \cite{refnVoganZuckerman} based on work by Kumaresan
\cite{refnKumaresan}; this implies a vanishing theorem for $L^2$-cohomology
(at least the Hausdorff part) which one can show
\cite{refnSaperLetterCasselman} is equivalent with Theorem
~\ref{thmVanishingLTwoCohomology}.

\section{Proof of the vanishing theorem for the cohomology of an
  $\L$-module}
\label{sectSketchProofVanishingTheorem}
In the following pages we sketch a proof of Theorem
~\ref{thmVanishingTheorem}, the vanishing theorem for the cohomology
$H(\Xhat;\M)$ of an $\L$-module $\M=(E_\cdot,f_{\cdot\cdot})$.  Most steps
of the proof are accompanied by an illustrative figure.

\newpage

\begin{figure}[t]
\begin{equation*}
\begin{xy}
0;<.54in,0in>:<0in,.7in>::
(0,2)*i{\oplus},(9,2)*i{\oplus},(0,6)*i{\oplus},(9,6)*i{\oplus}
,(2.2,4.5)="bone",(6.8,4.5)="btwo"
,(.5,2);(8.5,2)**@{}
?!{"bone";"bone"+a(240)}="a"?!{"btwo";"btwo"+a(-60)}="b"
,(.25,2);(1.2,6)**@{}?="mid",?(1)
**\crv~c{"mid"+(.5,0)}?!{"bone";"bone"+(-.866025404,.5)}="c"
,(1.2,6);(7.8,6)**@{}?="mid",?(1)**\crv~c{"mid"+(0,-.5)}?!{"bone";"bone"+(0,1)}="d"?!{"btwo";"btwo"+(0,1)}="e"
,(7.8,6);(8.75,2)**@{}?="mid",?(1)**\crv~c{"mid"+(-.5,0)}?!{"btwo";"btwo"+(.866025404,.5)}="f"
,(-.25,0)="l",(0,.25)="u",(.25,0)="r"
,(.25,2)+"l";(1.2,6)+"l"**@{}?="mid",?(1)
**\crv~c{"mid"+(.5,0)}?!{"bone";"bone"+(-.866025404,.5)}="cout"
,(1.2,6)+"l"+"u"*!C{\textstyle\bullet}
,(1.2,6)+"u";(7.8,6)+"u"**@{}?="mid",?(1)**\crv~c{"mid"+(0,-.5)}?!{"bone";"bone"+(0,1)}="dout"?!{"btwo";"btwo"+(0,1)}="eout"?*++!D{\textstyle \Xhat_{Q_1}}
,(7.8,6)+"u"+"r"*!C{\textstyle\bullet}*++!LD{\textstyle \Xhat_P}
,(7.8,6)+"r";(8.75,2)+"r"**@{}?="mid",?(1)**\crv~c{"mid"+(-.5,0)}?!{"btwo";"btwo"+(.866025404,.5)}="fout"?*++!L{\textstyle \Xhat_{Q_2}}
,(4.5,3)*{\textstyle \Xhat}
\end{xy}
\end{equation*}
\caption{The complex $\Gamma(\Xhat,\Sheaf(\M))$ lives on the disjoint union
  $\coprod_P \Xhat_P$}
\label{figBoundaryStrata}
\end{figure}

The cohomology $H(\Xhat;\M)$ may be computed as the cohomology of the
complex of global sections of the fine realization sheaf $\Sheaf(\M)$
defined in \S\ref{ssectRealization}.  Ignoring the differentials, this
complex is
\begin{equation*}
\Gamma(\Xhat,\Sheaf(\M)) = \bigoplus_{P\in \Pl} A\sp(\Xhat;\EE_P).
\end{equation*}
Thus a ``form'' $\phi$ is really a collection $(\phi_P)_{P\in\Pl}$ of
special differential forms on the various strata as pictured in Figure
~\ref{figBoundaryStrata}.  The differential is
\begin{equation*}
d\phi = d (\phi_P)_{P\in \Pl} = \bigl(\sum_{P\le Q} d_{PQ}\phi_Q\bigr)_{P\in
\Pl} = \bigl(d_P\phi_P + \sum_{P\le Q} f_{PQ}(\phi_Q|_{\Xhat_P})
\bigr)_{P\in \Pl}
\end{equation*}
where $d_P$ denotes the de Rham differential and we write
$\phi_Q|_{\Xhat_P}$ for $k_{PQ}(\phi_Q)$.

\clearpage

We would like to use Hodge theory to calculate the cohomology of this
complex and for that we need to consider the $L^2$-norm
\begin{equation*}
\|\phi\|^2 \equiv \sum_P \|\phi_P\|^2.
\end{equation*}
However in the locally symmetric metric on $X_P$, the boundary is at
infinite distance.  Thus the condition $\|\phi\|< \infty$ imposes
nontrivial $L^2$-growth conditions on the $\phi_P$ which would change the
cohomology.

\begin{figure}[t]
\begin{equation*}
\begin{xy}
0;<.54in,0in>:<0in,.7in>::
(0,2)*i{\oplus},(9,2)*i{\oplus},(0,6)*i{\oplus},(9,6)*i{\oplus}
,(2.2,4.5)="bone",(6.8,4.5)="btwo"
,(.5,2);(8.5,2)**@{}
?!{"bone";"bone"+a(240)}="a"?!{"btwo";"btwo"+a(-60)}="b"
,(.25,2);(1.2,6)**@{}?="mid",?(1)
**\crv~c{"mid"+(.5,0)}?!{"bone";"bone"+(-.866025404,.5)}="c"
,(1.2,6);(7.8,6)**@{}?="mid",?(1)**\crv~c{"mid"+(0,-.5)}?!{"bone";"bone"+(0,1)}="d"?!{"btwo";"btwo"+(0,1)}="e"
,(7.8,6);(8.75,2)**@{}?="mid",?(1)**\crv~c{"mid"+(-.5,0)}?!{"btwo";"btwo"+(.866025404,.5)}="f"
,(-.25,0)="l",(0,.25)="u",(.25,0)="r"
,(.25,2)+"l";(1.2,6)+"l"**@{}?="mid",?(1)
**\crv~c{"mid"+(.5,0)}?!{"bone";"bone"+(-.866025404,.5)}="cout"
,(1.2,6)+"l"+"u"*!C{\textstyle\bullet}
,(1.2,6)+"u";(7.8,6)+"u"**@{}?="mid",?(1)**\crv~c{"mid"+(0,-.5)}?!{"bone";"bone"+(0,1)}="dout"?!{"btwo";"btwo"+(0,1)}="eout"?*++!D{\textstyle \Xhat_{Q_1}}
,(7.8,6)+"u"+"r"*!C{\textstyle\bullet}*++!LD{\textstyle \Xhat_P}
,(7.8,6)+"r";(8.75,2)+"r"**@{}?="mid",?(1)**\crv~c{"mid"+(-.5,0)}?!{"btwo";"btwo"+(.866025404,.5)}="fout"?*++!L{\textstyle \Xhat_{Q_2}}
,"bone";"a"**@{-},"c"**@{-},"cout"*!C\dir{-},"d"**@{-},"dout"*!C\dir{-}
,"btwo"**@{-}
,"btwo";"e"**@{-},"eout"*!C\dir{-},"f"**@{-},"fout"*!C\dir{-},"b"**@{-}
\end{xy}
\end{equation*}
\caption{The Arthur-Langlands Partition}
\label{figArthurLanglands}
\end{figure}
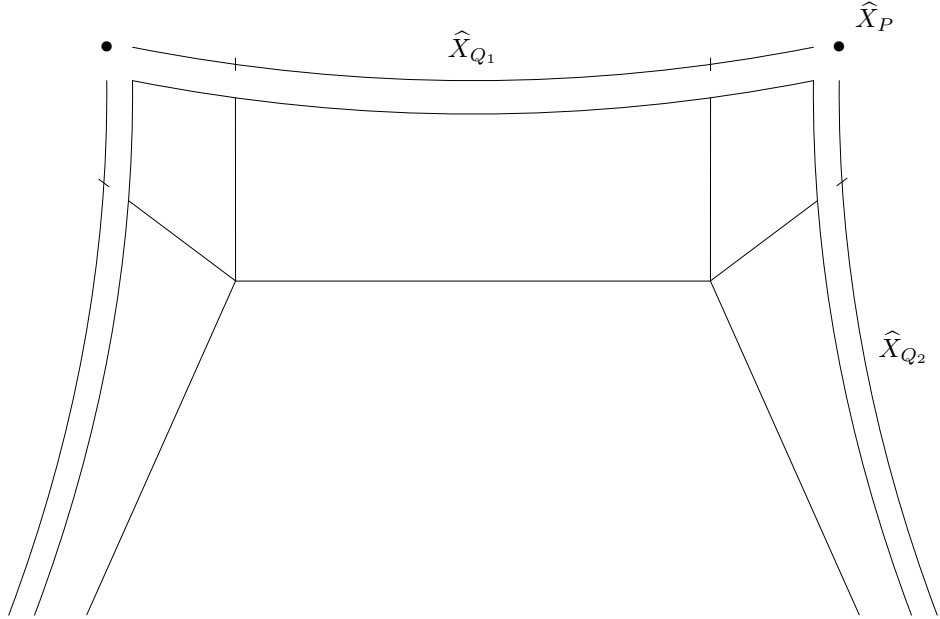

The solution is to replace each $X_P$ with a diffeomorphic copy which is
embedded as a domain in $X_P$ with compact closure.  If the diffeomorphic
embedding extended to the manifold with corners $\Xbar_P$ the new boundary
faces would now be at finite distance.  To do this naturally, consider the
Arthur-Langlands partition \cite{refnArthurTraceFormula} of $X$. This was
extended by the author to a partition of $\Xbar$ (and hence $\Xhat$) in
\cite{refnSaperTilings}; see Figure ~\ref{figArthurLanglands}.

\clearpage

The Arthur-Langlands partition divides $X$ into disjoint regions indexed by
$\Pl$; it depends on a certain parameter $\mathbf b =(b_P)_{P\in\Pl}$,
where $b_P\in \mathscr A_P^G$.  We will only need here the central region
$\Xbar_0$, which is the one indexed by $G\in\Pl$.  Recall that although the
decomposition $D= \mathscr A_P^G\times e_P$ does not descend to a
decomposition $\mathscr A_P^G\times Y_P$ of $X$, such a decomposition does
hold within an appropriate cylindrical set $W_P$ (see
\eqref{eqnCylindricalSet} in \S\ref{ssectBorelSerreThirdStep}).  We can use
this decomposition to describe $\Xbar_0$:
\begin{equation*}
\Xbar_0 \cap W_P = ((\lsp-A_P^G \cdot b_P)\times Y_P)\cap W_P.
\end{equation*}
In \cite{refnSaperTilings} it was shown that there exists a natural
piecewise analytic diffeomorphism (depending on another parameter) of
$\Xbar$ with $\Xbar_0$.  So the next step is to throw away everything
outside of $\Xbar_0$ and start over.

We endow the central region $\Xbar_0$ with the metric induced by restriction
from the locally symmetric metric on $X$.  Since $\Xbar_0$ is then a compact
Riemannian manifold with corners, a smooth form on $\Xbar_0$ or any of its
boundary faces will automatically have finite $L^2$-norm.

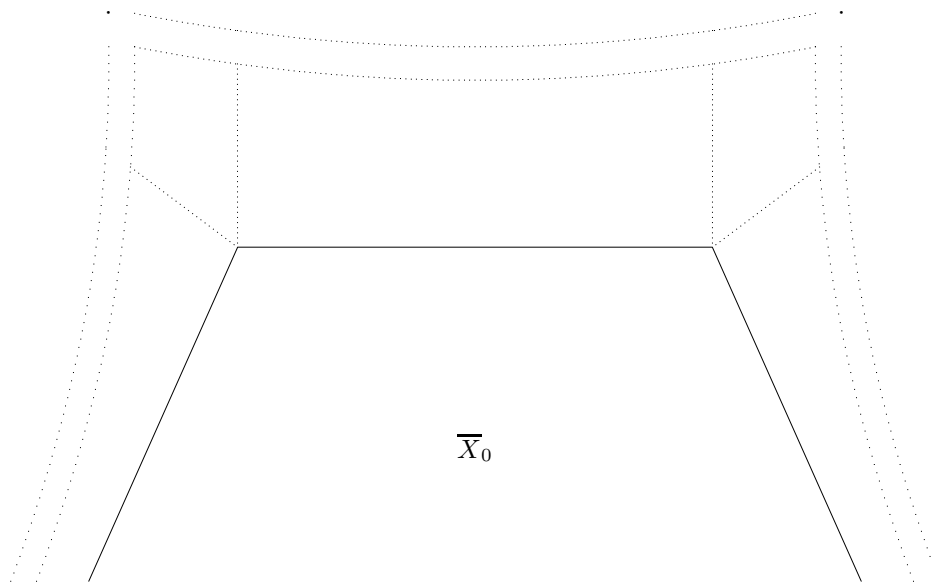
\begin{figure}[t]
\begin{equation*}
\begin{xy}
0;<.54in,0in>:<0in,.7in>::
(0,2)*i{\oplus},(9,2)*i{\oplus},(0,6)*i{\oplus},(9,6)*i{\oplus}
,(2.2,4.5)="bone",(6.8,4.5)="btwo"
,(.5,2);(8.5,2)**@{}
?!{"bone";"bone"+a(240)}="a"?!{"btwo";"btwo"+a(-60)}="b"
,(.25,2);(1.2,6)**@{}?="mid",?(1)
**\crv~c{~*=<1mm>{}~**@{.} "mid"+(.5,0)}?!{"bone";"bone"+(-.866025404,.5)}="c"
,(1.2,6);(7.8,6)**@{}?="mid",?(1)**\crv~c{~*=<1mm>{}~**@{.} "mid"+(0,-.5)}?!{"bone";"bone"+(0,1)}="d"?!{"btwo";"btwo"+(0,1)}="e"
,(7.8,6);(8.75,2)**@{}?="mid",?(1)**\crv~c{~*=<1mm>{}~**@{.} "mid"+(-.5,0)}?!{"btwo";"btwo"+(.866025404,.5)}="f"
,(-.25,0)="l",(0,.25)="u",(.25,0)="r"
,(.25,2)+"l";(1.2,6)+"l"**@{}?="mid",?(1)
**\crv~c{~*=<1mm>{}~**@{.} "mid"+(.5,0)}?!{"bone";"bone"+(-.866025404,.5)}="cout"
,(1.2,6)+"l"+"u"*!C{\textstyle\cdot}
,(1.2,6)+"u";(7.8,6)+"u"**@{}?="mid",?(1)**\crv~c{~*=<1mm>{}~**@{.}
  "mid"+(0,-.5)}?!{"bone";"bone"+(0,1)}="dout"?!{"btwo";"btwo"+(0,1)}="eout"?*++!D{}
,(7.8,6)+"u"+"r"*!C{\textstyle\cdot}*++!LD{}
,(7.8,6)+"r";(8.75,2)+"r"**@{}?="mid",?(1)**\crv~c{~*=<1mm>{}~**@{.}"mid"+(-.5,0)}?!{"btwo";"btwo"+(.866025404,.5)}="fout"?*++!L{}
,(4.5,3)*{\textstyle \Xbar_0}
,"bone";"a"**@{-},"c"**@{.},"cout"*!C\dir{.},"d"**@{.},"dout"*!C\dir{.}
,"btwo"**@{-}
,"btwo";"e"**@{.},"eout"*!C\dir{.},"f"**@{.},"fout"*!C\dir{.},"b"**@{-}
\end{xy}
\end{equation*}
\caption{The Central Region}
\label{figCentralRegion}
\end{figure}

\clearpage

Since $\Xbar_0$ is diffeomorphic with $\Xbar$ we can mimic the construction
of $\Xhat$ and $\Sheaf(\M)$ to obtain a realization $\Sheaf_0(\M)$ on
$\Xhat_0$ whose global sections also compute $H(\Xhat;\M)$; see Figure
~\ref{figCentralBoundaryStrata}.  Since we have removed the difficulty with
$L^2$-growth conditions, Hodge theory shows that classes in $H(\Xhat;\M)$
have harmonic representatives in $\Gamma(\Xhat_0, \Sheaf_0(\M))$.  This is
the analogue of the classical Hodge theorem for a compact Riemannian
manifold.

We now follow an argument from \cite{refnSaperSternTwo} however there is an
added combinatorial complication since our ``forms'' are actually
collections of forms and there are interaction terms in the differential.
The aim is to show that there exists $c>0$ such that the estimate
\begin{equation}
\|d\phi\|^2 + \|d^*\phi\|^2 \ge c \|\phi\|^2
\label{eqnEstimate}
\end{equation}
holds for all $\phi$ with degree $j\notin [c(\M),d(\M)]$.  This will
suffice to prove the theorem by (a generalization of) Proposition
~\ref{propEstimate}.

\begin{figure}[t]
\begin{equation*}
\begin{xy}
0;<.54in,0in>:<0in,.7in>::
(0,2)*i{\oplus},(9,2)*i{\oplus},(0,6)*i{\oplus},(9,6)*i{\oplus}
,(2.2,4.5)="bone",(6.8,4.5)="btwo"
,(.5,2);(8.5,2)**@{}
?!{"bone";"bone"+a(240)}="a"?!{"btwo";"btwo"+a(-60)}="b"
,"bone";"a"**@{-}
,"btwo"**@{-}
,"btwo";"b"**@{-}
,(-.216506351,.125)="l",(-.125,.25)="lu",(0,.25)="u",(.216506351,.125)="r",(.125,.25)="ru"
,"bone"+"l";"a"+"l"**@{-},
,"bone"+"lu"*!C{\textstyle\bullet}
,"bone"+"u";"btwo"+"u"**@{-}?*++!D{\textstyle \Xhat_{0Q_1}}
,"btwo"+"ru"*!C{\textstyle\bullet}*++!LD{\textstyle \Xhat_{0P}}
,"btwo"+"r";"b"+"r"**@{-}?*++!L{\textstyle \Xhat_{0Q_2}}
,(4.5,3)*{\textstyle \Xhat_0}
\end{xy}
\end{equation*}
\caption{The complex $\Gamma(\Xhat_0,\Sheaf_0(\M))$ lives on the disjoint
  union $\coprod_P \Xhat_{0P}$}
\label{figCentralBoundaryStrata}
\end{figure}

\clearpage

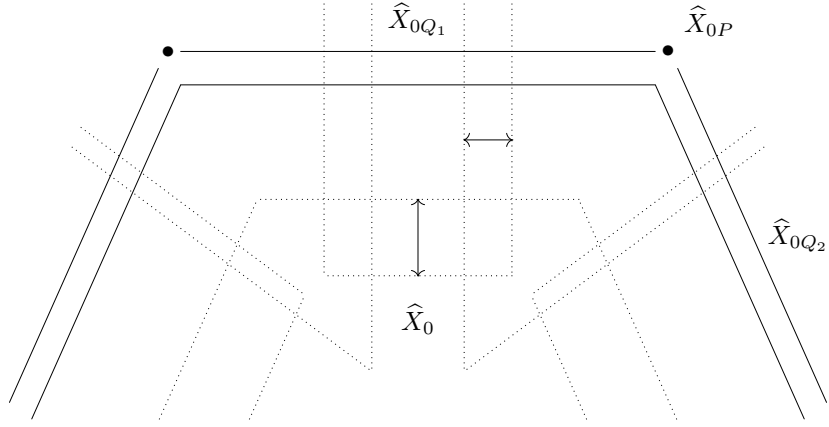
\begin{figure}[t]
\begin{equation*}
\begin{xy}
0;<.54in,0in>:<0in,.7in>::
(0,2)*i{\oplus},(9,2)*i{\oplus},(0,6)*i{\oplus},(9,6)*i{\oplus}
,(2.2,4.5)="bone",(6.8,4.5)="btwo"
,"bone"+/r1in/+/d1.5in/="bbone"
,"btwo"+/l1in/+/d1.5in/="bbtwo"
,"bone"+/r.75in/+/d1in/="bbbone"
,"btwo"+/l.75in/+/d1in/="bbbtwo"
,"bone"+/r.65in/+/d1.1in/="bbbbone"
,"btwo"+/l.65in/+/d1.1in/="bbbbtwo"
,"bone"+/r.4in/+/d.6in/="Gone"
,"btwo"+/l.4in/+/d.6in/="Gtwo"
,(.5,2);(8.5,2)**@{}
?!{"bone";"bone"+a(240)}="a"?!{"btwo";"btwo"+a(-60)}="b"
?!{"bbbbone";"bbbbone"+a(240)}="aaaa"?!{"bbbbtwo";"bbbbtwo"+a(-60)}="bbbb"
?!{"Gone";"Gone"+a(240)}="Ga"?!{"Gtwo";"Gtwo"+a(-60)}="Gb"
,"bone";"a"**@{-}
,"btwo"**@{-}
,"btwo";"b"**@{-}
,(-.216506351,.125)="l",(-.125,.25)="lu",(0,.25)="u",(.216506351,.125)="r",(.125,.25)="ru"
,"bone"+"l";"a"+"l"**@{-}?!{"bbone";"bbone"+(-.866025404,.5)}+/va(150).25in/="c"?!{"bbbbone";"bbbbone"+a(150)}+/va(150).25in/="cccc"
,"bone"+"lu"*!C{\textstyle\bullet}
,"bone"+"u";"btwo"+"u"**@{-}?*++!D{\textstyle \Xhat_{0Q_1}}?!{"bbone";"bbone"+(0,1)}+/u.25in/="d"?!{"bbbone";"bbbone"+(0,1)}+/u.25in/="ddd"?!{"bbtwo";"bbtwo"+(0,1)}+/u.25in/="e"?!{"bbbtwo";"bbbtwo"+(0,1)}+/u.25in/="eee"
,"btwo"+"ru"*!C{\textstyle\bullet}*++!LD{\textstyle \Xhat_{0P}}
,"btwo"+"r";"b"+"r"**@{-}?*++!L{\textstyle \Xhat_{0Q_2}}?!{"bbtwo";"bbtwo"+a(30)}+/va(30).25in/="f"?!{"bbbbtwo";"bbbbtwo"+a(30)}+/va(30).25in/="ffff"
,(4.5,2.75)*{\textstyle \Xhat_0}
,"bbone";"c"**@{.},"d"**@{.}
,"bbbone";"ddd"**@{.},"bbbtwo"**@{.}?(.5)="wGQb"
,"bbbtwo";"eee"**@{.}?(.5)="wQPb"
,"bbtwo";"e"**@{.}?!{"wQPb";"wQPb"+/l1in/}="wQPa","f"**@{.}
,"bbbbone";"cccc"**@{.},"aaaa"**@{.}
,"bbbbtwo";"bbbb"**@{.},"ffff"**@{.}
,"Gone";"Ga"**@{.},"Gtwo"**@{.}?(.5)="wGQa"
,"Gtwo";"Gb"**@{.}
,\ar@{<->} "wQPa";"wQPb"
,\ar@{<->} "wGQa";"wGQb"
\end{xy}
\end{equation*}
\caption{Covering by cylindrical sets $W_P$}
\label{figCylindricalCover}
\end{figure}

Let $\{W_P\}_{P\in\Pl}$ be a cylindrical cover of $\Xbar$ as in
\eqref{eqnCylindricalCover} of \S\ref{ssectBorelSerreThirdStep}; we can
arrange that for all $P\in \Pl$, the cylindrical set $W_P =
(\Abar_P^{G+}\cdot s_P) \times O_P$ projects to a relatively compact subset
of $X_{0P}$.  Such a cover induces a cover $\{W_P\}_{P\le Q}$ of each
$\Xhat_{0Q}$; for simplicity we denote the elements of such a cover again
by $W_P$.  The situation is pictured in Figure ~\ref{figCylindricalCover}.

Assume that for each $P\in \Pl$ we can establish an estimate
\begin{equation}
\|d\phi\|^2 + \|d^*\phi\|^2 \ge c_P \|\phi\|^2
\label{eqnEstimateP}
\end{equation}
for all $\phi$ with compact support in $\Xhat_0 \cap W_P$ and degree
$j\notin [c(\M),d(\M)]$.  Here $c_P>0$ is to be a constant independent of
$O_P$.  If this is possible, and if there exists a partition of unity
$\{\eta_P\}_{P\in\Pl}$ with $|d\eta_P|<\epsilon$ for all $P$ (where
$\epsilon>0$ depends on $\{c_P\}_{P\in\Pl}$), then the above estimates may
be patched together to yield \eqref{eqnEstimate}.  Such a partition of
unity always exists.  The point is that the magnitude $|d\eta_P|$ is
roughly inversely proportional to the ``width'' of the overlaps $W_P\cap
W_Q$ for all $Q\neq P$; two examples of these ``widths'', for $W_{Q_1}\cap
W_P$ and $W_G\cap W_{Q_1}$, are indicated in Figure
~\ref{figCylindricalCover}.  However it is possible to choose
$\{W_P\}_{P\in\Pl}$ so that these widths are arbitrarily large; this may
involve enlarging $\Xhat_0$ through manipulation of the parameter $\mathbf
b$ but there is plenty of room to do this since the boundary of $\Xhat$ is
infinitely far away.  Thus $\{\eta_P\}_{P\in\Pl}$ can be chosen with
arbitrarily small derivative.

It remains to restrict our attention to a single $W_P$ and prove
\eqref{eqnEstimateP}.

\clearpage

\begin{figure}[t]
\begin{equation*}
\begin{xy}
0;<.54in,0in>:<0in,.7in>::
(0,6)*i{\oplus},(9,6)*i{\oplus}
,(2.2,4.5)="bone",(6.8,4.5)="btwo"
,"bone"+/r1in/+/d1.5in/="bbone"
,"btwo"+/l1in/+/d1.5in/="bbtwo"
,"bone"+/r.75in/+/d1in/="bbbone"
,"btwo"+/l.75in/+/d1in/="bbbtwo"
,"bone"+/r.65in/+/d1.1in/="bbbbone"
,"btwo"+/l.65in/+/d1.1in/="bbbbtwo"
,(.5,2);(8.5,2)**@{}
?!{"bone";"bone"+a(240)}="a"?!{"btwo";"btwo"+a(-60)}="b"
?!{"bbbbone";"bbbbone"+a(240)}="aaaa"?!{"bbbbtwo";"bbbbtwo"+a(-60)}="bbbb"
,"bone";"btwo"**@{}?!{"bbtwo";"bbtwo"+(0,1)}="ein"
,"btwo";"b"**@{}?!{"bbtwo";"bbtwo"+a(30)}="fin"
,(-.216506351,.125)="l",(-.125,.25)="lu",(0,.25)="u",(.216506351,.125)="r",(.125,.25)="ru"
,"bone"+"u";"btwo"+"u"**@{}?!{"bbtwo";"bbtwo"+(0,1)}+/u.25in/="e"
,"btwo"+"r";"b"+"r"**@{}?!{"bbtwo";"bbtwo"+a(30)}+/va(30).25in/="f"
,"ein";"btwo"**@{-}
,"btwo";"fin"**@{-}
,(-.216506351,.125)="l",(-.125,.25)="lu",(0,.25)="u",(.216506351,.125)="r",(.125,.25)="ru"
,"ein"+"u";"btwo"+"u"**@{-}?(.5)*+!D{\textstyle \Xhat_{0Q_1}\cap W_P}
,"btwo"+"ru"*!C{\textstyle\bullet}*++!LD{\textstyle \Xhat_{0P}\cap W_P}
,"btwo"+"r";"fin"+"r"**@{-}?(.5)*+!LD{\textstyle \Xhat_{0Q_2}\cap W_P}
,"bbtwo";"btwo"**@{}?(.65)*++!C{\textstyle \Xhat_{0}\cap W_P}
,"bbtwo";"e"**@{.},"f"**@{.}
\end{xy}
\end{equation*}
\caption{A single cylindrical set $W_P$}
\label{figSingleCylindricalSet}
\end{figure}
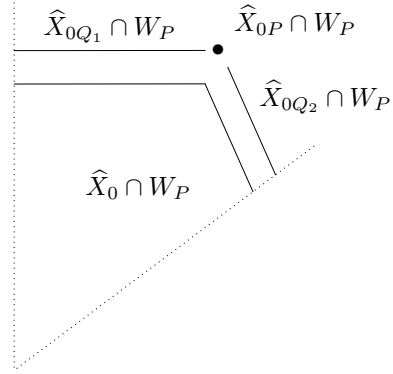

A section $\phi=(\phi_Q)_{Q\ge P}$ of $\Sheaf_0(\M)$ over $\Xhat_0\cap W_P$
is composed of $\EE_Q$-valued forms $\phi_Q$ which live on $\Xhat_{0Q}\cap
W_P$ as in Figure ~\ref{figSingleCylindricalSet} above.  We can apply
harmonic projection to each $\phi_Q$ along the nilmanifold fibers $\mathscr
N_P^{\prime Q}$ \cite{refnZuckerWarped} and obtain an $(\HH(\n_P^Q;E_Q)\otimes
\CC_{\r})$-valued form on $((\lsp-A_P^Q\cdot b_P)\times X_P)\cap W_P$; the
factor $\CC_\r$ accounts for the volume of the nilmanifold fibers.  However
since we are seeking to prove the estimate \eqref{eqnEstimateP} for
compactly supported forms, it certainly suffices to consider more generally
forms whose components have compact support anywhere in $(\lsp-A_P^Q\cdot
b_P)\times X_P$, which for simplicity we write as $\lsp-A_P^Q\times X_P$.
The result is that we have unfolded the situation of Figure
~\ref{figCentralBoundaryStrata}, which was a quotient by $\G$, to the
following quotient by $\G_P$:
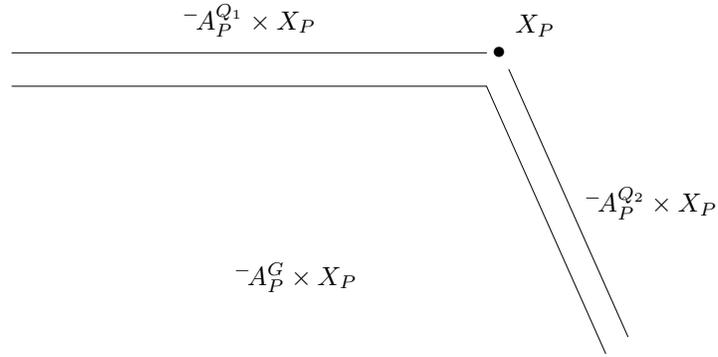
\begin{equation*}
\begin{xy}
0;<.54in,0in>:<0in,.7in>::
,(2.2,4.5)="bone",(6.8,4.5)="btwo"
,(.5,2.5);(8.5,2.5)**@{}
?!{"bone";"bone"+a(240)}="a"?!{"btwo";"btwo"+a(-60)}="b"
,"bone";"btwo"**@{-}
,"btwo";"b"**@{-}
,(-.216506351,.125)="l",(-.125,.25)="lu",(0,.25)="u",(.216506351,.125)="r",(.125,.25)="ru"
,"bone"+"u";"btwo"+"u"**@{-}?*++!D{\textstyle \lsp-A_P^{Q_1}\times X_P}
,"btwo"+"ru"*!C{\textstyle\bullet}*++!LD{\textstyle X_{P}}
,"btwo"+"r";"b"+"r"**@{-}?*++!L{\textstyle \lsp-A_P^{Q_2}\times
  X_P}
,"btwo"+/l1in/+/d1in/*{\textstyle \lsp-A_P^G\times X_P}
\end{xy}
\end{equation*}

\clearpage

We need to establish the estimate \eqref{eqnEstimateP} for compactly
supported forms in the following unfolded complex:
\begin{equation*}
\xymatrix @R=5.196cm @C=2.8cm @!0{
{\vbox{\hbox{\small$A_{(2)}(\lsp-A_P^{Q_1}\times X_P;$}
    \hbox{\small\hphantom{$A_{(2)}($}$\HH(\n_P^{Q_1};E_{Q_1})\otimes
      \CC_\r)$}}}
\ar@(u,l)[]!UC;[]!CL_{d_{Q_1Q_1}} 
\ar[rr]^-{d_{PQ_1}} &
{} &
{A_{(2)}(X_{P};\EE_P\otimes \CC_\r)} \ar@(r,u)[]!CR;[]!UC_{d_{PP}} \\
{} & {\vbox{\hbox{\small$A_{(2)}(\lsp-A_P^G\times X_P;$}
    \hbox{\small\hphantom{$A_{(2)}($}$\HH(\n_P;E_{G})\otimes \CC_\r)$}}}
\ar@(l,d)[]!CL;{[]!DC}_{d_{GG}} \ar[rr]^-{d_{Q_2G}}
\ar[ul]^-{d_{Q_1G}} \ar[ur]^-{d_{PG}} & {} &
   {\vbox{\hbox{\small$A_{(2)}(\lsp-A_P^{Q_2}\times X_P;$}
    \hbox{\small\hphantom{$A_{(2)}($}$\HH(\n_P^{Q_2};E_{Q_2})\otimes
      \CC_\r)$}}}
   \ar@(r,u)[]!CR;[]!UC_{d_{Q_2Q_2}} 
   \ar[ul]^-{d_{PQ_2}}
}
\end{equation*}
By Proposition ~\ref{propEstimate} it suffices to show that the cohomology
of this complex vanishes in degrees outside of $[c(\M),d(\M)]$ and that the
cohomology is Hausdorff in degree $c(\M)$.  It also suffices to do this
with the coefficients replaced by their $V\otimes \CC_\r$-isotypical
component for every irreducible representation $V$ of $L_P$.

Begin the calculation of the cohomology by forming the cohomology of each
summand (for simplicity we assume that $f_{QQ}=0$ for all $Q$, that is,
$\M$ is reduced as in \S\ref{ssectReducedLmodules}):
\begin{equation}
\vcenter{\hbox to 0pt{\hss
\xymatrix @R=5.196cm @C=2.8cm @!0{
{\vbox{\hbox{\small$H_{(2)}(\lsp-A_P^{Q_1}\times X_P;\VV\otimes
    \CC_\r)\mathop\otimes$}
  \hbox{\small$\qquad\Hom_{L_P}(V,H(\n_P^{Q_1};E_{Q_1}))$}}}
\ar[rr]^-{d_{PQ_1}} &
{} &
{\vbox{\hbox{\small$H_{(2)}(X_P;\VV\otimes \CC_\r)\mathop\otimes$}
   \hbox{\small$\qquad\Hom_{L_P}(V,E_{P})$}}}
\\
{} & {\vbox{\hbox{\small$ H_{(2)}(\lsp-A_P^G\times X_P;\VV\otimes
    \CC_\r)\mathop\otimes$}
    \hbox{\small$\qquad \Hom_{L_P}(V,H(\n_P;E_{G}))$}}}
\ar[ul]^-{d_{Q_1G}} \ar[ur]^-{d_{PG}} \ar[rr]^-{d_{Q_2G}} & {} &
{\vbox{\hbox{\small$ H_{(2)}(\lsp-A_P^{Q_2}\times X_P;\VV\otimes
      \CC_\r)\mathop\otimes$}
  \hbox{\small$\qquad\Hom_{L_P}(V,H(\n_P^{Q_2};E_{Q_2}))$}}}
  \ar[ul]^-{d_{PQ_2}}
}\hss}}
\label{eqnEOneTerm}
\end{equation}

Consider a typical $L^2$-cohomology group that occurs in
\eqref{eqnEOneTerm}, say
\begin{equation}
\label{eqnTypicalLTwoCohomology}
H_{(2)}(\lsp-A_P^{Q}\times X_P;\VV\otimes\CC_\r)
\end{equation}
where $Q\ge P$.  The locally constant sheaf $\VV\otimes\CC_\r$ is constant
along the $\lsp-A_P^Q$ factor however its metric is not; the coefficients
add an exponential weight of $a^{-2(\xi_V+\r)}$ to the $L^2$-norm integral.
This weight factors as
\begin{align*}
a^{-2(\xi_V+\r)} &= \prod_{\al \in \D_P^Q}(a^{-\b_\al})^{2(\xi_V+\r,\al)}\\
\intertext{corresponding to the product decomposition}
\lsp-A_P^Q &\cong [1,\infty)^{\D_P^Q},\qquad a \mapsto
(a^{-\b_{\al_1}},\dots,a^{-\b_{\al_r}}).
\end{align*}
However the weighted $L^2$-cohomology $H_{(2)}([1,\infty); s^k)$ (where $s
= a^{-\b_{\al}}$ and $[1,\infty)$ has the multiplicatively invariant metric
$ds^2/s^2$) vanishes for $k>0$ and is $\CC$ for $k<0$.  Thus whenever there
exists $\al \in \D_P^Q$ such that $(\xi_V+\r,\al)>0$, a K\"unneth-type
argument \cite{refnZuckerWarped} shows that
\eqref{eqnTypicalLTwoCohomology} vanishes.  Similarly if $(\xi_V+\r,\al)<0$
for all $\al$ we can drop the $\lsp-A_P^{Q}$ factor in
\eqref{eqnTypicalLTwoCohomology}.  The situation where $(\xi_V+\r,\al)=0$
for some $\al$ can be handled by a slightly more delicate argument.

In order to illustrate the remainder of the proof, suppose that $\xi_V$
satisfies
\begin{equation}
\label{eqnSampleV}
(\xi_V+\r,\al_1) < 0, \qquad (\xi_V+\r,\al_2) > 0,
\end{equation}
where $\D_P=\{\al_1,\al_2\}$ and $\D_P^{Q_i}=\{\al_i\}$ for $i=1,2$.  By
the above arguments, equation \eqref{eqnEOneTerm} becomes
\begin{equation*}
\xymatrix  @R=5.196cm @C=2.8cm @!0{
{\vbox{\hbox{\small$H_{(2)}(X_P;\VV)\mathop \otimes$}
  \hbox{\small$\qquad\Hom_{L_P}(V,H(\n_P^{Q_1};E_{Q_1}))$}}}
\ar[rr]^-{d_{PQ_1}} &
{} &
{\vbox{\hbox{\small$ H_{(2)}(X_P;\VV)\mathop\otimes$}
    \hbox{\small $\qquad\Hom_{L_P}(V,E_{P})$}}} \\
{} & *++{\txt{\small$0$}}
\ar[ul]^-{d_{Q_1G}} \ar[ur]^-{d_{PG}} \ar[rr]^-{d_{Q_2G}} & {} &
*++{\txt{\small$0$}} \ar[ul]^-{d_{PQ_2}}
}
\end{equation*}
However $Q_V=Q_V'=Q_1$ in this case so the total cohomology is (compare
\eqref{eqnExampleMicroSupportExpression} in \S\ref{ssectFunctorsLModules})
\begin{equation}
\label{eqnTotalCohomology}
H_{(2)}(X_P;\VV)\otimes\Hom_{L_P}(V,H(i_P^* \ihat_{Q_V}^! \M)).
\end{equation}
 
The vanishing theorem for $L^2$-cohomology, Theorem
~\ref{thmVanishingLTwoCohomology}, shows that the first factor of
\eqref{eqnTotalCohomology} vanishes unless $(V|_{M_P})^* \cong
\overline{V|_{M_P}}$ and the degree is in
\begin{equation*}
[\tfrac12(\dim D_P - \dim D_P(V)),  \tfrac12(\dim D_P + \dim D_P(V))];
\end{equation*}
the cohomology is also Hausdorff in degree $\tfrac12(\dim D_P - \dim
D_P(V))$.  On the other hand, if $(V|_{M_P})^* \cong \overline{V|_{M_P}}$
then the second factor of \eqref{eqnTotalCohomology} vanishes unless $V\in
\mS(\M)$ and the degree is in
\begin{equation*}
[c(V;\M),d(V;\M)].
\end{equation*}
The sum of these two degree ranges gives the desired result.

This proves the vanishing theorem for the cohomology of an $\L$-module.

\section{Micro-support of intersection cohomology}
When an $\L$-module is known explicitly the micro-support is not difficult
to compute; for example, consider our calculation of $\mS(\i_{G*}E)$ in
\S\ref{ssectExampleMicroSupport}.  Unfortunately $\IpC(\Xhat;E)$ is defined
inductively and we do not have an explicit closed formula for its local
cohomology.  Furthermore the conjugate-self-contragredient condition in the
definition of micro-support does not behave well under induction.
Nonetheless we can prove the following combinatorial theorem
\cite[Corollary ~17.2]{refnSaperLModules}:

\begin{thm}
\label{thmMicroSupportIC}
Let $E$ be an irreducible $G$-module and let $p=m$ or $n$ be a
middle-perversity.  Assume the $\QQ$-root system of $G$ does not contain a
factor of type $D_n$, $E_n$, or $F_4$.  If $(E|_{\lsp0\,G})^* \cong
\overline{E|_{\lsp0\,G}}$, then $\emS(\IpC(\Xhat;E)) = \{E\}$.
\end{thm}

Actually we can describe the entire micro-support, not just the essential
micro-support; see \cite[Theorem ~17.1]{refnSaperLModules} and
\eqref{eqnFundamentalP}, \eqref{eqnFundamentalG} below.  In fact
\cite{refnSaperLModules} treats the more general case in which
$(E|_{\lsp0\,G})^* \cong \overline{E|_{\lsp0\,G}}$ is not assumed.  We
expect that the condition on the $\QQ$-root system will be able to be
removed; the condition is satisfied in the Hermitian case or in the case of
a real equal-rank Satake compactification.

What does the theorem say?  Recall that an $L_P$-module $V$ belongs to
$\emS(\M)$ if and only if
\begin{align}
\tag{i}&(V|_{M_P})^* \cong \overline{V|_{M_P}},\text{ and}\\
\tag*{(ii)$'$} &H(i_P^* \ihat_{Q_V}^! \M)_V \longrightarrow H(i_P^*
\ihat_{Q'_V}^! \M)_V \text{ is nonzero,}
\end{align}
where we have defined $Q_V$ and $Q'_V$ by
\begin{align*}
\D_P^{Q_V} &=\{\,\al\in \D_P\mid (\xi_V+\r,\al)<0\,\},\\
\D_P^{Q_V'} &=\{\,\al\in \D_P\mid (\xi_V+\r,\al)\le 0\,\}.
\end{align*}
Since $\i_G^*\ihat_G^!\IpC(\Xhat;E) = E$, we certainly have $E\in
\emS(\IpC(\Xhat;E))$ if $(E|_{\lsp0\,G})^* \cong \overline{E|_{\lsp0\,G}}$.
The theorem is asserting that this is all: for every proper $P\in \Pl$ with
parabolic rank $r(P)=1,2,\dotsc$, there is no irreducible $L_P$-module
$V\in \emS(\IpC(\Xhat;E))$.

The only candidates for $V$ are the irreducible constituents of
$H(\n_P;E)$, namely $V_{w(\lambda+\r)-\r}$ for $w\in W_P$.  The assertion
of the theorem amounts to a subtle relationship between the combinatorics
of the Weyl element $w$, namely the minimal lengths in the cosets $W^Rw$
for all $R\ge P$, and the geometry of $w(\lambda+\r)$ \emph{vis \`a vis}
the roots $\al \in \D_P$.  Since these matters play an important role in
Goresky and MacPherson's topological trace formula
\cite{refnGoreskyMacPhersonTopologicalTraceFormula} and in the
representation of cohomology classes by Eisenstein series
\cite{refnSchwermerGeneric}, it is likely that Theorem
~\ref{thmMicroSupportIC} will have applications beyond the
Rapoport/Goresky-MacPherson conjecture.

The use of micro-support and of intersection cohomology is an efficient way
to encode this relationship; a description of the result without these
tools becomes unmanageable rapidly as the parabolic rank increases.
However it is not difficult for low parabolic rank.  In the following
subsections we will sketch proofs for the parabolic rank $1$ and $2$ cases
that will give some indication of the above relationship.  Unfortunately
these arguments do not generalize to higher parabolic rank; in particular,
the condition on the $\QQ$-root system begins to play a role with parabolic
rank $4$.  In the last subsection, we indicate some of the difficulties and
outline how the general proof in \cite{refnSaperLModules} works.

From now on, we assume that $E$ is irreducible with highest weight
$\lambda$.  We will use Kostant's theorem and the notation introduced in
\S\ref{ssectExampleMicroSupport}.

\subsection{Parabolic rank \boldmath$1$}
The simplest case is when $P$ is a maximal parabolic $\QQ$-subgroup.  Write
$\D_P = \{\al\}$.  For any $\L$-module $\M$, the expression $H(i_P^*
\ihat_Q^! \M)$ is either the local cohomology at $P$ or the local
cohomology supported at $P$ depending on whether $Q=G$ or $Q=P$
respectively.  Given the definition of intersection cohomology (see in
particular \eqref{eqnEPDefinitionClean} in
\S\ref{ssectTruncationFunctors}), the cohomology of the link is $H(\n_P;E)$
and we have
\begin{equation}
H(i_P^* \ihat_Q^! \IpC(\Xhat;E)) = \begin{cases}
  \t^{\le p(P)}H(\n_P;E) & \text{if $Q=G$,}\\
  \t^{> p(P)}H(\n_P;E)[-1] & \text{if $Q=P$.}
				   \end{cases}
\label{eqnPQCohomology}
\end{equation}
Recall that Kostant's theorem says that
\begin{equation}
\label{eqnKostant}
H(\n_P; E) = \bigoplus_{w\in W_P} H(\n_P;E)_w = \bigoplus_{w\in
W_P}V_{w(\lambda+\r)-\r}[-\l(w)].
\end{equation}

Thus fix $V= V_{w(\lambda+\r)-\r}$ for some $w\in W_P$.  Since $\xi_V + \r
= w(\lambda+\r)$, the values of $Q_V$ and $Q'_V$ are given by
\begin{itemize}
\item $(w(\lambda+\r),\al)<0 \quad\Longrightarrow\quad Q_V = Q'_V =G$,
\item $(w(\lambda+\r),\al)=0 \quad\Longrightarrow\quad Q_V=P,\ Q'_V=G$,
\item $(w(\lambda+\r),\al)>0 \quad\Longrightarrow\quad Q_V=Q'_V=P$.
\end{itemize}
On the other hand, it follows from \eqref{eqnPQCohomology} that the map in
\ref{eqnMicroSupportTwoPrime} is $0$ precisely in the following cases:
\begin{itemize}
\item when $\l(w) > p(P)$ and $Q_V=Q'_V=G$,
\item when $Q_V = P$, $ Q'_V = G$, and 
\item when $\l(w) \le p(P)$ and $Q_V=Q'_V=P$.
\end{itemize}
By \eqref{eqnMiddlePerversities} the value of $p(P)$ is
\begin{equation*}
p(P) = p(\dim \n_P + \#\D_P) = 
\begin{cases}
  \left\lfloor(\dim \n_P -1)/2\right\rfloor & p=m, \\
  \left\lfloor(\dim \n_P)/2\right\rfloor    & p=n.
\end{cases}
\end{equation*}

From the above facts, it is easy to see that the theorem in this case is
equivalent to following proposition, which is a generalization of an old
observation of Casselman \cite{refnCasselman} in the $\RRrank$ $1$ case;
our proof uses Lemma ~\ref{lemBasicLemma}.

\begin{prop}
\label{propQRankOne}
Let $E$ be an irreducible $G$-module with highest weight $\lambda$.  Let
$P$ be a maximal parabolic $\QQ$-subgroup with $\D_P=\{\al\}$ and let $w\in
W_P$.  Assume that $(V_{w(\lambda+\r)-\r}|_{M_P})^* \cong
\overline{V_{w(\lambda+\r)-\r}|_{M_P}}$.
\begin{enumerate}
\item\label{itemNegative} If $(w(\lambda+\r),\al)\le0$ then $\l(w) \ge
(\dim\n_P)/2$.
\item\label{itemPositive} If $(w(\lambda+\r),\al)\ge0$ then $\l(w) \le
(\dim\n_P)/2$.
\end{enumerate}
Furthermore assume that $(E|_{\lsp0\,G})^* \cong \overline{E|_{\lsp0\,G}}$.
Then
\begin{equation*}
\l(w) = (\dim \n_P)/2 \quad\Longrightarrow\quad (w(\lambda+\r),\al) =0
\text{ for all $\al\in\D_P$.}
\end{equation*}
\end{prop}
\begin{proof}
The first two parts are simply Lemma ~\ref{lemBasicLemma} with the $\dim
\n_P(V)$ term omitted.  For the final assertion, note that the proof of
Lemma ~\ref{lemBasicLemma} together with the hypothesis $\l(w) = (\dim
\n_P)/2$ shows that for each root $\g$ in $\n_{P\CC}$, exactly one of
$w^{-1}\g$ and $w^{-1}\bar\g = -w^{-1}(\theta \g)$ is negative, and that
there are no real roots.  The first statement implies that the positive
system $w\Phi^+$ is $\theta$-stable \cite[Lemma ~8.6]{refnSaperLModules}.
The second statement, together with the Kostant-Sugiura classification of
Cartan subalgebras \cite{refnKostantConjugacyCartan}, \cite{refnSugiura},
\cite{refnSugiuraCorrection}, implies that $\h$ is a fundamental Cartan
subalgebra for $\mathfrak g$.  Thus the operator $\t_G$ equals $\theta$ as
well, and hence the hypothesis $(E|_{\lsp0\,G})^* \cong
\overline{E|_{\lsp0\,G}}$ is equivalent to
$w\lambda|_{\h_{M_P,-1}+\sa_P^G}=0$.  (This is because $w\lambda$ is the
highest weight of $E$ with respect to $w\Phi^+$.)  Since
$w\r|_{\h_{M_P,-1}+\sa_P^G}=0$ due to $w\Phi^+$ being $\theta$-stable, we
find that $w(\lambda+\r)|_{\sa_P^G}=0$ as desired.
\end{proof}

\subsection{Bidegree in \boldmath$H(\n_P;E)$}
For $P$ with parabolic rank $\ge 2$ we need a refinement of Kostant's
theorem.  Recall that for $P\le Q$ we have an isomorphism
\begin{equation}
H(\n_P;E) \cong H(\n_P^Q;H(\n_Q;E));
\label{eqnBigradedNilpotentCohomology}
\end{equation}
in particular $H(\n_P;E)$ is bigraded given $Q\ge P$.  We write the
bidegree as $(\l^Q,\l_Q)$ and define \emph{truncation by bidegree}
functors $\t^{\l_Q\le n}$ and $\t^{\l_Q >n}$ on such bigraded modules.  How
is this structure reflected in the description of $H(\n_P;E)$ given by
Kostant's theorem \eqref{eqnKostant}?

Write 
\begin{equation*}
W=W^PW_P \qquad\text{and}\qquad W=W^QW_Q
\end{equation*}
where $W^P$ is the Weyl group of $L_P$ and $W_P$ is the set of minimal
length coset representatives of $W^P\back W$, and similarly for $Q$.  Since
$P/N_Q$ is a parabolic subgroup of $L_Q$, we also have a decomposition
\begin{equation*}
W^Q = W^P W_P^Q
\end{equation*}
and thus $W = W^P W_P^Q W_Q$.  One can show in fact that 
\begin{equation*}
W_P = W_P^Q W_Q
\end{equation*}
and hence for $w\in W_P$ we may write $w= w^Q w_Q$ corresponding to this
decomposition.  It is easy then to verify that
\begin{equation*}
H(\n_P;E)_w \cong H(\n_P^Q;H(\n_Q;E)_{w_Q})_{w^Q}.
\end{equation*}
So the bidegree of $H(\n_P;E)_w$ in \eqref{eqnBigradedNilpotentCohomology}
is $(\l^Q(w),\l_Q(w))\equiv (\l(w^Q),\l(w_Q))$.

Recall that $\l(w)$ may be interpreted as the number of positive
$\CC$-roots which $w^{-1}$ sends to negative roots.  For $w\in W_P$, all
such roots occur in $\n_{P\CC}$; the bidegree $\l_Q(w)$ is the number of
those roots occurring in $\n_{Q\CC}$ and the bidegree $\l^Q(w)$ is the
number of roots occurring in (a lift of) $\n_{P\CC}^Q$.  In Figure
~\ref{figUnipotentRadicals} this is illustrated for two examples in the
case $G=\GL_n(\RR)$, where positive roots may be identified with matrix
positions above the diagonal.  As one compares the various $\l_Q(w)$,
pictures such as these give helpful insight even for other groups, however
one must be aware of the implicit assumption in the pictures that the
Dynkin diagram of the root system is linear.
\begin{figure}
  \begin{center}
    \mbox{%
      \vbox{%
	\setbox\subfigbox=\hbox{%
	  $\left(\vcenter{
	    \xymatrix @=0pt@M=0pt@*=<1.5cm,1.5cm>{
	      {} & {\n_P^Q}  \ar@{-}[]!LD+0;[r]!LD+0 \ar@{-}!LD+0;{[]!LU}+0
	      \save []!LD*+!UR{\scriptstyle \al_1} \restore & {} \\
		    {}  & {} & {} \ar@{-}[]!LD+0;[]!RD+0 \ar@{-}[]!LD+0;[u]!LU+0 
		    \save []!U*{\n_Q} \restore 
		    \save []!LD*+!UR{\scriptstyle \al_2} \restore \\
			  {}  & {} & {}
	  }}
	  \right)$}%
	\dimen\subfigwd=\wd\subfigbox%
	\box\subfigbox%
	\vskip\baselineskip%
	\hbox to \dimen\subfigwd{\hfil\footnotesize(a) $\D_P=\{\al_1,\al_2\}$, $\D_P^Q=\{\al_1\}$\hfil}%
      }%
      \qquad
      \vbox{%
	\setbox\subfigbox=\hbox{
	  $\left(\vcenter{
	    \xymatrix @=0pt@M=0pt@*=<.75cm,.75cm>{
	      {} & {} \save []!LD*+!UR{\scriptstyle \al_1} \restore
	      \ar@{-}[]!LD+0;[r]!RD+0 \ar@{-}!LD+0;{[]!LU}+0 & {} & {} & {} \\ 
		  {} & {} & {\n_P^Q} \save []!LD*+!UR{\scriptstyle \al_2} \restore
		  \ar@{-}[]!LD+0;[]!RD+0 \ar@{-}[]!LD+0;[]!LU+0 & {} & {\n_Q} & {} \\ 
		      {} & {} & {} & {} \save []!LD*+!UR{\scriptstyle \al_3} \restore
		      \ar@{-}[]!LD+0;[rr]!RD+0 \ar@{-}!LD+0;{[u]!LU}+0 & {} & {} \\
			  {} & {} & {} & {} & {} \save []!LD*+!UR{\scriptstyle \al_4} \restore
			  \ar@{-}[]!LD+0;[]!RD+0 \ar@{-}[]!LD+0;[]!LU+0 & {\n_P^Q} \\ 
			      {} & {} & {} & {} & {} & {} \save []!LD*+!UR{\scriptstyle \al_5} \restore
			      \ar@{-}[]!LD+0;[]!RD+0 \ar@{-}[]!LD+0;[]!LU+0 \\
				  {}
	  }}
	  \right)$}%
	\dimen\subfigwd=\wd\subfigbox%
	\box\subfigbox%
	\vskip\baselineskip%
	\hbox{\vtop{\hbox to \dimen\subfigwd{\hfil\footnotesize%
	      (b) $\D_P=\{\al_1,\dotsc,\al_5\}$,\hfil}%
	    \hbox to \dimen\subfigwd{\hfil\footnotesize$\D_P^Q=
	      \{\al_2,\al_4,\al_5\}$\hfil}%
	}}%
      }%
}%
\end{center}
\caption{Unipotent radicals for various $P\le Q$ in $\GL_n(\RR)$}
\label{figUnipotentRadicals}
\end{figure}

\subsection{Parabolic rank \boldmath$2$}
Now assume $r(P)=2$ and write $\D_P=\{\al_1,\al_2\}$ with
$\D_P^{Q_i}=\{\al_i\}$ for $i=1$, $2$.  The intersection cohomology
$\L$-module is illustrated in Figure
~\ref{figIntersectionCohomologyPrankTwo} of
\S\ref{ssectionIntersectionCohomology}.  In order to find a complex
computing the link intersection cohomology at $P$, we start with
$H(\n_P;E)$ and truncate it via a mapping cone in bidegree $\l_{Q_1}>
p(Q_1)$ and independently in bidegree $\l_{Q_2} > p(Q_2)$.  This yields the
complex \eqref{eqnLinkIntersectionCohomology}.  We find that $H(\n_P;E)_w$
contributes a nonvanishing component to the link cohomology if neither
truncation applies or if both truncations apply (in which case its degree
is shifted up by $1$):
\begin{equation}
\label{eqnLinkICPRankTwo}
  \t^{\l_{Q_1}\le p(Q_1)}\t^{\l_{Q_2} \le p(Q_2)}H(\n_P;E) \oplus 
  \t^{\l_{Q_1}> p(Q_1)}\t^{\l_{Q_2} > p(Q_2)}H(\n_P;E)[-1].
\end{equation}
See \S\ref{ssectGraphicalCohomology} below for a further discussion.

This enables us to calculate $H(i_P^* \ihat_Q^!  \IpC(\Xhat;E))$ if $Q=G$
or $Q=P$.  For $Q=G$ we obtain the local intersection cohomology at $P$,
\begin{equation*}
\t^{\le p(P)}\bigl(\t^{\l_{Q_1}\le p(Q_1)}\t^{\l_{Q_2} \le p(Q_2)}H(\n_P;E)
    \oplus \t^{\l_{Q_1}> p(Q_1)}\t^{\l_{Q_2} > p(Q_2)}H(\n_P;E)[-1]\bigr).
\end{equation*}
For $Q=P$ we obtain the local intersection cohomology supported at $P$,
\begin{equation*}
 \t^{> p(P)}\bigl(\t^{\l_{Q_1}\le p(Q_1)}\t^{\l_{Q_2} \le p(Q_2)}H(\n_P;E)
    \oplus \t^{\l_{Q_1}> p(Q_1)}\t^{\l_{Q_2} >
    p(Q_2)}H(\n_P;E)[-1]\bigr)[-1].
\end{equation*}

In order to calculate micro-support, we need to gather some calculations
and a generalization of Proposition ~\ref{propQRankOne}.  First, given that
we are working with a middle perversity, we have
\begin{equation*}
p(P) = \begin{cases}
  \left\lfloor(\dim \n_P)/2\right\rfloor & p=m, \\
  \left\lfloor(\dim \n_P + 1)/2\right\rfloor    & p=n,
\end{cases}
\quad\text{and}\quad
p(Q_i) = \begin{cases}
  \left\lfloor(\dim \n_{Q_i} -1)/2\right\rfloor & p=m, \\
  \left\lfloor(\dim \n_{Q_i})/2\right\rfloor    & p=n.
\end{cases}
\end{equation*}
Next it is helpful to note that if $k$ and $l$ are integers, then
\begin{equation*}
l \le   \left\lfloor k/2\right\rfloor \quad\Longleftrightarrow\quad
l \le    k/2 \qquad\text{and}\qquad
l >   \left\lfloor k/2\right\rfloor \quad\Longleftrightarrow\quad
l >    k/2;
\end{equation*}
the analogous formulas for $<$ and $\ge$ do not hold.  Finally the
generalization of Proposition ~\ref{propQRankOne} is:
\begin{prop}
\label{propArbPrank}
Let $E$ be an irreducible $G$-module with highest weight $\lambda$.  Let
$P$ be a parabolic $\QQ$-subgroup and let $w\in W_P$.  Assume that
$(V_{w(\lambda+\r)-\r}|_{M_P})^* \cong
\overline{V_{w(\lambda+\r)-\r}|_{M_P}}$.
\begin{enumerate}
\item\label{itemNegativeGeneral} If $(w(\lambda+\r),\al) \le0$ for all
  $\al\in \D_P$, then $\l_Q(w) \ge (\dim\n_Q)/2$ for all $Q\ge P$.
\item\label{itemPositiveGeneral} If $(w(\lambda+\r),\al) \ge0$ for all
  $\al\in \D_P$, then $\l_Q(w) \le (\dim\n_Q)/2$ for all $Q\ge P$.
\end{enumerate}
Furthermore assume that $(E|_{\lsp0\,G})^* \cong \overline{E|_{\lsp0\,G}}$
and that the hypotheses of \itemref{itemNegativeGeneral} or
\itemref{itemPositiveGeneral} hold.  Then
\begin{equation*}
\l(w) = (\dim \n_P)/2 \quad\Longrightarrow\quad (w(\lambda+\r),\al) =0
\text{ for all $\al\in\D_P$.}
\end{equation*}
\end{prop}
\begin{proof}[Sketch of proof]
For every $A_P$-weight $\al$ in $\n_Q$ (not necessarily simple), the
argument we gave for Lemma ~\ref{lemBasicLemma} may be applied to the
$\CC$-roots $\g$ for \pagebreak which $\g|_{\sa_P} = \al$.  The result is an estimate
on the contribution to $\l_Q(w)$ from such roots.  The sum of these
estimates over all $A_P$-weights in $\n_Q$ yields the proposition.  The
last assertion follows as in Proposition ~\ref{propQRankOne}.
\end{proof}

Consider $V=V_{w(\lambda+\r)-\r}$ for some $w\in W_P$ and assume that
$(V|_{M_P})^* \cong \overline{V|_{M_P}}$.  In the case $Q_V=P$, we have
$(w(\lambda+\r),\al)\ge 0$ for all $\al\in \D_P$.  It follows from the
proposition and the preceding calculations that $H(i_P^* \ihat_{Q_V}^!
\IpC(\Xhat;E))_V$ is $0$ unless $\l(w)=(\dim \n_P)/2$, in which case
$(w(\lambda+\r),\al) =0$ for all $\al\in \D_P$ and thus $Q'_V=G$.
Conversely if $Q'_V=G$ then $(w(\lambda+\r),\al)\le 0$ for all $\al\in
\D_P$ and it follows that either $H(i_P^* \ihat_{Q'_V}^!  \IpC(\Xhat;E))_V
= 0$ or $\l(w)=(\dim \n_P)/2$ and $Q_V=P$.

Thus for any $V$ which satisfies $Q_V$ or $Q'_V$ equal to either $P$ or $G$
and which can conceivably contribute to micro-support, we have $Q_V=P$,
$Q'_V=G$, $\l(w)=(\dim \n_P)/2$, and
\begin{align}
  H(i_P^* \ihat_{Q_V}^! \IpC(\Xhat;E))_V &= 
  \begin{cases}
    V[-\l(w)-2] & \text{if $p=m$,} \\
    0           & \text{if $p=n$.}
  \end{cases} \label{eqnFundamentalP}\\
  H(i_P^* \ihat_{Q'_V}^! \IpC(\Xhat;E))_V &= 
  \begin{cases}
    0         & \text{if $p=m$,} \\
    V[-\l(w)] & \text{if $p=n$.}
  \end{cases} \label{eqnFundamentalG}
\end{align}
Such $V$ belong to $\mS(\IpC(\Xhat;E))$ but obviously do not contribute to
the essential micro-support; we call them \emph{fundamental} since they
occur if and only if $P$ contains a fundamental parabolic $\RR$-subgroup of
$G$ \cite{refnBorelCasselman}, \cite[Lemma ~8.8]{refnSaperLModules}.
(Recall that a parabolic $\RR$-subgroup is \emph{fundamental} if a
fundamental Cartan subalgebra for its Levi quotient lifts to a fundamental
Cartan subalgebra for $G$.)

Unlike the parabolic rank $1$ case, there are two more cases of $H(i_P^*
\ihat_Q^!  \IpC(\Xhat;E))$ to consider, namely when $Q=Q_1$ or $Q=Q_2$.
Since the two cases are identical after relabeling, we focus on $Q=Q_1$.
Use the complex $\ihat_P^*\ihat_Q^!\M$ of
\eqref{eqnExampleMicroSupportExpression} in \S\ref{ssectFunctorsLModules}
applied to the $\L$-module $\M=\IpC(\Xhat;E)$ pictured in Figure
~\ref{figIntersectionCohomologyPrankTwo}.  We compute that
\begin{equation*}
\begin{split}
H(i_P^* \ihat_{Q_1}^! \IpC(\Xhat;E)) =  &\t^{\l_{Q_1}> p(Q_1)}\t^{\l_{Q_2}
  \le p(Q_2)}H(\n_P;E)[-1] \oplus \\
&\t^{\le p(P)}\bigl(\t^{\l_{Q_1}> p(Q_1)}\t^{\l_{Q_2} > p(Q_2)}
  H(\n_P;E)[-1]\bigr) \oplus \\
&\t^{> p(P)}\t^{\l_{Q_1}\le p(Q_1)}\t^{\l_{Q_2} \le p(Q_2)}
  H(\n_P;E)[-1].
\end{split}
\end{equation*}
The first term are those classes in $\t^{\l_{Q_1} > p(Q_1)}H(\n_P;E)[-1]$
which map to $0$ in the link cohomology, the next term are those classes
which map to nonzero elements in the link cohomology which are not
truncated at $P$, and the third term are those link cohomology classes
being truncated at $P$ which do not come from $\t^{\l_{Q_1} >
p(Q_1)}H(\n_P;E)[-1]$.

We can now handle $V$ such that $Q_V=Q'_V=Q_1$.  In this case
\begin{equation*}
H(i_P^* \ihat_{Q_1}^! \IpC(\Xhat;E))_V = 0
\end{equation*}
by the proposition below (which is actually slightly stronger than
necessary).  Consequently $V$ cannot be in $\mS(\IpC(\Xhat;E))$, much less
the essential micro-support.  The case where $Q_V=Q'_V=Q_2$ is treated
similarly.
\begin{prop}
\label{propQRankTwo}
Let $E$ be an irreducible $G$-module with highest weight $\lambda$.  Let
$P$ be a parabolic $\QQ$-subgroup with $\D_P=\{\al_1,\al_2\}$ and let $w\in
W_P$.  Assume that $(V_{w(\lambda+\r)-\r}|_{M_P})^* \cong
\overline{V_{w(\lambda+\r)-\r}|_{M_P}}$ and that
\begin{equation*}
(w(\lambda+\r),\al_1)\le 0 \quad\text{and}\quad (w(\lambda+\r),\al_2)\ge0.
\end{equation*}
\begin{enumerate}
\item\label{itemQOne} If $\l_{Q_1}(w) \ge (\dim \n_{Q_1})/2$ then $\l(w)
\ge (\dim \n_P)/2$.
\item\label{itemQTwo} If $\l_{Q_2}(w) \le (\dim \n_{Q_2})/2$ then $\l(w)
\le (\dim \n_P)/2$.
\end{enumerate}
If either hypothesis is a strict inequality, the corresponding conclusion
is also strict.
\end{prop}
\begin{proof}[Sketch of proof]
Apply Proposition ~\ref{propQRankOne} to the group $L_{Q_1}$ with the
representation $E_{Q_1}= H^{\l_{Q_1}(w)}(\n_{Q_1};E)_{w_{Q_1}}$, the
maximal parabolic subgroup $P/N_{Q_1}$, and the Weyl element $w^{Q_1}\in
W_P^{Q_1}$.  If $\u$ is the highest weight of $E_{Q_1}$, the representation
$H^{\l^{Q_1}(w)}(\n_P^{Q_1};E_{Q_1})_{w^{Q_1}} = V_{w^{Q_1}(\u + \r)-\r}$
is identical with $V_{w(\lambda+\r)-\r}$ and hence satisfies the conjugate
self-contragredient condition.  It follows from the proposition that
$\l^{Q_1}(w) \ge (\dim \n_P^{Q_1})/2$ which proves \itemref{itemQOne} since
$\l(w) = \l_{Q_1}(w) + \l^{Q_1}(w)$.  Part \itemref{itemQTwo} follows
similarly.
\end{proof}

\begin{rem}
Whereas the conjugate self-contragredience of $H(\n_P;E)_w$ trivially
implied the conjugate self-contragredience of
$H(\n_P^{Q_1};E_{Q_1})_{w^{Q_1}}$ in the proof above, it does not in
general imply the conjugate self-contragredience of
$H(\n_{Q_1};E)_{w_{Q_1}}$.  As we will indicate in
\S\ref{ssectPRankGEThree} below, this is a fundamental difficulty in
formulating a proof that applies to arbitrary parabolic rank.  The natural
approach would be to obtain information regarding the cohomology of the
link by applying the theorem inductively to larger parabolic subgroups;
this fails since the self-contragredience condition is not preserved under
the induction which replaces $w$ by $w_R$ for $R\ge P$.
\end{rem}

\subsection{Geometric interpretation  of link cohomology}
\label{ssectGraphicalCohomology}
Although the calculation of link intersection cohomology preceding
\eqref{eqnLinkICPRankTwo} is trivial, it is helpful when considering higher
parabolic rank to view it geometrically.  Thus we digress to indicate this
point of view.  Recall that $|\D_P|$ denotes the closed $(r(P)-1)$-simplex
with vertices $\al\in \D_P$ ; this is a stratified space with strata
$|\D_P^Q|^\circ$ indexed by $Q>P$.  For $w\in W_P$, define the
\emph{$w$-shifted perversity}
\begin{equation*}
p_w(Q) \equiv p(Q) - \l_Q(w) = p(\dim \n_Q + \#\D_Q) - \l_Q(w), \quad
Q\ge P;
\end{equation*}
we view this as a generalized perversity for the space $|\D_P|$ and so the
intersection cohomology $I_{p_w}H(|\D_P|;\ZZ)$ is defined.  The
$H(\n_P;E)_w$-isotypical component of the complex
$\i_P^*j_{P*}^{\vphantom{*}}j_P^* \IpC(\Xhat;E)$ computing link
intersection cohomology can be expressed as
\begin{equation*}
H(\n_P;E)_w \otimes I_{p_w}C(|\D_P|;\ZZ);
\end{equation*}
we take this as defining $I_{p_w}C(|\D_P|;\ZZ)$.  The corresponding
isotypical component of the link intersection cohomology is
\begin{equation}
\label{eqnLinkIHIsotypicalComponent}
H(\n_P;E)_w \otimes I_{p_w}H(|\D_P|;\ZZ).
\end{equation}

The isotypical components of $H(\i_P^*\IpC(\Xhat;E))$, the local
intersection cohomology, are expressed similarly using
$I_{p_w}H(c(|\D_P|);\ZZ)$.  Here $c(|\D_P|)$ is the cone on $|\D_P|$, a
stratified space with strata $c(|\D_P^Q|)^\circ$ indexed by $Q\ge P$.  (The
stratum associated to $Q=P$ is the vertex of the cone.)

For example, in the case of parabolic rank two, $|\D_P|$ is a closed
$1$-simplex with vertices $ \{\al_1,\al_2\}$.  Let $B_w\subseteq
\{\al_1,\al_2\}$ denote the set of vertices where a degree $0$ class on the
interior $|\D_P|^\circ$ would be truncated according to $p_w$, that is,
\begin{equation*}
B_w \equiv
\begin{cases}
  \emptyset & \text{if $\l_{Q_1}(w) \le p(Q_1)$,  $\l_{Q_2}(w)\le
    p(Q_2)$,} \\
  \{\al_1\} & \text{if $\l_{Q_1}(w) > p(Q_1)$,  $\l_{Q_2}(w)\le
    p(Q_2)$,} \\ 
  \{\al_2\} & \text{if $\l_{Q_1}(w) \le p(Q_1)$,  $\l_{Q_2}(w) >
    p(Q_2)$,} \\
  \{\al_1,\al_2\} & \text{if $\l_{Q_1}(w) > p(Q_1)$,  $\l_{Q_2}(w) >
    p(Q_2)$.}
\end{cases}
\end{equation*}
Then $I_{p_w}C(|\D_P|;\ZZ)$ is the complex $\ZZ \to \ZZ^{\#B_w}$ (compare
\eqref{eqnLinkIntersectionCohomology} in
\S\ref{ssectionIntersectionCohomology}).

We graphically represent $I_{p_w}C(|\D_P|;\ZZ)$ in Figure
~\ref{figRelativeCohomologyPRankTwo} along with the value of
$I_{p_w}H(|\D_P|;\ZZ)$ by placing a dot near $\al_i$ if truncation is
taking place there.
\newcommand{\invrootthree}{.577350269}%
\newdimen\sidelength%
\newdimen\gap%
\sidelength=.9cm%
\gap=.2cm%
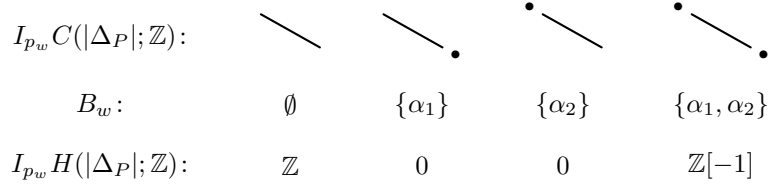
\begin{figure}%
\begin{equation*}%
\xymatrix @R-.6cm{
*{I_{p_w}C(|\D_P|;\ZZ)\colon} &
*!{\xybox{/va(-30).5\sidelength/;/va(150).5\sidelength/**[|(2)]@{-}}} &
*!{\xybox{/va(-30).5\sidelength/;/va(150).5\sidelength/**[|(2)]@{-},
    p-/\gap/*{\scriptstyle\bullet}}} &
*!{\xybox{/va(-30).5\sidelength/;/va(150).5\sidelength/**[|(2)]@{-},
    c+/\gap/*{\scriptstyle\bullet}}} &
*!{\xybox{/va(-30).5\sidelength/;/va(150).5\sidelength/**[|(2)]@{-},
    c+/\gap/*{\scriptstyle\bullet},p-/\gap/*{\scriptstyle\bullet}}} 
\\
*{B_w\colon} &
{\emptyset} & {\{\al_1\}} &
{\{\al_2\}} & {\{\al_1,\al_2\}} \\
*{I_{p_w}H(|\D_P|;\ZZ)\colon} & {\ZZ} & {0} & {0} & {\ZZ[-1]}
}
\end{equation*}%
\caption{The intersection cohomology $I_{p_w}H(|\D_P|;\ZZ)$ in parabolic
rank $2$}%
\label{figRelativeCohomologyPRankTwo}%
\end{figure}%
The corresponding diagrams for parabolic rank $3$ are given in Figure
~\ref{figRelativeCohomologyPRankThree}; here $|\D_P|$ is a $2$-simplex and
we draw a line near an edge if truncation of the degree $0$ class on the
interior is occurring, and we draw a dot near a vertex if the degree $0$ or
$1$ class on its link is being truncated.
\begin{figure}%
\begin{equation*}%
\xymatrix @R-.8cm{ *{I_{p_w}C(|\D_P|;\ZZ)\colon} &
  *!{\xybox{\xypolygon3"T"{~:{<\invrootthree\sidelength,0cm>:}
	~><{@{-}@*{[|(2)]}}}}} &
  *!{\xybox{{\xypolygon3"T"{~:{<\invrootthree\sidelength,0cm>:}
	  ~><{@{-}@*{[|(2)]}}}}, "T1"+/va(90)2\gap/*{\scriptstyle\bullet}}} &
  *!{\xybox{{\xypolygon3"T"{~:{<\invrootthree\sidelength,0cm>:}
	  ~><{@{-}@*{[|(2)]}}}}, "T1"+/va(90)2\gap/*{\scriptstyle\bullet},
      "T2"+/va(210)2\gap/*{\scriptstyle\bullet}}} &
  *!{\xybox{{\xypolygon3"T"{~:{<\invrootthree\sidelength,0cm>:}
	  ~><{@{-}@*{[|(2)]}}}}, "T1"+/va(90)2\gap/*{\scriptstyle\bullet},
      "T2"+/va(210)2\gap/*{\scriptstyle\bullet}, "T3"+/va(330)2\gap/*{\scriptstyle\bullet}}} \\
  *{I_{p_w}H(|\D_P|;\ZZ)\colon} & {\ZZ} & {0} & {\ZZ[-1]} & {\ZZ^2[-1]} \\
  *{\vbox to 1.25cm{}} \\ *{I_{p_w}C(|\D_P|;\ZZ)\colon} &
  *!{\xybox{{\xypolygon3"T"{~:{<\invrootthree\sidelength,0cm>:}
	  ~><{@{-}@*{[|(2)]}}}},
      "T2"+/va(150)\gap/;"T1"+/va(150)\gap/**[|(2)]@{-}}} &
  *!{\xybox{{\xypolygon3"T"{~:{<\invrootthree\sidelength,0cm>:}
	  ~><{@{-}@*{[|(2)]}}}},
      "T2"+/va(150)\gap/;"T1"+/va(150)\gap/**[|(2)]@{-},
      "T3"+/va(330)2\gap/*{\scriptstyle\bullet}}} &
  *!{\xybox{{\xypolygon3"T"{~:{<\invrootthree\sidelength,0cm>:}
	  ~><{@{-}@*{[|(2)]}}}},
      "T1"+/va(30)\gap/;"T3"+/va(30)\gap/**[|(2)]@{-},
      "T2"+/va(150)\gap/;"T1"+/va(150)\gap/**[|(2)]@{-}}} &
  *!{\xybox{{\xypolygon3"T"{~:{<\invrootthree\sidelength,0cm>:}
	  ~><{@{-}@*{[|(2)]}}}},
      "T1"+/va(30)\gap/;"T3"+/va(30)\gap/**[|(2)]@{-},
      "T2"+/va(150)\gap/;"T1"+/va(150)\gap/**[|(2)]@{-},
      "T1"+/va(90)2\gap/*{\scriptstyle\bullet}}} \\ *{I_{p_w}H(|\D_P|;\ZZ)\colon} & {0}
  & {\ZZ[-1]} & {\ZZ[-1]} & {0} \\ *{\vbox to 1.25cm{}} \\
  *{I_{p_w}C(|\D_P|;\ZZ)\colon} &
  *!{\xybox{{\xypolygon3"T"{~:{<\invrootthree\sidelength,0cm>:}
	  ~><{@{-}@*{[|(2)]}}}},
      "T1"+/va(30)\gap/;"T3"+/va(30)\gap/**[|(2)]@{-},
      "T2"+/va(150)\gap/;"T1"+/va(150)\gap/**[|(2)]@{-},
      "T3"+/va(270)\gap/;"T2"+/va(270)\gap/**[|(2)]@{-}}} &
  *!{\xybox{{\xypolygon3"T"{~:{<\invrootthree\sidelength,0cm>:}
	  ~><{@{-}@*{[|(2)]}}}},
      "T1"+/va(30)\gap/;"T3"+/va(30)\gap/**[|(2)]@{-},
      "T2"+/va(150)\gap/;"T1"+/va(150)\gap/**[|(2)]@{-},
      "T3"+/va(270)\gap/;"T2"+/va(270)\gap/**[|(2)]@{-},
      "T1"+/va(90)2\gap/*{\scriptstyle\bullet}}} s&
  *!{\xybox{{\xypolygon3"T"{~:{<\invrootthree\sidelength,0cm>:}
	  ~><{@{-}@*{[|(2)]}}}},
      "T1"+/va(30)\gap/;"T3"+/va(30)\gap/**[|(2)]@{-},
      "T2"+/va(150)\gap/;"T1"+/va(150)\gap/**[|(2)]@{-},
      "T3"+/va(270)\gap/;"T2"+/va(270)\gap/**[|(2)]@{-},
      "T1"+/va(90)2\gap/*{\scriptstyle\bullet}, "T2"+/va(210)2\gap/*{\scriptstyle\bullet}}} &
  *!{\xybox{{\xypolygon3"T"{~:{<\invrootthree\sidelength,0cm>:}
	  ~><{@{-}@*{[|(2)]}}}},
      "T1"+/va(30)\gap/;"T3"+/va(30)\gap/**[|(2)]@{-},
      "T2"+/va(150)\gap/;"T1"+/va(150)\gap/**[|(2)]@{-},
      "T3"+/va(270)\gap/;"T2"+/va(270)\gap/**[|(2)]@{-},
      "T1"+/va(90)2\gap/*{\scriptstyle\bullet}, "T2"+/va(210)2\gap/*{\scriptstyle\bullet},
      "T3"+/va(330)2\gap/*{\scriptstyle\bullet}}} \\ *{I_{p_w}H(|\D_P|;\ZZ)\colon} &
  {\ZZ^2[-1]} & {\ZZ[-1]} & {0} & {\ZZ[-2]} }
\end{equation*}%
\caption{The intersection cohomology $I_{p_w}H(|\D_P|;\ZZ)$ in parabolic
rank $3$.  (Configurations differing from previous ones by a rotation are
omitted.)}%
\label{figRelativeCohomologyPRankThree}%
\end{figure}
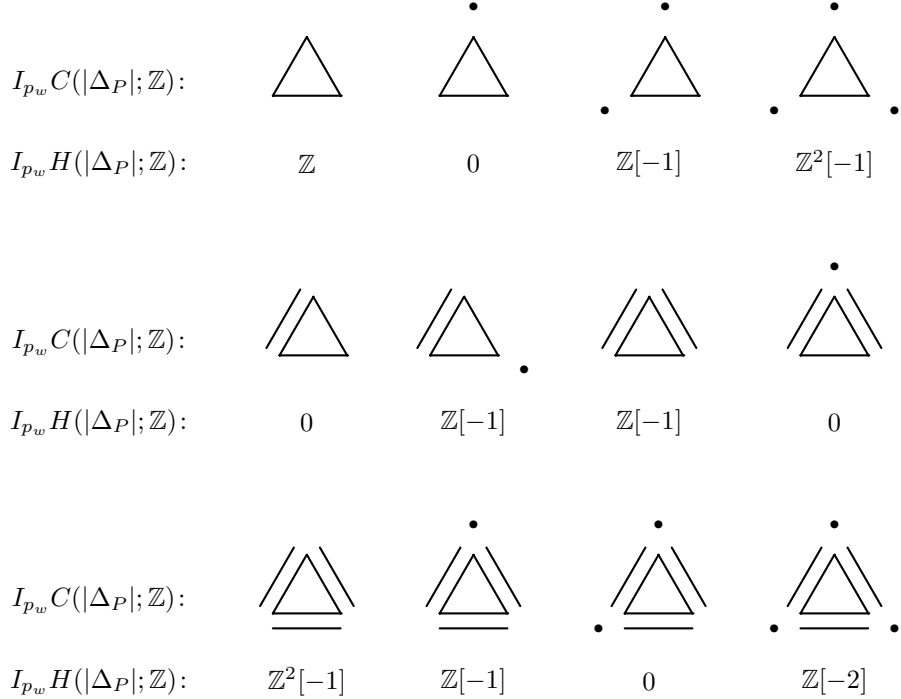%

Properly speaking, we should decorate the dots and lines in these figures
with the degree $p_w(Q)$ above which we truncate.  However the cohomology
potentially being truncated occurs only in one degree (namely $0$ or $1$),
so simply indicating whether or not it is truncated is sufficient.  The
same shortcut can be taken in parabolic rank $4$, but can fail beginning
with parabolic rank $5$.  The point is that if $P$ has parabolic rank $5$,
a vertex in $|\D_P|$ corresponds to a subgroup of parabolic rank $4$ and,
as we will indicate in the next paragraph, the isotypical components of
link intersection cohomology at a parabolic rank $4$ subgroup can live in
several degrees.

Incidentally, the figures show that in both parabolic rank $2$ and $3$, the
isotypical components \eqref{eqnLinkIHIsotypicalComponent} of link
intersection cohomology \eqref{eqnLinkICPRankTwo} occur in a single degree.
This phenomenon does not persist however; it can fail for parabolic rank
$\ge 4$.  For example, if $P$ has parabolic rank $4$ and thus $|\D_P|$ is a
$3$-dimensional simplex, there are two configurations that lead to
cohomology in more one degree: one can have truncation at all $4$ faces as
well as at $3$ non-coplanar edges or dually one can have truncation at all
$4$ vertices and at $3$ edges bounding a face.  The resulting complex in
the first case is $0\to\ZZ\to \ZZ^4 \to \ZZ^3 \to \ZZ$ which yields
cohomology in degrees $1$ and $2$.  These configurations can actually arise
from Weyl group elements.%
\footnote{Among the real symplectic groups one must wait until
$G=\Sympl_{20}(\RR)$ to see them.  For this group, the parabolic $P$ with
$\D^P = \{\al_2,\al_4,\al_5,\al_7,\al_8,\al_9\}$ has three elements $w\in
W_P$ which yield the first configuration above; one of them is
\begin{equation*}
\begin{split}
w =& s_{3} s_{2} s_{1} s_{6} s_{5} s_{4} s_{3} s_{2} s_{7} s_{6} s_{5}
s_{4} s_{3} s_{8} s_{10} s_{9} s_{8} s_{7} s_{6} s_{5} s_{4} \\
 & \qquad s_{10} s_{9} s_{8} s_{7} s_{6} s_{5} s_{10} s_{9} s_{8} s_{7}
s_{6} s_{10} s_{9} s_{8} s_{7} s_{10} s_{9} s_{8} s_{10} s_{9} s_{10}.
\end{split}
\end{equation*}
Unfortunately in this example, $\l(w)=42$, $\dim \n_P = 90$, and hence
$p_w(P) = 4$ or $4\frac12$.  This means that both of the link cohomology
classes in degrees $1$ and $2$ will not be truncated.  We expect that there
are examples in higher rank where only part of the link cohomology is
truncated.}

Finally, note that in all the illustrated examples, the intersection
cohomology groups $I_{p_w}H(|\D_P|;\ZZ)$ can be viewed simply as relative
cohomology groups $H(|\D_P|,C_w)$ for a certain subset $C_w \subseteq
\partial|\D_P|$.  Namely, $C_w$ is the union of the boundary strata at
which truncation is taking place.  Again this interpretation continues to
hold in parabolic rank $4$ but fails (at least with the above definition of
$C_w$) beginning with parabolic rank $5$.

\subsection{Parabolic rank \boldmath$\ge 3$}
\label{ssectPRankGEThree}
If $Q_V$ or $Q'_V$ is equal to $P$ or $G$, the argument already given using
Proposition ~\ref{propArbPrank} applies to $P$ with any parabolic rank to
show that $V=H^{\l(w)}(\n_P;E)_w$ can only contribute to the micro-support
when it is fundamental and can never contribute to the essential
micro-support \cite[Proposition ~17.7]{refnSaperLModules}.

It remains to consider $V$ such that $P<Q_V\le Q'_V < G$.  In this case
there is a direct though tedious proof for parabolic rank $3$ that $V$ is
not in the micro-support; it uses Figure
~\ref{figRelativeCohomologyPRankThree} and the same ideas as in the
parabolic rank $2$ case.  Obviously this approach does not lend itself to a
general argument; it requires a case-by-case calculation of $H(i_P^*
\ihat_Q^!  \IpC(\Xhat;E))_V$ for all the situations illustrated in Figure
~\ref{figRelativeCohomologyPRankThree}.

We will now outline the ingredients for a general proof.  In the case of
parabolic rank $3$, a complete proof can be constructed from what we
present here; in general one must consult \cite{refnSaperLModules}.  Recall
that we wish to show that
\begin{equation}
\label{eqnPurityVanishing}
H(i_P^* \ihat_Q^!\IpC(\Xhat;E))_V = 0
\end{equation}
under the hypotheses
\begin{enumerate}
\item\label{itemConjugateSelfContragredient} $V=H^{\l(w)}(\n_P;E)_w$ is
conjugate self-contragredient when restricted to $M_P$.
\item\label{itemQLocation} $P < Q_V\le Q \le Q'_V < G$.
\end{enumerate}
We can re-express
\begin{equation*}
H(i_P^* \ihat_Q^!\IpC(\Xhat;E))_V \cong H(\n_P;E)_w \otimes
I_{p_w}H_{c(|\D_P^Q|)}
\end{equation*}
\cite[Proposition ~17.4]{refnSaperLModules}, where we use the convenient
shorthand
\begin{equation*}
I_{p_w}H_{c(|\D_P^Q|)} \equiv I_{p_w}H\bigl(c(|\D_P|), c(|\D_P|)\setminus
c(|\D_P^Q|);\ZZ\bigr)
\end{equation*}
for the generalized perversity intersection cohomology supported on
$c(|\D_P^Q|)$.  (We will henceforth omit the coefficients $\ZZ$ from the
notation.)  Thus \eqref{eqnPurityVanishing} may be replaced by
\begin{equation}
\label{eqnCombinatorialPurityVanishing}
I_{p_w}H_{c(|\D_P^Q|)} = 0.
\end{equation}

We exploit the local calculation of intersection cohomology
\eqref{eqnLocalCalculation} at the vertex of the cone,
\begin{equation}
\label{eqnConePointVanishing}
I_{p_w}H^j(c(|\D_P|)) = \begin{cases}
  I_{p_w}H^j(|\D_P|) & \text{for $j\le p_w(P)$,}\\
  0                         & \text{for $j >  p_w(P)$.}
			    \end{cases}
\end{equation}
From this and the exact sequence of the pair, one can show that
\eqref{eqnCombinatorialPurityVanishing} is equivalent to
\begin{subequations}
\begin{align}
  &I_{p_w}H^{j-1}(|\D_P|\setminus |\D_P^Q|) = 0  
     &&\text{for $j > p_w(P) + 1$,} \label{eqnVanishingA} \\
  &\Im \left(I_{p_w}H^{p_w(P)}(|\D_P|\setminus |\D_P^Q|) \longrightarrow
     I_{p_w}H^{p_w(P)+1}_{|\D_P^Q|}\right) = 0
     &&\text{for $j = p_w(P) + 1$,} \label{eqnVanishingB} \\
  &I_{p_w}H^j_{|\D_P^Q|} = 0 &&\text{for $j <
     p_w(P)+1$.} \label{eqnVanishingC}
\end{align}
\end{subequations}
We will focus on proving \eqref{eqnVanishingA}.  It turns out that the
argument also shows that the domain of the map in \eqref{eqnVanishingB}
vanishes in half of the cases.  Then a dual argument proves
\eqref{eqnVanishingC} and also that the range of the map in
\eqref{eqnVanishingB} vanishes in the other half of the cases.

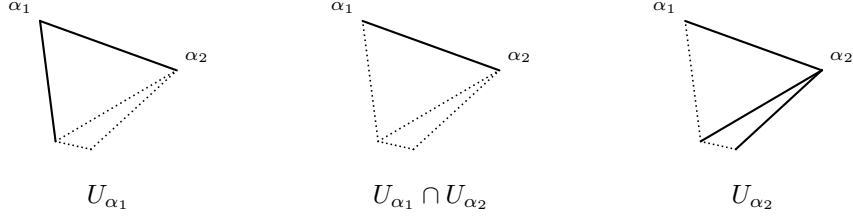
\begin{figure}
\begin{equation*}
\overset
{\begin{xy}
0;<2in,0in>:
(.1,.2)*i{\cdot},
(.0606,.5838)="a1"*+!RD{\scriptstyle \al_1},
(.4188,.4543)="a2"*+!LD{\scriptstyle \al_2},
(.1013,.2679)="a3",
(.1936,.2473)="a4",
"a1";"a2"**[|(2)]@{-},
"a1";"a3"**[|(2)]@{-},
"a2";"a3"**[|(1.5)]@{.},
"a2";"a4"**[|(1.5)]@{.},
"a3";"a4"**[|(1.5)]@{.}
\end{xy}}{U_{\al_1}}
\qquad\qquad
\overset
{\begin{xy}
0;<2in,0in>:
(.1,.2)*i{\cdot},
(.0606,.5838)="a1"*+!RD{\scriptstyle \al_1},
(.4188,.4543)="a2"*+!LD{\scriptstyle \al_2},
(.1013,.2679)="a3",
(.1936,.2473)="a4",
"a1";"a2"**[|(2)]@{-},
"a1";"a3"**[|(1.5)]@{.},
"a2";"a3"**[|(1.5)]@{.},
"a2";"a4"**[|(1.5)]@{.},
"a3";"a4"**[|(1.5)]@{.}
\end{xy}}{U_{\al_1}\cap U_{\al_2}}
\qquad\qquad
\overset
{\begin{xy}
0;<2in,0in>:
(.1,.2)*i{\cdot},
(.0606,.5838)="a1"*+!RD{\scriptstyle \al_1},
(.4188,.4543)="a2"*+!LD{\scriptstyle \al_2},
(.1013,.2679)="a3",
(.1936,.2473)="a4",
"a1";"a2"**[|(2)]@{-},
"a1";"a3"**[|(1.5)]@{.},
"a2";"a3"**[|(2)]@{-},
"a2";"a4"**[|(2)]@{-},
"a3";"a4"**[|(1.5)]@{.}
\end{xy}}{U_{\al_2}}
\end{equation*}
\caption{The Mayer-Vietoris spectral sequence for $I_{p_w}H(|\D_P|\setminus
  |\D_P^Q|)$, where $\D_P^Q=\{\al_3,\al_4\}$}
\label{figMayerVietoris}
\end{figure}

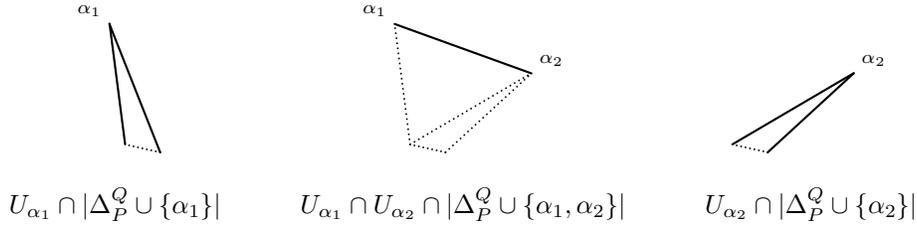
\begin{figure}
\begin{equation*}
\overset
{\begin{xy}
0;<2in,0in>:
(.1,.2)*i{\cdot},
(.0606,.5838)="a1"*+!RD{\scriptstyle \al_1},
(.4188,.4543)="a2",
(.1013,.2679)="a3",
(.1936,.2473)="a4",
"a1";"a4"**[|(2)]@{-},
"a1";"a3"**[|(2)]@{-},
"a3";"a4"**[|(1.5)]@{.}
\end{xy}}{U_{\al_1}\cap|\D_P^Q\cup\{\al_1\}|}
\qquad\quad
\overset
{\begin{xy}
0;<2in,0in>:
(.1,.2)*i{\cdot},
(.0606,.5838)="a1"*+!RD{\scriptstyle \al_1},
(.4188,.4543)="a2"*+!LD{\scriptstyle \al_2},
(.1013,.2679)="a3",
(.1936,.2473)="a4",
"a1";"a2"**[|(2)]@{-},
"a1";"a3"**[|(1.5)]@{.},
"a2";"a3"**[|(1.5)]@{.},
"a2";"a4"**[|(1.5)]@{.},
"a3";"a4"**[|(1.5)]@{.}
\end{xy}}{U_{\al_1}\cap U_{\al_2}\cap|\D_P^Q \cup \{\al_1,\al_2\}|}
\qquad\quad
\overset
{\begin{xy}
0;<2in,0in>:
(.1,.2)*i{\cdot},
(.0606,.5838)="a1",
(.4188,.4543)="a2"*+!LD{\scriptstyle \al_2},
(.1013,.2679)="a3",
(.1936,.2473)="a4",
"a2";"a3"**[|(2)]@{-},
"a2";"a4"**[|(2)]@{-},
"a3";"a4"**[|(1.5)]@{.}
\end{xy}}{U_{\al_2}\cap|\D_P^Q\cup\{\al_2\}|}
\end{equation*}
\caption{The Fary spectral sequence for $I_{p_w}H(|\D_P|\setminus
  |\D_P^Q|)$, where $\D_P^Q=\{\al_3,\al_4\}$}
\label{figFary}
\end{figure}

There are two spectral sequences \cite[Lemma ~3.7]{refnSaperLModules} we
use that abut to
\begin{equation}
\label{eqnArgumentAB}
I_{p_w}H^{j-1}(|\D_P|\setminus |\D_P^Q|).
\end{equation}
First some notation.  Let $S>P$ be the parabolic $\QQ$-subgroup
\emph{complementary} to $Q$ relative to $P$, that is, such that $\D_P^S =
\D_P \setminus \D_P^Q$.  For each $\al\in \D_P$, let $U_\al$ be the open
star neighborhood of the corresponding vertex of $|\D_P|$; note that $U_\al
= \bigcup_{\al\in \D_P^R\subseteq \D_P} |\D_P^R|^\circ$.  Set $U_R \equiv
\bigcap_{\al\in\D_P^R} U_\al $ for $R>P$, the open star neighborhood of
$|\D_P^R|^\circ$, and define a parabolic $\QQ$-subgroup $Q\vee R$ by
$\D_P^{Q\vee R} = \D_P^Q \cup \D_P^R$.  The \emph{Mayer-Vietoris spectral
sequence} for the open cover $\{U_\al\}_{\al\in \D_P^S}$ of
$|\D_P|\setminus |\D_P^Q|$ has%
\footnote{We are using $p$ both as an index and a perversity; this should
  not cause confusion.}
\begin{equation}
\label{eqnEOneMayerVietoris}
E_1^{p,j-1-p} = 
   \bigoplus_{\substack{P < R \le S \\ \#\D_P^R = p+1}}
   I_{p_w}H^{j-1-p}(U_R) \cong 
   \bigoplus_{\substack{P < R \le S \\ \#\D_P^R = p+1}}
   I_{p_w}H^{j-1-p}(c(|\D_R|)). 
\end{equation}
See Figure ~\ref{figMayerVietoris}.
On the other hand, the projection from $|\D_P^Q|$ defines a fibration
$|\D_P|\setminus |\D_P^Q| \to |\D_P^S|$.  For $P<R\le S$, the inverse image
of the stratum $|\D_P^R|^\circ \subseteq |\D_P^S|$ is $ U_R \cap
|\D_P^{Q\vee R}| \subseteq |\D_P|\setminus |\D_P^Q|$.  The associated
\emph{Fary spectral sequence} has
\begin{equation}
\label{eqnEOneFary}
E_1^{-p,j-1+p} =
   \bigoplus_{\substack{P < R \le S \\ \#\D_P^R = p}}
   I_{p_w}H^{j-1}_{U_R \cap |\D_P^{Q\vee R}|} \cong 
   \bigoplus_{\substack{P < R \le S \\ \#\D_P^R = p}}
   I_{p_w}H^{j-1}_{c(|\D_R^{Q\vee R}|)}. 
\end{equation}
See Figure ~\ref{figFary}.

The term of \eqref{eqnEOneMayerVietoris} associated to $R>P$ vanishes by
the truncation condition \eqref{eqnConePointVanishing} at the vertex when
$j - 1 - p > p_w(R)$, that is, when
\begin{equation}
\label{eqnWeakVanishingRange}
j > p_w(R) + \#\D_P^R.
\end{equation}
Since $p$ is a middle perversity, it is easy to verify that
\begin{equation*}
p_w(P) = p_w(R) + \lfloor (\dim\n_P^R + \#\D_P^R)/2\rfloor - \l^R(w) +
\delta,
\end{equation*}
where $\delta = 0$ or $1$ depending on the choice of $p$ and the parity of
$(\dim \n_P^R + \#\D_P^R)$ and $(\dim \n_P + \#\D_P)$.  However we know by
Proposition ~\ref{propArbPrank} (applied to $w^R$ and the parabolic $P/N_R$
of $L_R$) that
\begin{equation}
\label{eqnWeakInequality}
\l^R(w) \le (\dim\n_P^R)/2.
\end{equation}
Thus the Mayer-Vietoris spectral sequence yields the vanishing of
\eqref{eqnArgumentAB} for
\begin{equation}
\label{eqnWeakResult}
j >  p_w(P) + \lceil\#\D_P^R/2\rceil - \delta.
\end{equation}

Assume that $r(P)=3$.  Then $\#\D_P^R\le 2$ and thus \eqref{eqnArgumentAB}
vanishes for $j>p_w(P) + 1$ and even for $j=p_w(P)+1$ if $\delta=1$.  Thus
the Mayer-Vietoris spectral sequence argument yields a complete proof for
the parabolic rank $3$ case.

For $r(P)\ge 4$ the argument leading to \eqref{eqnWeakResult} is not
sufficient.  It would be sufficient, however, if \eqref{eqnWeakInequality}
were far enough away from equality.  It could also be improved if we had
more precise knowledge about $I_{p_w}H(|\D_R|)$, which suggests the utility
of an inductive argument.  On the other hand, the Fary spectral sequence
explicitly calls out for an inductive argument.  In fact, if we could apply
\eqref{eqnCombinatorialPurityVanishing} by induction, with $P$ replaced by
$R$, $w$ by $w_R$, and $Q$ by $Q\vee R$, all the terms in
\eqref{eqnEOneFary} would vanish.  Of course this is not possible (and
would lead to too strong a result) since the hypotheses
\itemref{itemConjugateSelfContragredient} and \itemref{itemQLocation} are
not preserved under such an induction.

The approach in \cite{refnSaperLModules} is to replace
\itemref{itemConjugateSelfContragredient} and \itemref{itemQLocation} by
other hypotheses that are more suited to induction and use a combination of
the two spectral sequences.  Recall that
\itemref{itemConjugateSelfContragredient} and \itemref{itemQLocation} were
used as input to Proposition ~\ref{propArbPrank}; more precisely they imply
that for all $A_P$-weights $\al$ occurring in $\n_P$,
\begin{equation}
\label{eqnOldInductionHypothesis}
\begin{cases}
\l_\al(w) \ge (\dim \n_\al)/2 &\text{if $\al$ occurs in $\n_P^Q$,}\\
\l_\al(w) \le (\dim \n_\al)/2 &\text{if $\al$ occurs in $\n_P^S$.}
\end{cases}
\end{equation}
Here $\dim \n_\al$ is the number of $\CC$-roots $\g$ such that $\g|_{\sa_P}
= \al$ and $\l_\al(w)$ is the number of such roots which in addition
satisfy $w^{-1}\g<0$.  The conditions \eqref{eqnOldInductionHypothesis} are
still not preserved by induction from $P$ to $R$.  However we can prove
\cite[Lemma ~17.8]{refnSaperLModules}%
\footnote{This is where the condition on the $\QQ$-root system in Theorem
~\ref{thmMicroSupportIC} comes into play.  Recall that a
\emph{quasi-minuscule representation} is one in which the non-zero weights
form a single Weyl group orbit.  The root systems of type $A_n$, $B_n$,
$C_n$, and $G_2$ are precisely the irreducible root systems such that the
corresponding simply connected $\QQ$-split group has a quasi-minuscule
representation whose weights under the usual ordering are totally ordered;
thus our hypothesis allows us to express the $\QQ$-roots of $G$ as
differences of weights $\varepsilon_i - \varepsilon_j$ analogously to the
usual expression of roots in $\GL_n$.  We expect that a further study of
quasi-minuscule representations will allow us to treat the cases $D_n$,
$E_n$, and $F_4$, although we may have to allow $r(T)=3$ or $4$.}
that there exists a parabolic $\QQ$-subgroup $T\ge P$ (depending on $w$)
with $r(T)=1$ or $2$ and $r(Q\vee T)$, $r(S\vee T)\le 1$ such that for all
$A_P$-weights $\al$ occurring in $\n_P$,
\begin{equation}
\label{eqnNewInductionHypothesis}
\begin{cases}
  \l_\al(w) \ge (\dim \n_\al)/2 &\text{if $\al$ occurs in $\n_T^{Q\vee T}$,}
  \\ 
  \l_\al(w) \le (\dim \n_\al)/2 &\text{if $\al$ occurs in $\n_T^{S\vee T}$.}
\end{cases}
\end{equation}
We illustrate the difference between conditions
\eqref{eqnOldInductionHypothesis} and \eqref{eqnNewInductionHypothesis} in
Figure ~\ref{figOldVersusNewInductionHypotheis} (here $P\le Q$ are as in
Figure ~\ref{figUnipotentRadicals}(b) and only one possible value of $T$ is
considered).
\begin{figure}
  \begin{center}
    \mbox{%
      \vbox{%
	\setbox\subfigbox=\hbox{%
	  $\left(\vcenter{
	    \xymatrix @=0pt@M=0pt@*=<.75cm,.75cm>{
	{} & {+} \save []!LD*+!UR{\scriptstyle \al_1} \restore
	  \ar@{-}[]!LD+0;[rrrr]!RD+0 \ar@{-}!LD+0;{[]!LU}+0 & {} & {} &
	  {} & {} \\
	{} & {} & {-} \save []!LD*+!UR{\scriptstyle \al_2} \restore
	  \ar@{-}[]!LD+0;[rrr]!RD+0 \ar@{-}[]!LD+0;[u]!LU+0 & {} &
	  {} & {} \\  
	{} & {} & {} & {+} \save []!LD*+!UR{\scriptstyle \al_3} \restore
	  \ar@{-}[]!LD+0;[rr]!RD+0 \ar@{-}!LD+0;{[uu]!LU}+0 & {} & {} \\
	{} & {} & {} & {} & {-} \save []!LD*+!UR{\scriptstyle \al_4} \restore
	  \ar@{-}[]!LD+0;[r]!RD+0 \ar@{-}[]!LD+0;[uuu]!LU+0 & {-} \\ 
        {} & {} & {} & {} & {} & {-} \save []!LD*+!UR{\scriptstyle \al_5}
	  \restore \ar@{-}[]!LD+0;[]!RD+0 \ar@{-}[]!LD+0;[uuuu]!LU+0 \\
        {}
	  }}
	  \right)$}%
	\dimen\subfigwd=\wd\subfigbox%
	\box\subfigbox%
	\vskip\baselineskip%
	\hbox{\vtop{\hbox to \dimen\subfigwd{\hfil\footnotesize(a) \quad $-
	      = \n_P^Q$\hfil}%
	    \hbox to \dimen\subfigwd{\hfil\footnotesize\hphantom{(a) \quad
	      }$+ = \n_P^S$\hfil}%
	}}%
      }%
      \qquad
      \vbox{%
	\setbox\subfigbox=\hbox{
	  $\left(\vcenter{
	    \xymatrix @=0pt@M=0pt@*=<.75cm,.75cm>{
	{} & {} \save []!LD*+!UR{\scriptstyle \al_1} \restore
	  \ar@{-}[]!LD+0;[rrrr]!RD+0 \ar@{-}!LD+0;{[]!LU}+0 & {-} & {} &
	  {} & {} \\
	{} & {} & {-} \save []!LD*+!UR{\scriptstyle \al_2} \restore
	  \ar@{-}[]!LD+0;[rrr]!RD+0 \ar@{-}[]!LD+0;[u]!LU+0 & {} &
	  {} & {} \\  
	{} & {} & {} & {+} \save []!LD*+!UR{\scriptstyle \al_3} \restore
	  \ar@{-}[]!LD+0;[rr]!RD+0 \ar@{-}!LD+0;{[uu]!LU}+0 & {+} & {+} \\
	{} & {} & {} & {} & {} \save []!LD*+!UR{\scriptstyle \al_4} \restore
	  \ar@{-}[]!LD+0;[r]!RD+0 \ar@{-}[]!LD+0;[uuu]!LU+0 & {} \\ 
        {} & {} & {} & {} & {} & {} \save []!LD*+!UR{\scriptstyle \al_5}
	  \restore \ar@{-}[]!LD+0;[]!RD+0 \ar@{-}[]!LD+0;[uuuu]!LU+0 \\
        {}
	  }}
	  \right)$}%
	\dimen\subfigwd=\wd\subfigbox%
	\box\subfigbox%
	\vskip\baselineskip%
	\hbox{\vtop{\hbox to \dimen\subfigwd{\hfil\footnotesize(b) \quad $-
	      = \n_T^{Q\vee T}$\hfil}%
	    \hbox to \dimen\subfigwd{\hfil\footnotesize\hphantom{(b) \quad
	      }$+ = \n_T^{S\vee T}$\hfil}%
	}}%
      }%
    }%
  \end{center}
  \caption{An example of $\n_P^Q$ versus $\n_T^{Q\vee T}$ and $\n_P^S$
    versus $\n_T^{S\vee T}$ in $\GL_n(\RR)$.  Here
    $\D_P=\{\al_1,\dotsc,\al_5\}$, $\D_P^Q= \{\al_2,\al_4,\al_5\}$, and
    $\D_P^T= \{\al_1,\al_4,\al_5\}$.}
  \label{figOldVersusNewInductionHypotheis}
\end{figure}
It is clear from the figure that significant information is gained in
passing from \eqref{eqnOldInductionHypothesis} to
\eqref{eqnNewInductionHypothesis} even though some information is lost as
well.  Clearly the conditions \eqref{eqnNewInductionHypothesis} are
preserved by induction from $P$ to $R$ provided $R\le T$.

The proof of  \eqref{eqnVanishingA}--\eqref{eqnVanishingC}, and hence
$I_{p_w}H_{c(|\D_P^Q|)} = 0$, in \cite[Theorem ~17.9]{refnSaperLModules}
uses  the existence of $T$ satisfying \eqref{eqnNewInductionHypothesis}
(plus one additional hypothesis).  The proof of \eqref{eqnVanishingA} is
separated into two parts.  The vanishing of $I_{p_w}H^{j-1}_{|\D_P^{Q\vee
T}|\setminus |\D_P^Q|}$, is proved using the Fary spectral sequence and
induction from $P$ to $R\le T$.  The other part, the vanishing of
$I_{p_w}H^{j-1}(|\D_P|\setminus|\D_P^{Q\vee T}|)$, is proved with the
Mayer-Vietoris spectral sequence (note that $\#\D_P^R\le 1$ in this case).

\section{Functoriality of micro-support}
\label{sectFunctorialityMicroSupport}
It is important to understand how the micro-support of an $\L$-module $\M$
changes when the various functors $k_*$, $k^*$, and $k^!$ are applied to
$\M$.  An example of this for the case of a direct image functor has
already been illustrated in \S\ref{ssectExampleMicroSupport} with the
calculation of $\mS(\i_{G*}E)$.  It remains to consider inverse image
functors.  Specifically one would like understand how $\mS(k^*\M)$, for
example, is derived from $\mS(\M)$.  Furthermore, in view of Theorem
~\ref{thmVanishingTheorem}, it is important also to estimate how $c(k^*\M)$
and $d(k^*\M)$ relate to the corresponding quantities for $\M$.  A number
of general theorems of this sort are found in \cite{refnSaperLModules};
here we will instead simply illustrate these ideas by sketching the proof
of the following result needed for Theorem ~\ref{thmRapoportConjecture} and
illustrated in Figure ~\ref{figProjection}.
\begin{figure}%
\begin{center}%
\begin{equation*}%
\begin{xy}
<0mm,-10mm>;<2mm,-10mm>:
<10mm,-2mm>="c",
(15,0)="adown";(0,10)="bdown" **\crv{(10,5)&(5,8)}
?(.5)="xdown"*-{\scriptstyle \bullet},
"adown";+"c" **\crv{(17.5,0)},
"bdown";{"bdown"+"c"+<.8mm,1.6mm>} **\crv{~*=<5pt>{.} "bdown"+<5mm,0mm>},
"adown"+"c";{"bdown"+"c"+<.8mm,1.6mm>} **\crv{~*\dir{} (10,5)+"c" & (5,8)+"c"}
?(.33)="e";"adown"+"c" **@{-},
"e";{"bdown"+"c"+<.8mm,1.6mm>} **\crv{~*=<5pt>{.} (10,5)+"c" & (5,8)+"c"},
"e"*++!LD{\scriptstyle X},
<0mm,10mm>;<2mm,10mm>:
<8mm,8mm>="d",
(15,0)="aup";(0,10)="bup" **\crv{(10,5)&(5,8)}
?(.5)="xup"*-{\scriptstyle \bullet},
"aup";+"d" **\crv{(17,0)},
"bup";{"bup"+"d"} **\crv{"bup"+<4mm,0mm>},
"aup"+"d";{"bup"+"d"} **\crv{"d"+(10,5) & "d"+(5,8)},
"adown";"aup" **@{-},
"bdown";"bup" **@{-},
"xdown";"xup" **@{-} ?(.5)*+!RD{\scriptstyle \Xhat_{R,\l}},
<0mm,-50mm>;<2mm,-50mm>:
<10mm,-2mm>="c",<8mm,8mm>="d",
(15,0)="a";(0,10)="b"
**\crv{(10,5)&(5,8)}
?(.5)="x"*-{\scriptstyle \bullet}*+!UR{\scriptstyle x}
?(.2)*!UR{\scriptstyle X_{R,h}},
"a";+"c" **\crv{(17.5,0)},
"a";+"d" **\crv{(17,0)},
"b";{"b"+"c"+<.4mm,1.6mm>} **\crv{~*=<5pt>{.} "b"+<5mm,0mm>},
"b";{"b"+"d"} **\crv{"b"+<4mm,0mm>},
"a"+"d";{"b"+"d"} **\crv{"d"+(10,5) & "d"+(5,8)},
"a"+"c";{"b"+"c"+<.8mm,1.6mm>} **\crv{~*\dir{} "c"+(10,5) & "c"+(5,8)}
?(.17)="e";"a"+"c" **@{-},
"e";{"b"+"c"+<.8mm,1.6mm>} **\crv{~*=<5pt>{.} "c"+(10,5) & "c"+(5,8)},
"e"*+++!L{\scriptstyle X},
"xdown"+<0mm,-8mm> \ar _{\pi} @{>} {"x"+<0mm,20mm>}
\end{xy}
\end{equation*}%
\end{center}%
\caption{The projection $\pi\colon \Xhat\to \Xstar$ over a point $x\in
X_{R,h}$ (adapted from \cite{refnSaperIHP})}%
\label{figProjection}%
\end{figure}%
The complete proof is found in \cite[Corollary 26.2]{refnSaperLModules}.

\begin{thm}
\label{thmFunctoriality}
Let $\M$ be an $\L$-module on $\Xhat$ with $\emS(\M)=\{E\}$ and
$c(E;\M)=d(E;\M)=0$.  Let $x\in X_{R,h}$ be a point on a stratum of a real
equal-rank Satake compactification $\Xstar$.  Let $\pi\colon \Xhat \to
\Xstar$ be Zucker's projection and let $k\colon \pi^{-1}(x)\cong
\Xhat_{R,\l}\hookrightarrow \Xhat$ denote the inclusion.  Then
\begin{equation*}
d(k^*\M) \le \tfrac12\codim X_{R,h} - \#\D_R
\qquad\text{and}\qquad
c(k^!\M) \ge \tfrac12\codim X_{R,h} + \#\D_R.
\end{equation*}
\end{thm}

We will make several simplifying assumptions along the way in order to
highlight the basic idea.  To begin with though, we do not make any
assumption about $\mS(\M)$ nor about the equal-rank nature of the boundary
components of $\Xstar$; these will be added later to emphasize where they
are needed.

\subsection{}
Consider the triple of subsets,
\begin{equation*}
\{x\} \subseteq X_{R,h} \subseteq U_{R,h} \subseteq \Xstar,
\end{equation*}
where
\begin{equation*}
U_{R,h} \equiv \coprod_{\cl{X_{S,h}}\supseteq X_{R,h}} X_{S,h}
\end{equation*}
is the open star neighborhood of $X_{R,h}$.  The inclusion $k$ may thus be
factored into the composition of three inclusions,
\begin{equation*}
\xymatrix @R-.5cm {
*+!{\pi^{-1}(x)} \xyhookrightarrow[r]^-{\ihat_{R,\ell}} & 
*+!{\pi^{-1}(X_{R,h})} \xyhookrightarrow[r]^-{\ihat_R} &
*+!{\pi^{-1}(U_{R,h})} \xyhookrightarrow[r]^-{j} & *+!{\Xhat} \\
{} & {\Xhat_R\cap W} \ar@2{-}[u] & {W\vphantom{\Xhat_R}} \ar@3{-}[u] },
\end{equation*}
and we will examine (in reverse order) how each affects micro-support.  For
simplicity we will only consider $k^*$; the functor $k^!$ is treated
similarly.

\subsection{An open embedding}
\enlargethispage*{\baselineskip}
The map
\begin{equation*}
\xymatrix {
  *+!{\displaystyle j\colon W =
    \coprod_{\substack{P\le R \\ P^\dag = R}}\,
    \coprod_{\substack{Q\ge P \\ Q\cap R = P^{\vphantom{\dag}}}} X_Q}
  \xyhookrightarrow[r] & *+!{\Xhat} }
\end{equation*}
is an embedding of an open admissible subspace.  The effect on
micro-support is simple \cite[Proposition ~22.2]{refnSaperLModules}:
\begin{equation}
\label{eqnMicroSupportOpenEmbedding}
\mS(j^* \M) = \mS(\M) \cap \Bigl(\coprod_{Q\in\Pl(W)}
\IrrRep(L_Q)\Bigr)
\end{equation}
and for $V\in \mS(j^*\M)$, $c(V;j^*\M)= c(V;\M)$ and $d(V;j^*\M)= (V;\M)$.

\subsection{The inclusion of the closure of a  stratum}
\begin{figure}
\entrymodifiers={}
\begin{align*}
&\xymatrix @!0 @M=0pt {{\txt{\llap{$G\colon$}}} &
{\circ} \ar@{-}[r] &
{\circ} \ar@{-}[r] &
{\circ} \ar@{-}[r] &
{\circ} \ar@{-}[r] &
{\circ} \ar@{-}[r] &
{\circ} \ar@{-}[r] &
{\circ} \save []+/u1ex/*!DC{\scriptstyle \al_1}\restore \ar@{-}[r] &
{\circ} \save [lllllll].[]!C *!/u1ex/\frm{_\}}
   *!/u3.25ex/{\scriptstyle \z(\D^R)}\restore \ar@{-}[r] &
{\circ} \save []+/u1ex/*!DC{\scriptstyle \al_0}\restore \ar@{-}[r] &
{\circ} \ar@{-}@<.15ex>[r] \ar@{-}@<-.15ex>[r] \ar@{}[r]|{<}&
{\circ} \save [l].[]!C *!/u1ex/\frm{_\}}
   *!/u3.25ex/{\scriptstyle \kap(\D^R)}\restore
   \save []="aaa",\POS "aaa"+/u1ex/+/r4ex/*!L{\scriptstyle
  \lsb\QQ\u} \ar @{.} "aaa"\restore} \\[1.5ex]
&\xymatrix @!0 @M=0pt {{\txt{\llap{$R\colon$}}} &
{\circ} \ar@{-}[r] &
{\circ} \ar@{-}[r] &
{\circ} \ar@{-}[r] &
{\circ} \ar@{-}[r] &
{\circ} \ar@{-}[r] &
{\circ} \ar@{-}[r] &
{\circ} \ar@{-}[r] &
{\circ} &
{} &
{\circ} \ar@{-}@<.15ex>[r] \ar@{-}@<-.15ex>[r] \ar@{}[r]|{<}&
{\circ} } \\[1.5ex]
&\xymatrix @!0 @M=0pt {{\txt{\llap{$P\colon$}}} &
{\circ}  &
{}  &
{\circ} \ar@{-}[r] &
{\circ} &
{} &
{\circ} &
{} &
{\circ} &
{} &
{\circ} \ar@{-}@<.15ex>[r] \ar@{-}@<-.15ex>[r] \ar@{}[r]|{<}&
{\circ}} \\[1.5ex]
&\xymatrix @!0 @M=0pt {{\txt{\llap{$Q\colon$}}} &
{\circ} \ar@{-}[r] &
{\circ} \ar@{-}[r] &
{\circ} \ar@{-}[r] &
{\circ} &
{} &
{\circ} &
{} &
{\circ} &
{} &
{\circ} \ar@{-}@<.15ex>[r] \ar@{-}@<-.15ex>[r] \ar@{}[r]|{<}&
{\circ}} \\[1.5ex]
&\xymatrix @!0 @M=0pt {{\txt{\llap{$T\colon$}}} &
{\circ} \ar@{-}[r] &
{\circ} \ar@{-}[r] &
{\circ} \ar@{-}[r] &
{\circ} &
{} &
{\circ} &
{} &
{\circ} \ar@{-}[r] &
{\circ} \ar@{-}[r] &
{\circ} \ar@{-}@<.15ex>[r] \ar@{-}@<-.15ex>[r] \ar@{}[r]|{<}&
{\circ}}
\end{align*}
\caption{An example of the Dynkin diagrams for $G$ and the Levi quotients
  of various parabolic $\QQ$-subgroups considered in the text.  (Here
  $G=\Sympl_{22}(\RR)$ and $\Xstar$ is the Baily-Borel-Satake
  compactification.)}
\label{figFunctorialityDynkinDiagrams}
\end{figure}
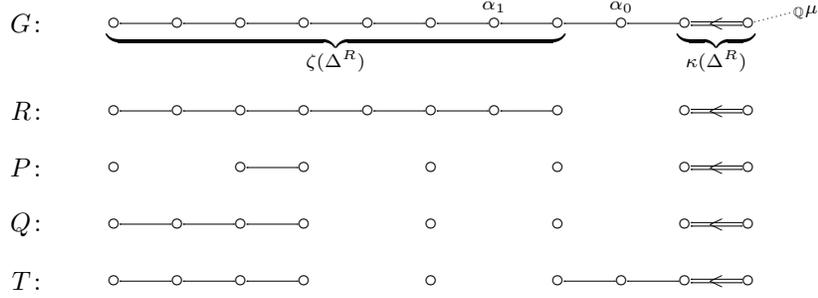

The next map
\begin{equation*}
\xymatrix {
  *+!{\displaystyle \ihat_R\colon \Xhat_R\cap W  =  \coprod_{\substack{P\le
	R \\ P^\dag = R}} X_P}
\xyhookrightarrow[r] & *+!{W} }
\end{equation*}
is the main issue; it is the embedding of a (relatively) closed stratum.
To treat it we outline a simplified version of the arguments of
\cite[Propositions ~22.6 and ~23.3]{refnSaperLModules}.  Set $\M_1 = j^*\M$
whose micro-support we understand by \eqref{eqnMicroSupportOpenEmbedding}.
Let $P\in\Pl(\Xhat_R\cap W)$.  If an $L_P$-module $V$ belongs to
$\mS(\ihat_R^*\M_1)$, then $V|_{M_P}$ is conjugate self-contragredient and
\begin{equation*}
H^j(\i_P^* \ihat_Q^! (\ihat_R^* \M_1))_V \neq 0
\end{equation*}
for some $Q\le R$ such that $Q_V^R\le Q \le Q_V^{\prime R}$ and some $j$.
Note that $\ihat_Q^! \ihat_R^* \M_1 = \ihat_T^! \M_1$ where $T\ge Q$ is the
unique parabolic $\QQ$-subgroup complementary to $R$ relative to $Q$, that
is, such that $\D_Q^T = \D_Q \setminus \D_Q^R$.  (See Figure
~\ref{figFunctorialityDynkinDiagrams} for an example of $P$, $Q$, $R$, and
$T$.)%
Thus we have
\begin{equation}
\label{eqnPTNonvanishing}
H^j(\i_P^* \ihat_T^! \M_1)_V \neq 0.
\end{equation}
We will show that if $V$ is not already in $\mS(\M_1)$, then it is
explicitly related to an element $\tilde V\in \mS(\ihat_{\tilde R}^*\M_1)$,
where $X_{\tilde R,h}$ is a larger stratum in $U_{R,h}$.  A repetition of
the argument will then lead us to an element in $\mS(\M_1)$.

For simplicity now assume:
\begin{enumerate}
\item\label{itemDynkinLinear} The Dynkin diagram of the $\QQ$-root system
for $G$ is linear.
\item\label{itemEndWeight} The highest $\QQ$-weight $\lsb\QQ\u$ of the
  representation defining the Satake compactification is orthogonal to all
  simple $\QQ$-roots except one at the end of the Dynkin diagram.
\item\label{itemEqualRank} All rational boundary components are equal-rank.
\end{enumerate}
These conditions are satisfied for practically all the real equal-rank
Satake compactifications (including the Baily-Borel-Satake
compactification) considered in Theorem ~\ref{thmRapoportConjecture}, and
in particular for the example illustrated in Figure
~\ref{figFunctorialityDynkinDiagrams}.  Assumptions
~\itemref{itemDynkinLinear} and \itemref{itemEndWeight} are simply for ease
of exposition; they ensure that
\begin{itemize}
\item All saturated parabolic $\QQ$-subgroups $R$ are maximal.
\end{itemize}
Assumption \itemref{itemEqualRank} will be needed to ensure that $\tilde V$
below is actually conjugate self-contragredient.

Since $R$ is maximal by the above assumptions, write $\D_P^R = \D_P
\setminus \{\al_0\}$ and let $\tilde P > P$ have $\D_P^{\tilde
P}=\{\al_0\}$.  We need to understand how $T$ relates to $Q_V$ and $Q_V'$.
By definition, $T=Q\vee \tilde P$, that is, $\D_P^T = \D_P^Q \cup
\{\al_0\}$.  For $Q_V$ and $Q_V'$ there are three cases:
\begin{itemize}
\item[(a)] $(\xi_V+\r,\al_0) < 0$ and thus $Q_V = Q_V^R\vee \tilde P$ and
  $Q_V' = Q_V^{\prime R}\vee \tilde P$.
\item[(b)] $(\xi_V+\r,\al_0) = 0$ and thus $Q_V = Q_V^R$ and $Q_V' =
  Q_V^{\prime R}\vee \tilde P$.
\item[(b)] $(\xi_V+\r,\al_0) > 0$ and thus $Q_V = Q_V^R$ and $Q_V' =
  Q_V^{\prime R}$.
\end{itemize}
In cases (a) and (b), $Q_V \le T \le Q_V'$ and thus equation
\eqref{eqnPTNonvanishing} implies that $V\in \mS(\M_1)$.  In case (c) on
the other hand, $T\not\le Q'_V$.

We investigate case (c) further.  Consider the long exact sequence of the
triple $(U, U\setminus (U\cap \Xhat_Q), U\setminus (U\cap \Xhat_T))$, where
$U$ is a small neighborhood of a point of $X_P$.  The $V$-isotypical part
of this sequence is
\begin{equation*}
\dotsb \longrightarrow H^j(\i_P^*\ihat_Q^!\M_1)_V
\longrightarrow H^j(i_P^* \ihat_T^! \M_1)_V 
\longrightarrow H^j(\i_P^*\i_{\tilde P*}(\i_{\tilde P}^*\ihat_T^!\M_1))_V
\longrightarrow \dotsb\ .
\end{equation*}
By \eqref{eqnPTNonvanishing}, the middle term is nonzero for some $j$, so
either the first or last term must also be nonzero.  If the first term is
nonzero, again $V\in \mS(\M_1)$ since $Q_V \le Q \le Q'_V$.  On the other
hand, it is easy to see that the last term can be rewritten as
$H^{\l}(\n_P^{\tilde P}; H^{j-\l}(\i_{\tilde P}^*\ihat_T^!\M_1))_V$, so if
it is nonzero we must have
\begin{equation}
\label{eqnVprimeNonvanishing}
H^{j-\l(w)}(\i_{\tilde P}^*\ihat_T^!\M_1)_{\tilde V}\neq 0
\end{equation}
for some irreducible $L_{\tilde P}$-module $\tilde V$ such that $V =
H^{\l(w)}(\n_P^{\tilde P};\tilde V)_w$ for some $w\in W_P^{\tilde P}$.  It
is not difficult to check \cite[Lemma ~23.2]{refnSaperLModules} that
$\tilde V|_{M_{\tilde P}}$ is conjugate self-contragredient from
\itemref{itemEqualRank} and the fact that any representation of an
equal-rank group is conjugate self-contragredient.

By assumptions \itemref{itemDynkinLinear} and \itemref{itemEndWeight}, the
roots in $\D_P$ have a linear ordering induced from that of $\D$ and
$\al_0$ will occur at one end; this means that $\al_0$ is orthogonal to all
$\al \in \D_P\setminus \D_P^{\tilde P}$ except for one, say $\al_1$.  Set
$\tilde \al = \al|_{A_{\tilde P}}$ for $\al\in \D_P\setminus \D_P^{\tilde
P}$.  One can show \cite[Lemmas ~22.5 and ~23.1]{refnSaperLModules} that
\begin{equation}
\label{eqnEqualityGreaterThan}
(\xi_{\tilde V}+\r,\tilde \al) 
\begin{cases}
= (\xi_V+\r,\al) & \text{for $\al \neq \al_1$,} \\
> (\xi_V+\r,\al_1) &  \text{for $\al = \al_1$}
\end{cases}
\end{equation}
(the inequality is since we are in case (c)).  Thus if $\al_1 \notin
\D_P^Q$ (as in Figure ~\ref{figFunctorialityDynkinDiagrams}), then
$Q_{\tilde V} \le T \le Q_{\tilde V}'$ and hence
\eqref{eqnVprimeNonvanishing} implies $\tilde V \in \mS(\M_1)$.

If $\al_1 \in \D_P^Q$, we want to apply to \eqref{eqnVprimeNonvanishing}
the reverse of the reasoning that led to \eqref{eqnPTNonvanishing}.  Set
$\tilde R$ to be the maximal $\QQ$-parabolic subgroup with $\D_P^{\tilde R}
= \D_P\setminus \{\al_1\}$ and set $\tilde Q\equiv \tilde R\cap T$ (see
Figure ~\ref{figFunctorialityDynkinDiagramsTilde}).%
\begin{figure}
\entrymodifiers={}
\begin{align*}
&\xymatrix @!0 @M=0pt {{\txt{\llap{$G\colon$}}} &
{\circ} \ar@{-}[r] &
{\circ} \ar@{-}[r] &
{\circ} \ar@{-}[r] &
{\circ} \ar@{-}[r] &
{\circ} \ar@{-}[r] &
{\circ} \ar@{-}[r] &
{\circ} \save []+/u1ex/*!DC{\scriptstyle \al_1}\restore \ar@{-}[r] &
{\circ} \save [lllllll].[]!C *!/u1ex/\frm{_\}}
   *!/u3.25ex/{\scriptstyle \z(\D^R)}\restore \ar@{-}[r] &
{\circ} \save []+/u1ex/*!DC{\scriptstyle \al_0}\restore \ar@{-}[r] &
{\circ} \ar@{-}@<.15ex>[r] \ar@{-}@<-.15ex>[r] \ar@{}[r]|{<}&
{\circ} \save [l].[]!C *!/u1ex/\frm{_\}}
   *!/u3.25ex/{\scriptstyle \kap(\D^R)}\restore
   \save []="aaa",\POS "aaa"+/u1ex/+/r4ex/*!L{\scriptstyle
  \lsb\QQ\u} \ar @{.} "aaa"\restore} \\[1.5ex]
&\xymatrix @!0 @M=0pt {{\txt{\llap{$P\colon$}}} &
{\circ}  &
{}  &
{\circ} \ar@{-}[r] &
{\circ} &
{} &
{\circ} &
{} &
{\circ} &
{} &
{\circ} \ar@{-}@<.15ex>[r] \ar@{-}@<-.15ex>[r] \ar@{}[r]|{<}&
{\circ}} \\[1.5ex]
&\xymatrix @!0 @M=0pt {{\txt{\llap{$Q\colon$}}} &
{\circ} \ar@{-}[r] &
{\circ} \ar@{-}[r] &
{\circ} \ar@{-}[r] &
{\circ} &
{} &
{\circ}  \ar@{-}[r] &
{\circ}  \ar@{-}[r] &
{\circ} &
{} &
{\circ} \ar@{-}@<.15ex>[r] \ar@{-}@<-.15ex>[r] \ar@{}[r]|{<}&
{\circ}} \\[1.5ex]
&\xymatrix @!0 @M=0pt {{\txt{\llap{$T\colon$}}} &
{\circ} \ar@{-}[r] &
{\circ} \ar@{-}[r] &
{\circ} \ar@{-}[r] &
{\circ} &
{} &
{\circ} \ar@{-}[r] &
{\circ} \ar@{-}[r] &
{\circ} \ar@{-}[r] &
{\circ} \ar@{-}[r] &
{\circ} \ar@{-}@<.15ex>[r] \ar@{-}@<-.15ex>[r] \ar@{}[r]|{<}&
{\circ}} \\[1.5ex]
&\xymatrix @!0 @M=0pt {{\txt{\llap{$\tilde P\colon$}}} &
{\circ}  &
{}  &
{\circ} \ar@{-}[r] &
{\circ} &
{} &
{\circ} &
{} &
{\circ} \ar@{-}[r] &
{\circ} \ar@{-}[r] &
{\circ} \ar@{-}@<.15ex>[r] \ar@{-}@<-.15ex>[r] \ar@{}[r]|{<}&
{\circ}} \\[1.5ex]
&\xymatrix @!0 @M=0pt {{\txt{\llap{$\tilde R\colon$}}} &
{\circ} \ar@{-}[r] &
{\circ} \ar@{-}[r] &
{\circ} \ar@{-}[r] &
{\circ} \ar@{-}[r] &
{\circ} \ar@{-}[r] &
{\circ}  &
{} &
{\circ} \ar@{-}[r] &
{\circ} \ar@{-}[r] &
{\circ} \ar@{-}@<.15ex>[r] \ar@{-}@<-.15ex>[r] \ar@{}[r]|{<}&
{\circ}} \\[1.5ex]
&\xymatrix @!0 @M=0pt {{\txt{\llap{$\tilde Q\colon$}}} &
{\circ} \ar@{-}[r] &
{\circ} \ar@{-}[r] &
{\circ} \ar@{-}[r] &
{\circ} &
{} &
{\circ} &
{}  &
{\circ} \ar@{-}[r] &
{\circ} \ar@{-}[r] &
{\circ} \ar@{-}@<.15ex>[r] \ar@{-}@<-.15ex>[r] \ar@{}[r]|{<}&
{\circ}}
\end{align*}
\caption{The case $\al_1 \in \D_P^Q$}
\label{figFunctorialityDynkinDiagramsTilde}
\end{figure}
Then $\ihat_T^! \M_1 = \ihat_{\tilde Q}^! \ihat_{\tilde R}^* \M_1$ and
equation \eqref{eqnEqualityGreaterThan} implies that
\begin{equation*}
Q_{\tilde V}^{\tilde R}\le \tilde Q \le Q_{\tilde V}^{\prime \tilde R}.
\end{equation*}
Thus we have 
\begin{equation*}
H^{j-\l(w)}(\i_{\tilde P}^* \ihat_{\tilde Q}^! (\ihat_{\tilde R}^*
\M_1))_{\tilde V} \neq 0
\end{equation*}
and hence $\tilde V\in \mS(\ihat_{\tilde R}^* \M_1)$.  In this case we may
repeat the argument with $V$ replaced by $\tilde V$ to eventually obtain an
element of $\mS(\M_1)$.

In conclusion, we see that for every $V\in \mS(\ihat_R^*\M_1)$, either
\begin{enumerate}
\item $V\in \mS(\M_1)$, or
\item\label{itemFunctorialInductionStep} $V = H^{\l(w)}(\n_P^{\tilde
P};\tilde V)_w$ for an $L_{\tilde P}$-module $\tilde V\in \mS(\M_1)$
satisfying
\begin{gather*}
\tilde P\cap R = P, \\
(\xi_V+\r,\al) > 0\text{ for all $\al\in \D_P^{\tilde P}$, and } \\
d(V;\ihat_R^*\M_1) \le   d(\tilde V;\M_1) + \l(w),
\end{gather*}
or
\item $V$ is obtained from an element of $\mS(\M_1)$ by a sequence of steps
  similar to \itemref{itemFunctorialInductionStep}.
\end{enumerate}

\subsection{The inclusion of a fiber}
The set $\Xhat_R\cap W = \pi^{-1}(X_{R,h})$ is a flat bundle over $X_{R,h}$
with typical fiber $\Xhat_{R,\ell}$.  Our final map
\begin{equation*}
\xymatrix { *+!{\ihat_{R,\l}\colon \pi^{-1}(x)\cong \Xhat_{R,\ell}}
  \xyhookrightarrow[r] & *+!{\Xhat_R\cap W} }
\end{equation*}
is the inclusion of such a fiber.  Thus $\Xhat_{R,\ell}$ is not an
admissible subspace of $\Xhat_R\cap W$ and the functors $\ihat_{R,\l}^*$
and $\ihat_{R,\l}^!$ are defined differently from those in
\S\ref{ssectFunctorsLModules}.

Instead, note that the spaces $\Xhat_{R,\ell}$ and $\Xhat_R\cap W$ have
strata indexed by the same set: $P\in \Pl$ such that $P\le R$ with $P^\dag
= R$.  For such $P$, the corresponding stratum of $\Xhat_{R,\ell}$ is
$X_{P,\ell}$ while the corresponding stratum of $\Xhat_R\cap W$ is $X_P$.
Let $\M_2 \equiv \ihat_R^*\, j^*\M = (E_\cdot,f_{\cdot\cdot})$ be an
$\L$-module on $\Xhat_R\cap W$ (whose micro-support we understand by the
preceding two subsections).  The $\L$-module $\ihat_{R,\l}^*\M_2 =
(E'_\cdot,f'_{\cdot\cdot})$ on $\Xhat_{R,\l}$ is defined
\cite[\S3.5]{refnSaperLModules} by
\begin{align*}
E'_{P_\l} & \equiv \Res_{L_{P,\l}}^{L_P} E_P, \\
f'_{P_\l Q_\l} & \equiv \Res_{L_{P,\l}}^{L_P} f_{PQ};
\end{align*}
the $\L$-module $\ihat_{R,\l}^!\M_2$ is defined similarly but with a degree
shift of $-\dim D_{R,h}$.

The effect of $\ihat_{R,\l}^*$ (and likewise $\ihat_{R,\l}^!$) on
micro-support is simple \cite[Proposition ~2.8]{refnSaperLModules}:
\begin{equation*}
\mS(\ihat_{R,\l}^*\M_2) = \{\, \Res_{L_{P,\l}}^{L_P} V \mid V\in
\mS(\M_2)\,\}.
\end{equation*}
(This is an abuse of notation since $\Res_{L_{P,\l}}^{L_P} V$ is not
actually an irreducible $L_{P,\l}$-module, merely isotypical.  In the above
formula we mean the unique irreducible $L_{P,\l}$-module that occurs in
$\Res_{L_{P,\l}}^{L_P} V$.)  

\subsection{Sketch of the proof of Theorem ~\ref{thmFunctoriality}}
Since $k^*\M = \ihat_{R,\ell}^*\, \ihat_R^*\, j^* \M$ and we assume that
$\emS(\M)=\{E\}$ with $c(E;\M)=d(E;\M)=0$, Theorem ~\ref{thmFunctoriality}
may be deduced from the $3$ preceding subsections.  The only issue is to
estimate $\l(w)$ which is accomplished by the following proposition.  This
is the point at which we need to assume that all real (as opposed to
rational) boundary components are equal-rank.

\begin{prop}
\label{propEqualRankBasicLemma}
Let $\lsb\RR\Dstar$ be a real equal-rank Satake compactification.  Let $P$
be a parabolic $\QQ$-subgroup and let $w\in W_P$.  Let $V =
H^{\l(w)}(\n_P;E)_w$ be the corresponding irreducible $L_P$-module and
assume that $(V|_{M_P})^* \cong \overline{V|_{M_P}}$.
\begin{enumerate}
\item\label{itemEqualRankNegative} If $(\xi_V+\r,\al) \le0$ for all $\al\in
\D_P$, then
\begin{equation*}
\l(w) \ge  (\dim\n_P + \#\D_P + \dim D_{P,\l}(V))/2.
\end{equation*}
\item\label{itemEqualRankPositive} If $(\xi_V+\r,\al) \ge0$ for all $\al\in
\D_P$, then
\begin{equation*}
\l(w) \le (\dim\n_P - \#\D_P - \dim D_{P,\l}(V))/2.
\end{equation*}
\end{enumerate}
\end{prop}

The proposition follows from Lemma ~\ref{lemBasicLemma} and the estimate
\cite[Lemma~25.2]{refnSaperLModules}
\begin{equation}
\dim \n_P(V) \ge \#\D_P + \dim D_{P,\l}(V).
\end{equation}
A similar estimate was proved in \cite[\S5.5]{refnBorelVanishingTheorem}
which only requires that $D$ and $D_{R,h}$ are equal-rank and that a
certain condition (B) holds which has to be verified case-by-case (it is
satisfied for most real equal-rank Satake compactifications.)  The proof of
the estimate in \cite{refnSaperLModules} is free from case-by-case analysis
but strongly uses the assumption that all real boundary components are
equal-rank.

\section{Proof of the Rapoport/Goresky-MacPherson conjecture}
To prove Theorem ~\ref{thmRapoportConjecture} we need to show that
$\pi_*\IpC(\Xhat;\EE) \cong \IpC(\Xstar;\EE)$.  This is done by verifying
that $\pi_*\IpC(\Xhat;\EE)$ satisfies the local characterization of
intersection cohomology on $\Xstar$, conditions
\itemref{itemCoefficientsOnSmooth}--\itemref{itemLocalAttaching} from
\S\ref{ssectSheafTheoreticIC}.  Condition
~\itemref{itemCoefficientsOnSmooth} is obvious.  Let $i_x\colon
\{x\}\hookrightarrow \Xstar$ denote the inclusion of a point $x\in \Xstar$.
The local vanishing condition ~\itemref{itemLocalVanishing} amounts to
\begin{equation}
\label{eqnIHVanishingCondition}
H^j(i_x^*\pi_*\IpC(\Xhat;\EE)) =0 \qquad \text{for $x\in X_{R,h}$,
  $j\ge\tfrac12\codim X_{R,h}$}
\end{equation}
since $p(k) = k/2 - 1$ for either middle perversity when $k$ is even.
Given the other conditions and constructibility, the attaching condition
~\itemref{itemLocalAttaching} is equivalent
\cite[V.4]{refnBorelIntersectionCohomology},
\cite{refnGoreskyMacPhersonIHTwo} to a local covanishing condition:
\begin{equation}
\label{eqnIHCoVanishingCondition}
H^j(i_x^!\pi_*\IpC(\Xhat;\EE)) =0 \qquad \text{for $x\in X_{R,h}$,
  $j\le\tfrac12\codim X_{R,h}$}.
\end{equation}

Let $k\colon \pi^{-1}(x) \hookrightarrow \Xhat$ be the inclusion.  We may
re-express
\begin{equation*}
H^j(i_x^*\pi_*\IpC(\Xhat;\EE)) \cong H^j(\pi^{-1}(x);k^*\IpC(\Xhat;\EE));
\end{equation*}
since $\IpC(\Xhat;\EE) = \Sheaf(\IpC(\Xhat;E))$, this is isomorphic to the
$\L$-module cohomology $H^j(\Xhat_{R,\l}; k^*\IpC(\Xhat;E))$.  Thus we can
use Theorem ~\ref{thmVanishingTheorem} to see this vanishes for $j >
d(k^*\IpC(\Xhat;E))$.  By Theorem ~\ref{thmMicroSupportIC}, we know that
$\emS(\IpC(\Xhat;E))=\{E\}$ and hence $d(k^*\IpC(\Xhat;E)) \le
\tfrac12\codim X_{R,h} -1$ by Theorem ~\ref{thmFunctoriality}.  This
establishes \eqref{eqnIHVanishingCondition}; an analogous argument treats
\eqref{eqnIHCoVanishingCondition}.

\section{A generalization of Goresky-Harder-MacPherson's theorem}
The weighted cohomology $\L$-module $\WnC(\Xhat;E)$ may be constructed
\cite[\S6]{refnSaperLModules} by a weight truncation functor similar to the
degree truncation functor that is used for $\IpC(\Xhat;E)$.  The following
analogue of Theorem ~\ref{thmMicroSupportIC} is proved in
\cite[Theorem~16.3]{refnSaperLModules}.  Note that there is no assumption
on the $\QQ$-root system.  Unlike Theorem ~\ref{thmMicroSupportIC}, this
theorem is not difficult to prove; this is because there is an explicit
non-inductive formula for the local weighted cohomology.
\begin{thm}
\label{thmMicroSupportWC}
Let $E$ be an irreducible $G$-module and let $\eta=\u$ or $\v$ be a
middle-weight profile.  If $(E|_{\lsp0\,G})^* \cong
\overline{E|_{\lsp0\,G}}$, then $\emS(\WnC(\Xhat;E)) = \{E\}$.
\end{thm}

Consequently, the argument in the preceding section also proves the
following generalization of Goresky, Harder, and MacPherson's result
\cite{refnGoreskyHarderMacPherson}, Theorem
~\ref{thmGoreskyHarderMacPherson}:
\begin{thm}
Let $\Xstar$ be a real equal-rank Satake compactification, let $p$ be a
middle-perversity, and let $\eta$ be a middle-weight profile.  There is a
natural quasi-isomorphism $\IpC(\Xstar;\EE) \cong \pi_*\WnC(\Xhat;\EE)$ and
hence an isomorphism $I_pH(\Xstar;\EE)\cong W^\eta H(\Xhat;\EE)$.
\end{thm}

\bibliographystyle{amsplain}
\bibliography{cdm02_paper}
\end{document}